
\magnification 1200
\input amstex
\documentstyle{amsppt}
\NoBlackBoxes
\NoRunningHeads

\hsize = 6.2 truein
\vsize = 9.1 truein

\define\vre{\varepsilon}
\define\hd{Hausdorff dimension}
\define\hs{homogeneous space}
\define\df{\overset\text{def}\to=}
\define\br{\Bbb R}
\define\bn{\Bbb N}
\define\bz{\Bbb Z}
\define\bq{\Bbb Q}

\define\vwa{very well approximable}
\define\vwma{very well multiplicatively approximable}
\define\di{Diophantine}
\define\da{Diophantine approximation}

\define\va{\bold a}
\define\vb{\bold b}
\define\ve{\bold e}
\define\vx{\bold x}
\define\vy{\bold y}

\define\vs{\bold s}
\define\vw{\bold w}

\define\vq{\bold q}

\define\vc{\bold c}
\define\vf{\bold f}
\define\vt{\bold t}

\define\vv{\bold v}
\define\vh{\bold h}

\define\nz{\smallsetminus \{0\}}

\define\be{Besicovitch}

\define\cag{$(C,\alpha)$-good}

\define\p{\Phi}
\define\vrn{\varnothing}
\define\ssm{\smallsetminus}
\define\dist{\operatorname{dist}}
\define\diam{\operatorname{diam}}
\define\supp{\operatorname{supp}}
\define\cov{\operatorname{cov}}
\define\Lie{\operatorname{Lie}}
\define\rk{\operatorname{rk}}
\define\rank{\operatorname{rank}}
\define\GL{\operatorname{GL}}
\define\Gr{\operatorname{Gr}}
\define\Int{\operatorname{Int}}
\define\SL{\operatorname{SL}}

\define\Ad{\operatorname{Ad}}
\define\const{\operatorname{const}}
\define\un#1#2{\underset\text{#1}\to#2}

\newif\ifdraft\drafttrue

\define\spr{Sprind\v zuk}



\topmatter
\title Flows on $S$-arithmetic homogeneous spaces and applications  to
metric Diophantine approximation
\endtitle

\author { Dmitry Kleinbock and George Tomanov}
\endauthor

            \address{Dmitry Kleinbock,  Department of
Mathematics, Brandeis University, Waltham, MA 02454-9110, USA}
          \endaddress

\email kleinboc\@brandeis.edu \endemail

            \address{ George Tomanov, Institut Girard Desargues,
Universit\` e Claude Bernard -- Lyon 1,
69622 Villeurbanne cedex,
France}
          \endaddress

\email Georges.Tomanov\@desargues.univ-lyon1.fr \endemail

\abstract The main goal of this work is to
establish quantitative nondivergence estimates for flows on
homogeneous spaces of products of real and $p$-adic Lie groups. 
These results have applications both to ergodic theory 
and to \da. Namely, earlier results of Dani (finiteness of 
locally finite ergodic
unipotent-invariant measures on real homogeneous
spaces) and Kleinbock-Margulis
(strong extremality of nondegenerate submanifolds of $\br^n$) 
are generalized to the 
$S$-arithmetic setting.

\comment . We also provide  proofs of $S$-adic
analogs of earlier results of S.\,G.\ Dani about locally finite
invariant under unipotent flows measures on real homogeneous
spaces.
establish $S$-arithmetic generalization of
\endabstract


\endtopmatter

\document

\heading{0. Introduction}\endheading

\subhead{0.1}\endsubhead Actions of unipotent subgroups on \hs s
of real Lie groups provide examples of important dynamical systems  with numerous
applications to geometry and number theory. A fundamental
phenomemon discovered by G.\,A.\ Margulis in 1971 \cite{M} showed
that, in sharp contrast to partially hyperbolic flows, orbits of
unipotent flows are never divergent.
Then in the papers \cite{D1--3}
   S.\,G.\ Dani generalized and
strengthened this result. In particular, for one-parameter flows
he proved that any unipotent orbit returns to big compact subsets with high
frequency. A very general explicit estimate for this frequency in
terms of the size of compact sets was given in 1998 in the paper
of Margulis and the first named author \cite{KM}. In fact it was
done in a bigger generality, i.e.\ for
a large class of maps
from $\br^d$ into $\SL(n,\br)$, making it possible to derive important
applications  to metric \da\ on manifolds.
The nondivergence theorem of  Dani and Margulis was used
in M.\ Ratner's  proof  \cite{Rt1--2} of Raghunathan's
topological conjecture, and
in Dani's proof \cite{D2} of  finiteness of  locally  finite  ergodic measures 
invariant under a
unipotent flow.

\comment

The nondivergence theorem of  Dani and Margulis was used
by M.\ Ratner in her proof  \cite{Rt1--2} of Raghunathan's
topological conjecture, and by  
Dani to  prove  that any locally finite ergodic measure invariant under a
unipotent flow is finite \cite{D2, Theorem 4.3}.

Earlier work of Margulis and Dani had played a crucial role in
M.\ Ratner's \cite{Rt1--2}  proof  of  conjectures of
Raghunathan, Dani and Margulis on orbit closures and invariant
measures of groups generated by one-parameter unipotent subgroups.
In particular, 
   Dani used his results to prove that any locally finite
measure invariant under a unipotent flow is finite \cite{D2,
Theorem 4.3}.

\endcomment

In \cite{BP}, A.\ Borel and G.\ Prasad initiated the study of
actions of $S$-adic unipotent groups on \hs s.
The generalization of Ratner's results to
the
$S$-adic setting was achieved in the papers \cite{MT1--2} and  \cite{Rt3--4}.
In the present paper we prove an $S$-adic analogue of explicit
quantitative estimates from \cite{KM}.  As one of the
applications, we  show the finiteness of  locally finite
unipotent-invariant ergodic measures on $S$-adic \hs s. We also
derive applications to number theory. Namely, generalizing the
correspondence between homogeneous dynamics and \da\ over $\br$,
we prove $p$-adic, and more generally,  $S$-arithmetic versions of
conjectures of A.\ Baker and V.\  \spr, in particular answering a
question posed by I.\ Shparlinski.

In order to state our results, let us introduce some notation. Let
$S$ be a finite set of normalized valuations of $\bq$ containing
the Archimedean one, $\bz_S$ the ring of $S$-integers of $\bq$,
and $\bq_S$ the direct product of completions $\bq_v$ of $\bq$
over $v\in S$. Put $\GL(m,\bq_S) \df \prod_{v\in S}\GL(m,\bq_v)$
and
$$ \GL^1(m,\bq_S) \df \left\{(g^{(v)})_{v\in S}\in \GL(m,\bq_S)\left|
\prod_{v\in S}|\det(g^{(v)})|_v = 1\right.\right\}\,.
$$
Then one can  interpret the homogeneous space $$\Omega_{S,m}^1 =
\GL^1(m,\bq_S)/\GL(m,\bz_S)
$$
       as the space of  lattices $\Lambda$ in $\bq_S^m$ of covolume
       $1$.
For any $v\in S$, the valuation $|\cdot|_v$ induces the norm
$\|\cdot\|_v$ on $\bq_v^m$, and we let 
$$c(x) \df \prod_{v \in
S}|x^{(v)}|_v\quad\text{for}\quad x = (x^{(v)})_{v\in S}\in \bq_S\,,\tag 0.1
$$ and $c(\vx)
\df \prod_{v \in S}\|\vx^{(v)}\|_v$ for $\vx = (\vx^{(v)})_{v\in
S}\in \bq_S$. The latter function plays the same role as the usual
norm in the case $S = \{\infty\}$, see Lemma 7.8.

Now define
$$
Q_\vre \df \left\{\Lambda \in \Omega_{S,m}^1  \bigm| \Lambda
\text{ contains no  }\vx \text{ with  }c(\vx) <
\vre\right\}\,.\tag 0.2
$$
It follows from the generalized Mahler's Compactness Criterion
(Theorem 7.10) that the sets $Q_\vre$ are compact.


       Let us now formulate a special case of our Theorem 8.4.

           \proclaim{Theorem}  Consider the space
  $X = \prod_{v\in S}\bq_v^{d_v}$,
$d_v \in\bn$, endowed with the  product metric and the measure
$$\lambda \df \prod_{v\in S}\lambda_v\,,\text{ where $\lambda_v$
is a Haar measure on $\bq_v^{d_v}$}\,,\tag 0.3
$$
and let $h = (h_v)_{v\in S}: X\to \GL^1(m,\bq_S)$
be a polynomial map (see \S 8.4 for the definition) with $h(0) =
e$.
Then there exists $\alpha > 0$
such that for every $\Lambda\in
\Omega_{S,m}^1 $ one can find
$C > 0$
   with the
following property: for any positive $\vre$ and any ball $B\subset
X $ centered at $0$
        one has
$$
\frac{\lambda\left(\big\{x\in B\bigm| h(x)\Lambda\notin  Q_\vre
\big\}\right)}{\lambda(B)}
\le C\vre^\alpha
\,.
$$
\endproclaim


We remark that Theorem 8.4 is in its turn a special case of a more
general result, Theorem 8.3, where both the polynomiality of the
map $h$ and the restriction that $\Lambda$ belongs to a compact
set $L$ of $\Omega_{S,m}^1$ are significantly relaxed.

\subhead{0.2}\endsubhead Theorem 8.4 makes it possible to
generalize the main results from \cite{D1--3} to the case of
products of real and $p$-adic algebraic groups. In particular,
here is a special case of our Theorem 9.1:

\proclaim{Theorem} Let $G = \prod_{v \in S} {G}_v$, where $G_v$ is
the group of $\bq_v$-rational points of an algebraic group
$\bold{G}_v$ defined over $\bq_v$, let  $\Gamma$ be a lattice in
$G$, $d$  a natural number, and let $X$ and $\lambda$ be as in
Theorem 0.1. Then for any  $\beta
>       0$ there exists a compact $M_\beta
\subset G/\Gamma$  such that
for any $y \in G/\Gamma$ and any  parametrization $\phi: X\to U$ of
degree $\leq d$ of a unipotent algebraic subgroup $U$ of $G$ (see
{\rm \S 9.1} for precise definitions)
one
of the following holds: either

{\rm(i)} there exists $R = R(y,\phi)$ such that
$$
\frac{\lambda \left(\big\{x\in B\bigm| \phi(x)y \in M_\beta
\big\}\right)}{\lambda(B)} \geq 1 -\beta 
$$
for any ball
$B$  in
$X$
centered at $0$
with radius at least $R$, or

{\rm(ii)}  there exists a closed
proper subgroup $H$ of $G$ containing $U$ such that the orbit
$Hy$ is closed and  carries a finite $H$-invariant Borel measure.
\endproclaim

In the particular case when $U$ is one-parameter, Theorem   0.2
has been announced, with indication of the proof, in \cite{MT2,
Theorem 11.4}, \cite{Rt4, Theorem 9.1} and \cite{To, Theorem 3.3},
and used for the proof of Ratner's uniform distribution theorem
\cite{Rt4, Theorem 3}  and for other related results.

\subhead{0.3}\endsubhead
The next theorem, 
which follows from Theorem 9.1,
has been proved for real Lie groups  by Dani \cite{D2, Theorem 4.3}, and for the general
case has been announced  in \cite{MT2, Theorem 11.5}.

\proclaim{Theorem} Let  $G$ and $\Gamma$ be as in Theorem 0.2, let
$H$ be a subgroup of $G$ generated by unipotent one-parameter
subgroups of $G$, and let $\mu$ be a locally finite $H$-invariant
measure on $G/\Gamma$. Then there exist Borel $H$-invariant
subsets $Y_i$, $1 \leq i < \infty$, such that $\mu(Y_i) < \infty$
for all $i$ and $\mu(G/\Gamma \ssm \bigcup_{i}Y_i) = 0$. In
particular, every locally finite $H$-invariant ergodic measure on
$G/\Gamma$ is finite.
\endproclaim

\subhead{0.4}\endsubhead Now let us turn to applications to number
theory. Let $S$ be a finite set consisting of distinct normalized
valuations of $\bq$, with or without the infinite valuation.
We will  interpret
elements
$$
\vy\ =\ (\vy^{(v)})_{v\in S}\ =\ (y_1,\dots,y_n)
$$
of $\bq_S^n = \prod_{v\in S}\bq_{v}^{n}$, where $\vy^{(v)}=
(y^{(v)}_1,\dots,y^{(v)}_n)\in \bq_{v}^n$ and $y_i =
(y^{(v)}_i)\in \bq_S$, as linear forms on $\bq_S^n$, and will
study their values  $\vy\cdot\vq = y_1q_1 +\dots+y_nq_n$ at
integer points $\vq = (q_1,\dots,q_n)$. Denote by $\ell$ the
cardinality of $S$, and say that $\vy\in \bq_S^n$ is {\sl \vwa\/},
or {\sl VWA\/},  if
for some $\vre> 0$ there are infinitely many solutions $\tilde\vq
= (q_0,q_1,\dots,q_n) \in\bz^{n+1}$ to
$$|q_0 + \vq\cdot\vy|^\ell \le \cases \|\tilde\vq\|^{-(n+1)
(1 + \vre)}\quad\text{ if }\infty\notin S\\
\|\vq\|^{-n(1 + \vre)}\qquad\text{\ \  \ if }\infty\in S\,.\endcases\
\tag
0.4
$$
Also, say that $\vy$ is  {\sl \vwma\/},  or {\sl VWMA\/}, if
   for some $\vre > 0$ there are infinitely many $\tilde \vq\in \bz^{n+1}$
       such that
$$
c(q_0 + \vq\cdot\vy) \le \cases \Pi_{\sssize +}(\tilde \vq)^{-
(1 + \vre)}\qquad\quad\text{\ \  \ if }\infty\notin S\\
\Pi_{\sssize +}(\vq)^{-(1 + \vre)}|q_0|_+^{- \vre}\quad\text{ if
}\infty\in S\,,\endcases \tag 0.5
$$
where we put $$|x|_+ \df \max(|x|, 1)\quad\text{for }x\in\br\,,\quad
\Pi_{\sssize +}(\vx) \df \prod_{i = 1}^n
|x_i|_+\quad\text{for }\vx\in\br^k\,,$$
and define $c(\cdot)$ as in (0.1).
This unifies the standard definitions in the real and $p$-adic
set-ups; we refer the reader to \S 10 for motivation, in particular
for an explanation of the term  $|q_0|^{- \vre}$ in the second
line of (0.5), and for a uniform way to write down the expressions
in the right hand sides of  (0.4) and (0.5).
It is easy to check using the Borel-Cantelli lemma that the set
of VWMA vectors (and hence  the smaller set of VWA vectors) has
zero Haar measure.

The subject of metric \da\ with dependent quantities originated
with a conjecture of Mahler 
proved by \spr\ in the 1960s
(see \cite{Sp1}), stating  that for almost every $x \in \br$, the
vector
$$
\vf(x) = (x,x^2,\dots,x^n)
\tag 0.6
$$
is not VWA.
   On the other hand, a similar statement with VWA replaced by VWMA,
   conjectured by A.~Baker in 1972
\cite{B}, has not been proved until the paper \cite{KM}, which
introduced  a dynamical approach to this class of problems. (See
also \cite{K1} for a  survey.)
Note that \spr\ had also proved
   the $p$-adic counterpart of Mahler's  conjecture, and the problem of
establishing its multiplicative version, that is, proving that   the
vector (0.6) is not VWMA for $\lambda$-a.e.\ $x\in\bq_p$, was
recently posed
by I.\ Shparlinski (V.~Bernik, private communication).

More generally, following \cite{KLW}, let us say
that a measure $\mu$ on $\bq_S^n$ is {\sl extremal\/}
(resp.~{\sl strongly extremal\/}) if $\mu$-almost every vector in
$\bq_S^{n}$ is not VWA (resp., not VWMA).
In these terms, the conjectures of Mahler (resp.~Baker) state that
the pushforward
$\vf_*\lambda$ of Lebesgue measure $\lambda$ on $\br$ by the  map
$\vf$ as in (0.6) is extremal (resp.~strongly extremal).


An important property of the curve (0.6) is that it does not
belong to any proper affine subspace of $\br^n$. More generally,
consider a $C^k$ map $\vf:U\to F^n$, where $F$ is any locally
compact valued field and  $U$ is an open subset of $F^d$, and say
that $\vf$ is {\it nondegenerate \/} at $\vx_0\in U$ if   the
space $F^n$ is spanned by partial derivatives of $\vf$ at $\vx_0$
   up to some finite order. (Note that the definition of $C^k$ functions
of an ultrametric variable is  more
involved than in the real case; see \S 3 for details.) 
One can view this condition as an infinitesimal version of not
lying in any proper affine hyperplane, i.e.~of the linear
independence of $1,f_1,\dots,f_n$ over $F$ (see \S 4 and \cite{KM, \S 1}
   for further discussion).

It was conjectured by \spr\ in 1980 \cite{Sp2, Conjecture H$_2$} and
proved in \cite{KM}
in 1998
that $\vf_*\lambda$ is strongly extremal for $\vf:U\to \br^n$,
$U\subset \br^d$,
which is nondegenerate at $\lambda$-a.e.\ point of $U$. Much less has
been known
for other fields. For example, the extremality of $\vf_*\lambda$ was shown
by E.\ Kovalevskaya
\cite{Ko1--2} for $\vf:\bz_p\to\bz_p^3$ which is normal in the
sense of Mahler
(a subclass of $p$-adic analytic functions) and nondegenerate $\lambda$-a.e.

In this paper we are able to prove much more general results.
The following theorem, which  we derive from Theorem 8.3, yields
an $S$-arithmetic version of
\spr's Conjecture H$_2$, in particular answering Shparlinski's question.

\proclaim{Theorem}  Let $S$ be a finite set of normailzed
valuations of $\bq$, for any $v\in S$ take $k_v,d_v\in\bn$ and an
open subset $U_v\subset \bq_{v}^{d_v}$, and let $\lambda$ be defined
as in {\rm (0.3)}.
Suppose that $\vf$ is of the
form $(\vf^{(v)})_{v\in S}$, where each $\vf^{(v)}$ is a $C^{k_v}$
map from $U_v$ into $\bq_{v}^n$ which is nondegenerate at
$\lambda_v$-a.e.\ point of $U_v$.
    Then
$\vf_*\lambda$ is strongly extremal.
\endproclaim

Note that the paper \cite{Z} considers the case when each
$\vf^{(v)}$ is of the form (0.6), and proves the extremality of
$\vf_*\lambda$.

\medskip


Theorem 0.4 is a special case of  Theorem 10.4, which 
requires a certain terminology so we do not state it
in the introduction.
In fact, Theorem 10.4 generalizes the main result of
\cite{KLW},
which, among other things, studies \di\ properties
of generic points on certain
fractal subsets of $\br^n$, in particular, the so-called
{\sl self-similar open set condition\/} fractals (see \S10.6 for details).
Following \cite{H}, those can be considered in
vector spaces over arbitrary locally compact fields, and Theorem
10.4, combined with
certain results of
\cite{KLW}, implies that natural measures supported on them are
strongly extremal.


\subhead{0.5}\endsubhead The structure of the paper is as follows.
  In
\S 1 and \S 2 we introduce
and discuss the so-called {\sl \be\/} metric spaces, {\sl Federer\/} measures
and {\sl \cag\/} functions, that is, the language in which
       our main results are stated. 
The most important examples of
good functions
are given by 
linear combinations of coordinate functions of smooth nondegenerate
maps $F^d\to F^n$.
This has been established in \cite{KM, \S 3} for $F = \br$, and  
in 
\S3 and \S4 we develop a similar theory in the ultrametric setting. 
The key ingredient of the proof (Proposition 3.3), which
hinges on combinatorics of higher order difference quotients
of $C^k$ functions of an ultrametric variable,
makes it possible to bypass the use of the Mean Value Theorem. 
 Then in \S 5 and
\S 6 we generalize a
combinatorial estimate of \cite{KM, \S 4} to the setting of  functions on a 
\be\ metric space $X$ which are good with respect to 
a Federer measure on $X$.  In \S 7  we prove
auxiliary results 
about discrete $\bz_S$-submodules of $\bq_S^m$.
The  quantitative $S$-arithmetic
nondivergence is discussed in \S 8 where, in particular, Theorem
0.1 is proved.  \S 9 is devoted to proving Theorems
0.2 and 0.3. Then in \S 10 we turn to the $S$-arithmetic
\da, giving all the definitions, stating the most general strong extremality
result (Theorem 10.4), and mentioning applications to fractal
measures. The proof
of Theorem 10.4 breaks into two special cases (when
$S$ does or does not contain the Archimedean valuation), which are treated
in \S 11 and \S 12 respectively. In both cases the argument is based on a
modification of the dynamical approach to real \da\ as developed in
\cite{KM}.
The last section of the paper lists several possible generalizations
and open questions.

\medskip

\noindent {\bf Acknowledgements:}
Part of the
work was done during the authors' collaboration at the  University
of Lyon 1,
Brandeis University
and the Max Planck Institute;
the
hospitality of these institutions is gratefully
acknowledged. A preliminary version of this paper \cite{KT}
appeared as a preprint of the MPI. The authors are thankful to G.\,A.\ Margulis
for his interest in this work and valuable comments.
Thanks are also due to Vasily Bernik, Yuri Bilu, 
Elon Lindenstrauss 
 and Barak
Weiss for helpful discussions. The first named author was supported
in part by NSF
Grants DMS-0072565 and DMS-0232463.

\heading{\S1. \be\ covering property}\endheading


\subhead{1.1}\endsubhead
For a metric space $X$, $x\in X$ and $r > 0$, we denote
by
$B(x,r)$
          the open  ball
$B(x,r) \df \{y\in X\mid \dist(x,y) < r\}$
of radius
$r$ centered at
$x$,
and by $B(x,r^+)$ the closed ball $B(x,r^+) \df \{y\in X\mid
\dist(x,y)
\le r\}\,.$
(Note that
$B(x,r^+)$ in general does not have to coincide with the closure
$\overline{B(x,r)}$
of $B(x,r)$.) 
For a
subset
$B$ of
$X$ and a
function $f:X\to {{F}}$, where $({{F}},|\cdot|)$ is a valued field,
we let
$\|f\|_{B} \df \sup_{x\in B}|f(x)|$.
If  $\mu$ is a  locally finite
Borel measure  on $X$ and
$B$ is a subset of  $X$ with $\mu(B)>0$, 
we define
         $\| f \|_{\mu,B}$ to be equal to $\|f\|_{B\,\cap\,\supp\mu}$,
which, in  case    $f$ is continuous and $B$ is
open, is the same as  the
$L^\infty(\mu)$ norm of $f|_B$, i.e.
$$\| f \|_{\mu,B} =\sup\big\{c  \bigm| \mu(\{x\in B : |f(x)| >
c\}) > 0\big\}\,.$$

We will say that a metric space $X$  is {\sl  \be\/} if there exists a
constant
$N_{\sssize X}$ such that the following holds: for any bounded subset $A$
of
$X$ and for any family $\Cal B$ of nonempty open balls in $X$ such that
$$
\forall\,x\in A \text{ is a center of some ball of }\Cal B\,,\tag 1.1
$$
there is a finite or countable subfamily $\{B_i\}$ of $\Cal B$ with
$$
1_A \le \sum_i 1_{B_i} \le N_{\sssize X}, \tag 1.2
$$
i.e.~$A\subset \bigcup_i B_i$, and the multiplicity of that subcovering is at
most $N_{\sssize X}$.

\example{Example} Suppose that  $X$ is {\sl ultrametric\/}, that is,
the non-Archimedean triangle inequality
$
\dist(x_1,x_2) \le \max_{i = 1,2}\dist(x,x_i)$ holds
for all $x,x_1,x_2\in X$.  Then any two balls in $X$ are either disjoint
or contain one another (this observation will be repeatedly used
throughout the sequel). This implies that any covering of any subset of
$X$ by balls has a subcovering of multiplicity $1$; thus any separable
ultrametric space is
\be\ with $N_{\sssize X} = 1$ (the separability of $X$ is equivalent to the
collection of all its balls being countable). \endexample

\example{1.2.\ Example} The fact that $\br^d$  is Besicovitch is the content
of Besicovitch's Covering Theorem \cite{Mt, Theorem 2.7} (the statement
``$N_\br = 2$" is known as Vitali's Covering Theorem). In fact, Besicovitch's
proof, see \cite{Mt, pp.~29--34}, can be easily generalized to give the
following
\endexample

\proclaim{1.3.\ Lemma} For a metric space $X$, define
$$
M_{\sssize X} \df \sup\left\{k\left|\aligned
&\exists\text{ balls }B_i = B(x_i,r_i),\ 1\le i \le k, \\
&\text{ such that\ }\bigcap_{i= 1}^k B_i\ne \varnothing\,,
\text{ and }x_i\notin\bigcup\Sb{j= 1}\\{j\ne i}\endSb ^k B_j\
\forall\,i\endaligned\right.\right\} 
\,,\tag 1.3
$$
and also, for $c > 1$,
$$
D_{\sssize X}(c) \df \sup\left\{k\left|\aligned
&\exists\,x\in X, \ r > 0\text{ and pairwise disjoint balls
}
\\ &B_1,\dots,B_k\text{ of radius }r\text{ contained in
}B(x,cr)\endaligned\right.\right\}
\,.\tag 1.4
$$
Then $N_{\sssize X} \le M_{\sssize X}D_{\sssize X}(8)$; hence $X$ is \be\ if
$M_{\sssize X}$ and $D_{\sssize X}(8)$ are finite.
\endproclaim

\demo{Proof} We follow the proof of Theorem 2.6 in \cite{Mt}. Take a
bounded $A\subset X$ and a family $\Cal B$ of nonempty open balls in $X$
satisfying (1.1). For each $x\in A$ pick one ball
$B\big(x,r(x)\big)\in\Cal B$. As $A$ is bounded, we may assume that
$$
R_1 = \sup_{x\in A}r(x) < \infty\,.
$$
Choose $x_1\in A$ with $r(x_1) \ge R_1/2$ and then inductively
$$
x_{j+1}\in A\ssm
\bigcup_{i = 1}^jB\big(x_i,r(x_i)\big) \quad\text{with } r(x_{j+1}) \ge
R_1/2
$$
as long as possible. Since $A$ is bounded,
the process terminates and we get a finite sequence $x_1,\dots,x_{k_1}$.

Next let
$$
R_2 = \sup\left\{r(x)\left| x\in A\ssm \bigcup_{i = 1}^{k_1}
B\big(x_i,r(x_i)\big)\right.\right\}\,.
$$
Choose $x_{k_1+1}\in A\ssm \bigcup_{i = 1}^{k_1}
B\big(x_i,r(x_i)\big)$ with $r(x_{k_1+1}) \ge R_2/2$ and again
inductively
$$
x_{j+1}\in A\ssm \bigcup_{i = 1}^jB\big(x_i,r(x_i)\big)
\quad\text{with } r(x_{j+1}) \ge R_2/2\,.
$$

Continuing this process we obtain an increasing sequence of integers
$0 = k_0 < k_1 < k_2 < \dots$, a decreasing sequence of positive
numbers $R_i$ with $2R_{i+1} \le R_i$, and a sequence of balls $B_i =
B\big(x_i,r(x_i)\big)$ with the following properties. For $j\in \bn$,
let
$
I_j = \{k_{j-1} + 1,\dots,k_j\}$. Then one has
$$
\align
r(x) < R_j/2&\qquad\text{for }j\in \bn,\ x\in A\ssm \bigcup_{i =
1}^{k_j} B_i\,,\tag 1.5a\\
R_{j}/2 \le r(x_i) \le R_i&\qquad\text{for }i\in I_j\,,\tag 1.5b\\
x_{j+1}\in A\ssm \bigcup_{i = 1}^j B_i&\qquad\text{for }j\in \bn\,, \tag 1.5c\\
x_{k}\in A\ssm \bigcup_{l\ne j}\bigcup_{i\in I_l} B_i&\qquad\text{for
}j\in \bn,\ k\in I_j\,, \tag 1.5d\\
A\subset \bigcup_{i = 1}^\infty B_i&\,.\tag 1.5e
\endalign
$$

The first three properties follow immediately from the construction.

To prove (1.5d), take  $j\in \bn$,
$k\in I_j$, $l\ne j$ and $i\in I_l$. If $l < j$, then $x_{k}\notin B_i$
by (1.5b). If $l > j$, then $x_{i}\notin B_k$ by (1.5b), and also $
r(x_i)  < r(x_k)$ by the construction, hence   $x_{k}\notin B_i$.
Finally, to verify (1.5e), observe that since $R_j\to 0$ as $j\to\infty$,
(1.5a) forces $r(x)$ to be equal to zero for any $x\in A\ssm \bigcup_{i =
1}^\infty B_i$.

Clearly (1.5e) proves the first inequality in (1.2). To establish the
second inequality, assume that a point $x\in X$ belongs to $k$ balls
$B_i$, say,
$$
x\in\bigcup_{i = 1}^k B_{m_i}\,.
$$
Using (1.5d) and (1.3), we see that the indices $m_i$ can belong to at
most $M_{\sssize X}$ different blocks $I_j$. We now claim that for any
$j\in\bn$,
$$
\#\big(I_j\cap \{m_i\mid i = 1,\dots,k\}\big) \le D_{\sssize X}(8)\,.
$$
Indeed, fix $j\in\bn$ and write $\{n_1,\dots,n_l\}\df I_j\cap
\{m_i\mid i = 1,\dots,k\}$. By (1.5b) and (1.5c), the balls
$B\big(x_{n_i},R_j/4)\big)$, $i = 1,\dots,l$, are disjoint and contained
in $B(x, 2R_j)\big)$, so the claim follows from (1.4). This proves that
$X$ is \be,  with $N_{\sssize X} \le M_{\sssize X} D_{\sssize X}(8)$.
\qed\enddemo

\subhead{1.4}\endsubhead We will use the above lemma to prove that the products
of $\br^d$ and certain ultrametric spaces are \be. Here and hereafter, the
product of two metric spaces $(X,\dist_{\sssize X})$ and $(Y,\dist_{\sssize
Y})$ will always be supplied with the product metric
$$\dist\big((x_1,y_1),(x_2,y_2)\big) =
\max\big(\dist_{\sssize X}(x_1,x_2), \dist_{\sssize
Y}(y_1,y_2)\big)\,,\tag 1.6$$ so
that balls in
$X\times Y$ are products of balls in $X$ and in $Y$. In particular, this
convention forces the product of two ultrametric spaces to be ultrametric as
well.

\proclaim{Lemma} If $Y$ is ultrametric, one has $M_{\sssize X\times
Y} = M_{\sssize X}$.
\endproclaim

\demo{Proof} Assume $M_{\sssize X\times
Y} > M_{\sssize X}$, and choose $k > M_{\sssize X}$ and balls $B_i =
B(x_i,r_i)$,
$1\le i
\le k$, in $X\times Y$  such that $\bigcap_{i= 1}^k B_i\ne \varnothing$ and
$x_i$ is not in $\bigcup_{j= 1,\,j\ne
i}^k B_j$ for each $i$. Write $B_i =  E_i
\times F_i$, where $E_i$ and
$F_i$ are projections of $B_i$ onto $X$ and $Y$ respectively. Without loss of
generality suppose that the sequence $\{r_i\}$ is
non-increasing. Since $Y$ is ultrametric and $\bigcap_{i=
1}^k F_i\ne \varnothing$, one has $F_k \subset F_i$ for all $i$. On the
other hand, since $M_{\sssize X} < k$, the center of $E_k$ mush lie in $E_i$
for some $i < k$, therefore $x_k$ mush lie in $B_i$, a
contradiction.

The converse inequality is straightforward (and not needed for our purposes).
\qed\enddemo

\proclaim{1.5.\ Corollary} If $Y$ is ultrametric, one has
$$
N_{\sssize X\times Y}
\le M_{\sssize X}D_{\sssize X}(8)D_{\sssize Y}(8)\,;\tag 1.7
$$ in
particular,  $X\times Y$ is \be\ if the three constants in the right hand side
of {\rm (1.7)} are finite.
\endproclaim

\demo{Proof} It suffices to check that $D_{\sssize X\times
Y}(8) \le D_{\sssize X}(8)D_{\sssize Y}(8)$, which is straightforward,
as, since $Y$ is ultrametric, any ball in $Y$ of radius $8r$ is a
disjoint union of at most $D_{\sssize Y}(8)$ balls of radius $r$.
\qed\enddemo

\example{1.6.\ Example} Let ${{F}}$  be a locally compact field
with a nontrivial ultrametric valuation, and let $p$ be the
number of elements in the residue class field of ${{F}}$, that is,
the number of representatives in the closed unit ball ({\sl the
ring of integers of ${{F}}$\/}) $\Cal O \df B(0,1^+)$ modulo  the
open unit ball ({\sl the valuation ideal of ${{F}}$\/}) $\Cal P
\df B(0,1)$. Without loss of generality we can, and will from now
on, {\sl normalize\/} the valuation so that the diameter of $\Cal
P$ is equal to $1/p$. (If ${F} = \bq_p$, this way one gets $ \Cal
O = \bz_p$, $\Cal P = p\bz_p$, and $|\cdot| = |\cdot|_p$, the
standard $p$-adic valuation.) Then it is easy to see that for any
$c \ge 1$,  any ball in ${{F}}$ of radius $cr$ is a disjoint union
of at most $p^{[\log_{p}c] + 1}$ balls of radius $r$.
         Therefore
Corollary 1.5, in  particular, implies that  the metric space
$$
X = \br^{d_0}\times {F}_1^{d_1}\times\dots\times {F}_\ell^{d_\ell}\tag 1.8
$$
is \be\ for any ultrametric locally compact fields
${F}_1,\dots,{F}_\ell$ and any
\linebreak $d_0,d_1,\dots,d_\ell\in
\bn$.
\endexample

\subhead{1.7}\endsubhead We close the section with a measure-theoretic
counterpart of the \be\ property.
        Namely,
say that a locally finite
Borel measure
$\mu$ on
$X$ is {\sl uniformly Federer
\/}
        if there exists
$D > 0$ such that
$$
\sup_{r > 0}
\frac{\mu\big(B(x,3r)\big)}{\mu\big(B(x,r)\big)} < D\quad\text{for
all ${x\in \supp\mu}$}\,.\tag 1.9
$$
Equivalently, one can replace ``$3$'' in (1.9) by any $c > 1$.
In other words, $\mu$ is uniformly Federer if and only if for
all
$c> 1$ one has
$$
D_{\mu}(c)\df \sup\Sb{x\in \supp\mu}\\{r > 0} \endSb
\frac{\mu\big(B(x,cr)\big)}{\mu\big(B(x,r)\big)} < \infty\,.
$$
To simplify notation, we are going to write $
D_{\mu}$ instead of
$D_{\mu}(3)$.
         Note that if $\mu$ is a uniformly Federer measure  on
$X$ with $\supp\mu = X$, for all $c> 1$  one automatically has
$D_{\sssize X}(c)
\le D_{\mu}(2c)$, or even $\le D_{\mu}(c)$ if $X$ is ultrametric.

It is often useful to have a non-uniform version of the above
definition: following
\cite{KLW}, we will say that $\mu$ as above  is {\sl Federer\footnote{See
\cite{S} and \cite{KLW, \S 6} for
an even weaker non-uniform version.}\/} if for
$\mu$-a.e.\ $x\in X$  there
exists a neighborhood
$U$ of
$x$  such that $\mu|_U$ is
        uniformly Federer.

\example{Example} Let  $X$ be as in (1.8), and denote by
$p_i$ the number of elements in the residue class field of
${F}_i$, $i = 1, \dots,\ell$.
%
It is clear
that any Haar
measure $\lambda$ on
$X$ is uniformly Federer, with $
D_{\lambda}(c)
\le c^{d_0}\prod_{i = 1}^\ell (c p_i)^{d_i}$.
        \endexample


\heading{\S 2. \cag\ functions}\endheading

\subhead{2.1}\endsubhead
Roughly speaking, a function is said to be {\sl good\/} if the set of points
where it takes small values has small measure. To simplify notation, it will
be convenient to introduce a special symbol for a set of points $x$
in a set $B$
such that the value of a function $f$ at $x$ has norm less than
$\vre$. Namely, let us define
$$
B^{f,\vre}\df \big\{ x\in B\bigm|
|f(x)| < \vre\}
$$
for any $f:B\to F$, where $({F},|\cdot|)$ is
a valued field.

   Now let $X$ be a  metric space and
$\mu$ a  Borel measure on $X$.
For a subset $U$ of $X$ and  $C,\alpha> 0$, say that
         a
function $f:X\to {F}$  is {\sl \cag\ on $U$ with respect to
$\mu$\/} if for any open ball
$B\subset U$ centered in $\supp\mu$ one has
$$
\forall\,\vre > 0\quad \mu\big(B^{f,\vre}\big) \le C\left(\frac\vre{\|f\|_{\mu,
B}}\right)^\alpha  \mu(B)\,.
\tag 2.1
$$

  In all the applications of our results, the metric
space
$X$ will be the normed ring as in  (1.8), and $\mu$  will be chosen to
be  a   Haar measure $\lambda$ on $X$, in which case we will
omit the reference to the measure and will simply say that the functions
are \cag\ on $U$.  In particular, $\mu$ will
be positive on open sets, so it will be always possible to replace
$\|f\|_{\mu,
B}$ in (2.1) by $\|f\|_{B}$ and not pay attention to the
restriction of the center of $B$ lying  in $\supp\mu$.
The above definition generalizes the one
         from \cite{KM}, which involved functions on $\br^d$,
         with $\mu$ being Lebesgue measure.
See however \cite{KLW} where measures supported on proper subsets of $\br^d$
are considered.

The following properties are immediate from the definition:

\proclaim{Lemma} Let a metric spaces $X$, a measure $\mu$ on $X$, a subset
$U$ of $ X$ and
$C,\alpha > 0$ be given.

\roster
\item"{\rm (a)}" $f$ is \cag\ on $U$  with respect to $\mu$ $\iff$ so is $|f|$;

\item"{\rm (b)}" $f$ is \cag\ on $U$ with respect to $\mu$
$\Longrightarrow$  so is $c f$
$\forall\,c\in {F}$;

\item"{\rm (c)}"  $f_i$, $i\in I$, are \cag\ on $U$ with respect to
$\mu$ $\Longrightarrow$ so is
$f = \sup_{i\in I}|f_i|$;

\item"{\rm (d)}"  $f$ is  \cag\ on $U$ with respect to $\mu$, and
$c_1\le \tfrac
{|f(x)|}{|h(x)|}\le c_2$
for all $x\in U$ $\Longrightarrow$   $h$ is
$\big(C(c_2/c_1)^\alpha,\alpha)$-good on $U$ with respect to $\mu$.

\endroster
\endproclaim

\subhead{2.2}\endsubhead The next lemma will be useful in dealing
with functions
on products of metric spaces:

\proclaim{Lemma} Let metric spaces $X,Y$ with measures $\mu,\nu$
be given. Suppose $f$ is a continuous function  on $U\times V$,
where $U\subset X$ and $V\subset Y$ are open subsets, and
suppose $C,D, \alpha, \beta$
are positive constants such that
$$
\aligned
\text{for  all $y\in V\cap\,\supp\nu$, the function }x\mapsto f&(x,y)\\
\text{ is
$(C,\alpha)$-good on $U$ with respect to $\mu$}\,,
\endaligned\tag 2.2a
$$
and
$$
\aligned
\text{for  all $x\in U\cap\,\supp\mu$, the function }y\mapsto f&(x,y)\\
\text{  is
$(D,\beta)$-good on $V$  with respect to $\nu$}
\,.\endaligned\tag 2.2b
$$
Then $f$ is $(E,\gamma)$-good on $U\times V$  with respect to $\mu\times \nu$,
         where
$$
\gamma = \frac{\alpha\beta}{\alpha + \beta}\quad\text{and}\quad E =
(\alpha +
\beta)\left(\left(\frac{C}{\beta}\right)^{\beta}\left(\frac{D}{\alpha}\right)^{\alpha}
\right)^{\frac1{\alpha +\beta}}\,.\tag 2.3
$$
\endproclaim

\demo{Proof} 
Fix a ball in
$U\times V$ of the form
$A\times B$,  where $A$ and $B$ are balls in $X$ and $Y$ intersecting
the supports
of $\mu$ and $\nu$ respectively. Without loss of generality let us
rescale $\mu|_A$,
$\nu|_B$ and $f$ so that $\mu(A) = \nu(B) = \|f\|_{\mu\times \nu, A\times
B} = 1$. Take an arbitrary $\vre > 0$; we  need to demonstrate that
$$ (\mu\times \nu)\big((A\times B)^{f,\vre} \big)  \le E
\vre^\gamma\,.\tag 2.4$$

For $y\in B$ let us denote by $f_y$ the function
$x\mapsto f(x,y)$. Also denote by $\varphi$ the function defined on $B$ by
$
\varphi(y)\df \|f_y\|_{\mu,
A}
$; note that $\|\varphi\|_{\nu,
B} = 1
$ because of the assumed normalization of $f$.
In view of (2.2a), for any $y\in B\,\cap\,\supp\nu$ one has
$$
\mu(A^{f_y,\vre}) \le
C\left(\frac{\vre}{\varphi(y)}\right)^\alpha  \iff \varphi(y) \le
\left(\frac{C}{\mu(A^{f_y,\vre})}\right)^{1/\alpha}\vre
\,. \tag 2.5
$$

Take an arbitrary $t > 0$ (to be fixed later), and denote
$$
B_t \df \big\{y\in B \mid \mu(A^{f_y,\vre}) \ge t\big\}\,.
$$
In view of (2.5),  $y\in B_t$ implies that $\varphi(y)$ is not bigger
than $\left({C}/{t}\right)^{1/\alpha}\vre$. Since it follows from
Lemma 2.1(c) and (2.2b)
that  $\varphi$
is $(D,\beta)$-good on $V$  with respect to $\nu$, one can write
$$
\nu(B_t) \le \nu\big(\big\{y\in B \mid \varphi(y) \le
\left({C}/{t}\right)^{1/\alpha}\vre\big\}\big) \le
D\left(\left({C}/{t}\right)^{1/\alpha}\vre
\right)^\beta  \,.\tag 2.6
$$

Now observe that one has $ \mu\big(\{x\in A \mid (x,y) \in (A\times
B)^{f,\vre} \}\big)  < t$
whenever $y\notin B_t$, therefore, by Fubini,
$$
\aligned
(\mu\times \nu)\big((A\times B)^{f,\vre} \big)  &<
(\mu\times \nu)\big(A\times B_t \big) + t \cdot \nu(B\ssm B_t ) \\
&\le \nu( B_t ) + t
\un{(2.6)}\le t + ( D{C}^{\beta/\alpha}\vre^\beta) \cdot
{t}^{-\beta/\alpha} \,.
\endaligned\tag 2.7
$$
The function in the right hand side of (2.7)
attains its minimum
when $$t = \left( D{C}^{\beta/\alpha}\vre^\beta
\tfrac\beta\alpha\right)^{\frac\alpha{\alpha+\beta}} = \left({C}^{\beta}
\left(\tfrac{D\beta}{\alpha}\right)^{\alpha}
\right)^{\frac1{\alpha +\beta}}\vre^{\frac{\alpha\beta}{\alpha + \beta}}\,; $$
substituting it into (2.7), one easily obtains (2.4) with $E$ and
$\gamma$ given by (2.3).
\qed \enddemo


\subhead{2.3}\endsubhead Applying the above lemma repeatedly, one
easily obtains

\proclaim{Corollary} For $j = 1,\dots,d$, let $X_j$ be a metric
space,
$\mu_j$ a measure on $X_j$, $U_j\subset X_j$ open, $C_j,\alpha_j >
0$, and let $f$ be a function
on $U_1\times\dots\times U_d$ such that for any $j = 1,\dots,d$ and
any $x_i \in U_i$ with $i \ne j$, the function
$$
y\mapsto f(x_1,\dots,x_{j-1},y,x_{j+1},\dots,x_d)\tag 2.8
$$
is $(C_j,\alpha_j)$-good on $U_j$  with respect to $\mu_j$. Then
$$f\text{ is $(\tilde C,\tilde \alpha)$-good  on
$U_1\times\dots\times U_d$ with respect to }\mu_1\times\dots\times
\mu_d\,,\tag 2.9
$$
where
$\tilde C$ and $\tilde \alpha$ are explicitly computable in terms of
$C_j$, $\alpha_j$. In particular, if each of the functions {\rm (2.8)\/}
is $(C,\alpha)$-good on $U_j$  with respect to $\mu_j$,
{\rm (2.9)\/} holds with
$$
\tilde \alpha = \alpha/d \quad\text{and}\quad\tilde C = dC\,.
$$
\endproclaim

\subhead{2.4}\endsubhead The papers \cite{KM} and \cite{BKM} describe various
classes of real valued functions on open subsets of $\br^d$ which are
\cag\ with
respect to  Lebesgue measure. For example, the fact that polynomials in
one real variable  of degree $\le k$ have
that property (with $\alpha = 1/k$) follows easily from Lagrange's
interpolation, see  \cite{DM, Lemma 4.1}  and \cite{KM,
Lemma 3.2}. Similarly, following \cite{To}, one can consider
polynomials over other locally compact fields:

\proclaim{Lemma} Let ${F}$ be either $\br$ or a locally compact ultrametric
valued
field.
Then
%
%
for any $d,k\in\bn$, any polynomial $f\in {F}[x_1,\dots,x_d]$
of degree not
greater than
$k$ is
$(C,1/dk)$-good on ${F}^d$ with respect to  $\lambda$, where
$C$ is a
constant depending only on $d$ and $k$.
\endproclaim

       \demo{Proof}
For ultrametric $F$ the case $d = 1$ is proved in  \cite{To,
Lemma 4.1}, and the general case  
immediately follows from the one-dimensional case 
and 
Corollary 2.3. Likewise, one can use \cite{KM,
Lemma 3.2} and Corollary 2.3
to establish the claim for real  polynomials of several variables.
\qed
\enddemo


\medskip

\subhead{2.5}\endsubhead Another result of \cite{KM}, which  can be
thought of as a generalization
of the real case of the previous lemma,  is that, roughly speaking,
a smooth real-valued
function is good on an
open subset of $\br^d$ (with respect to Lebesgue measure
$\lambda$) provided that its partial derivatives of some order do not vanish.

\proclaim{Lemma \rm \cite{KM, Lemma 3.3}} Given $k\in\bn$ and
an open subset $V$ of $\br^d$,   let $f\in C^k(V)$ be
such that for some constants
$\,0 <a \le A$ one has\footnote{The upper estimate in \cite{KM} was
stronger than
stated here, but in fact our weaker condition suffices for the proof.}
$$
a\le|\partial_i^{k}f(\vx)|\le A\quad \forall\,\vx\in V,\ i = 1,\dots,d\,.\tag
2.10
$$
Then $f$ is
$\,\left(dC_k\left({A}/{a}\right)^{1/k},1/{dk}\right)$-good
on $V$, where   $\br^d$ is  understood to be equipped with the $l^\infty$
metric  (induced by the norm
$\|\vx\| = \max_{i = 1}^d |x_i|$, i.e.\ the product metric on
$\br\times\cdots\times\br$
in the sense of {\rm (1.6)}), and $C_k$ is a constant dependent on $k$ only
(and explicitly estimated in \cite{KM}).
\endproclaim

Our goal in the next  section is to describe the class
of ultrametric $C^k$ functions and prove a non-Archimedean analogue
(Theorem 3.2) of Lemma 2.5.


\heading{3.
Ultrametric $C^k$ functions}
\endheading

\subhead{3.1}\endsubhead In this section we state and prove the
ultrametric analogue of the ``$d = 1$'' case of Lemma 2.5. We start by
introducing certain terminology, most of which is borrowed from
\cite{Sc}.  Here and until the end of the section
   $F$ is a complete field with a nontrivial ultrametric
valuation $|\cdot|$, and $f$  an $F$-valued function on a subset
$U$ of
$F$ without isolated points.
   The {\sl first difference quotient\/}
$\Phi^1f$ of
$f$ is the function of two variables given by
$$
\Phi^1f(x,y) \df \frac{f(x) - f(y)}{x-y} \quad(x,y\in U,\ x\ne y)\,,
$$
defined on $$\nabla^2 U\df \{(x,y)\in U\times  U\mid
x\ne y\}\,.
$$ We say that
$f$ is 
{\sl $C^1$ at $a$
\/}
if the limit $$\lim_{(x,y)\to(a,a)}\Phi^1f(x,y)$$
exists,
and that 
$f\in
C^1(U)$
if $f$ is
$C^1$ at every point of $U$.

More generally, for $k\in \bn$ set $$\nabla^k U \df
\{(x_1,\dots,x_k)\in U^{k}\mid  x_i\ne x_j \text{ for } i\ne j \}\,,$$
and define the {\sl $k$-th order difference quotient\/}
$\Phi^kf:\nabla^{k+1} U\to F$ of $f$ inductively by $\Phi^0f \df f$
and
$$\Phi^kf(x_1,x_2,\dots,x_{k+1}) \df
\frac{\Phi^{k-1}(x_1,x_3,\dots,x_{k+1}) -
\Phi^{k-1}(x_2,x_3,\dots,x_{k+1})}{x_1-x_2}\,.
$$
Note that one could equivalently take any pair of variables in place of
$(x_1,x_2)$, and that $\Phi^kf$ is a symmetric function of its $k+1$
variables. Then say that $f$ is
$C^k$ at $a$
    if the limit
$$\lim_{(x_1,\dots,x_{k+1})\to(a,\dots,a)}\Phi^kf(x_1,\dots,x_{k+1})$$
exists,  and that 
$f\in C^k(U)$
    if $f$ is
$C^k$ at every point of $U$. The latter is equivalent to  $\Phi^kf$ being
extendable to a continuous function $\bar\Phi^kf:U^{k+1}\to F$. Note that
$\nabla^{k+1} U$ is dense in $U^{k+1}$ if $U$ has no isolated points, so
the extension is unique if it exists.
%
We refer the reader to \cite{Sc, \S 27--29} for basic facts about $C^k$
functions. For instance, one can show that $C^k$ functions  $f$ are $k$
times differentiable, and in fact
$$
f^{(k)}(x)  = k!\bar\Phi^k(x,\dots,x)\,.\tag 3.1
$$
In
particular, $f\in C^k$ implies that  $f^{(k)}$ is continuous. However
the converse is not true, see  \cite{Sc, \S 27, Remark 1} for a
counterexample.  On the other hand, locally analytic functions are $C^k$
for every $k$.

\medskip

The definition of $C^k$ functions of  several
   ultrametric variables is a straighforward generalization of the one for
single-variable functions. If  $f$ is  an $F$-valued function on
$U_1 \times \cdots \times U_d$, where each $U_i$ is a subset
of
$F$ without isolated points, let us denote by $\Phi^k_if$ the $k$th
order difference quotient
    of
$f$ with respect to the variable $x_i$, and, more generally, for a
{\sl multiindex\/} $\beta = (i_1,\dots,i_d)$ let
$$\Phi_\beta f \df \Phi_1^{i_1}\circ\dots\circ \Phi_d^{i_d} f\,,$$
where it is not hard to
check that the composition can be taken in any order. The latter
``difference quotient
of order $\beta$'' is defined on $\nabla^{i_1}U_1 \times \cdots
\times \nabla^{i_d}U_d$,
and as before we say that $f$ belongs to $C^k(U_1 \times \cdots \times U_d)$
    if  for any multiindex $\beta$ with  $|\beta| \df
\sum_{j = 1}^d i_j$ at most $k$,  $\Phi_\beta f$ is
extendable to a continuous function $\bar\Phi_\beta f:
U_1^{i_1+1} \times \cdots \times U_d^{i_d + 1}\to F$.  As in the
one-variable case,
   one can show that partial derivatives $
\partial_\beta f \df \partial_1^{i_1}\circ\dots\circ \partial_d^{i_d}f$
of a  $C^k$ function  $f$ exist and are continuous
as long as $|\beta|\le k$.
Moreover, one has
$$
\partial_\beta f(x_1,\dots,x_d)  = \beta
!\bar\Phi_\beta(x_1,\dots,x_1,\dots,x_d,\dots,x_d)\,.
\tag 3.2
$$
where $\beta ! \df \prod_{j = 1}^d i_j!$, and each of the variables
$x_j$ in the right hand side of
(3.2) is repeated $i_j + 1$ times.

\subhead{3.2}\endsubhead An elementary observation, which will be
repeatedly used, is that if a function $f:U\to F$, where $U$ is an open
subset of $F^d$, is continuous at
$\vx_0\in U$ and $f(\vx_0) \ne 0$, then there exists a neighborhood $V$ of
$\vx_0$ such that
$|f(\vx)| = |f(\vx_0)|$ for all $\vx\in V$. Thus,
a natural ultrametric
replacement for
inequalities of type  (2.10)
would be assuming
that the absolute value of certain difference quotients of $f$
is identically equal to some $A > 0$  on some open set.

With this in mind, let us state 
an
ultrametric analogue of
   Lemma 2.5.

\proclaim{Theorem} Let $V_1,\dots,V_d$ be nonempty  open sets in $F$,
and let $k\in \bn$,
$A_1,\dots,A_d > 0$
and
$f\in C^k(V_1\times\dots\times V_d)$ be such that $$|\Phi_j^kf|
\equiv A_j \quad\text{on}\quad
\nabla^{k+1} V_j \times \prod_{i\ne j} V_i\,,\quad j = 1,\dots,d\,.\tag 3.3$$
Then
$f$ is
$\big(dk^{3-1/k},1/{dk}\big)$-good
on $V_1\times\dots\times V_d$.
    \endproclaim

One can immediately observe that  (3.3)  amounts to saying that
the absolute value of the $k$th order difference quotient of each of
the one-variable
functions (2.8),
$j = 1,\dots,d$, is equal to $A_j$ on $\nabla^{k+1} V_i$. Therefore one can use
Corollary 2.3 to easily derive the above theorem from its one-dimensional case.
In other words, it suffices to take an open subset $V$ of $F$, let
$k\in \bn$, $A > 0$
and
$f\in C^k(V)$ be such that
$$|\Phi^kf(x_1,x_2,\dots,x_{k+1})| = A \quad\forall\,
(x_1,x_2,\dots,x_{k+1})\in \nabla^{k+1} V\,,\tag 3.4$$
and prove that
$f$ is
$\big(k^{3-1/k},1/{k}\big)$-good
on $V$.

\medskip

The strategy of the proof will be similar to the one used in \cite{KM}
to prove  the one-dimensional case
of Lemma 2.5. However, we need to pay special attention to the following
implication of (3.4) which one gets for free in a similar situation
when $F = \br$:


\proclaim{3.3.\ Proposition} Let $V$ be a ball in $F$, and let $k\in
\bn$, $A > 0$
and
$f:V\to F$ be such that {\rm (3.4)}
holds.
Then for any $\vre
>    0$, the set
$
V^{f,\vre}
$
is a disjoint union of at most $k$ balls.
    \endproclaim

If in addition one assumes that $f\in C^k(V)$, (3.4), in
view of (3.1), would
imply that the absolute value of $
f^{(k)}(x)$ for $x\in V$ is a nonzero constant.
Note that  nonvanishing of the $k$th derivative of a real function $f$ on an
interval $V\subset \br$ immediately implies, due to the Mean Value
Theorem, that $
V^{f,\vre}$  consists of at most $k$ intervals. Unfortunately such a
theorem is not present in the ultrametric calculus, so one has to look
for alternative approaches.

\subhead{3.4}\endsubhead To prove the proposition, we will need the following
auxilliary lemma:

\proclaim{Lemma} Let $V$ be an open subset of $F$, $f$ a function
$V\to F$, $k\ge
2$, and let $x_1,\dots,x_k,y\in V$ be pairwise different. Also assume that
$$
| y - x_k| \le | x_i - x_k|\quad \forall\,i
< k\,,\tag 3.5a
$$
$$
|\Phi^{k-i}f(x_i,\dots,x_{k})| \ge |\Phi^{k-i}f(x_{i-1},\dots,x_{k-1})| \quad
\forall\,i
    = 2,\dots,k\,,\tag 3.5b
$$
and
$$
|\Phi^{k-1}f(x_1,\dots,x_{k})| \ge
|\Phi^{k-1}f(x_{1},\dots,x_{k-1},y)|\,.\tag 3.5c
$$
Then
$$
|f(y)| \le \max\big(|f(x_k)|, |f(x_{k-1})|\big) \,.\tag
3.5d
$$
\endproclaim

\demo{Proof} Note that (3.5a) implies that
$$
|y - x_1| = |y - x_k + x_k - x_1| \le |x_k - x_1|\,,
$$
and from (3.5c) one gets
$$
\aligned
&|\Phi^{k-2}f(x_2,\dots,x_{k}) -
\Phi^{k-2}f(x_1,\dots,x_{k-1})|= |x_k - x_1|\cdot|
\Phi^{k-1}f(x_1,\dots,x_{k})|\\
\ge \  &|x_k - x_1|\cdot|\Phi^{k-1}f(x_{1},\dots,x_{k-1},y)|\ge
|y - x_1|\cdot|\Phi^{k-1}f(x_{1},\dots,x_{k-1},y)| \\ = \
&|\Phi^{k-2}f(x_2,\dots,x_{k-1},y) -
\Phi^{k-2}f(x_1,\dots,x_{k-1})|\,.
\endaligned\tag 3.6
$$
Now let us use  induction on $k$. If $k = 2$, (3.6) says that
$|f(x_2) - f(x_1)| \ge |f(y) -
f(x_1)|$, which readily implies that  $
|f(y)| \le \max\big(|f(x_1)|, |f(x_{2})|\big)$. If $k > 2$ and the claim
is true with $k$ replaced by $k - 1$,  observe that (3.6) and the ``$i =
2$'' case of (3.5b) imply that
$|\Phi^{k-2}f(x_2,\dots,x_{k})| \ge
|\Phi^{k-2}f(x_2,\dots,x_{k-1},y)|$. Therefore the lemma can be applied
to
$x_2,\dots,x_k,y$, and (3.5d) follows.  \qed
\enddemo

\demo{{\bf 3.5.} Proof of Proposition 3.3} Replacing $f$ by $f/A$
without loss of
generality we may, and will, assume that $A = 1$.

Let $\vre > 0$. Note that it follows from the discreteness of the
valuation that
$
V^{f,\vre}$ is the union of finitely many balls. Assume, by
contradiction, that $
V^{f,\vre}= \cup_{i = 1}^n B_i$, where $n\ge k  + 1$ and $B_i$ are different
components of  $
V^{f,\vre}$. There exist $x_1,\dots,x_k\in
V^{f,\vre}$ such that each $x_i$ belongs to a different component (which
after changing the indices we can assume to be $B_i$) and
$$
|\Phi^{k-1}f(x_1,\dots,x_{k})|  = \sup\Sb y_i\in B_{\ell(i)} \\ i\ne j
\Rightarrow \ell(i)\ne\ell(j)  \endSb |\Phi^{k-1}f(y_1,\dots,y_{k})|\,.
$$
Next we rearrange $x_1,\dots,x_k$ in such a way that for all $\ell =
1,\dots,k-2$ one has
$$
|\Phi^{\ell}f(x_{k-\ell},\dots,x_{k})|\ge
|\Phi^{\ell}f(x_{k-\ell - 1},\dots,\check{x_i},\dots,x_{k})|\quad\text{for
all } i = k-\ell - 1,\dots,k\,,\tag 3.7a
$$
where $\check{x_i}$ means that the term ${x_i}$ is missing.

Denote
$$
R \df \min\big(|\Phi^{k-1}f(x_1,\dots,x_{k})|, |x_k -
x_1|,\dots,|x_k - x_{k-1}|)\,,\tag 3.7b
$$
and take $y\in B(x_k,R^+)$. Then, using (3.4) and $A = 1$, one writes
$$
|y - x_k| = |\Phi^{k-1}f(x_1,\dots,x_{k}) -
\Phi^{k-1}f(x_1,\dots,x_{k-1},y)| \le
|\Phi^{k-1}f(x_1,\dots,x_{k})|\,.\tag 3.7c
$$
It follows from (3.7abc) that conditions (3.5abc) hold, therefore, by
Lemma 3.4,
$$
\vre >  \max\big(|f(x_k)|, |f(x_{k-1})|\big) \ge |f(y)|\,.
$$
This proves that $B_k\supset B(x_k,R^+)$, and  from the fact
that balls $B_i$  are disjoint it follows that $|x_k - x_i| \ge
|\Phi^{k-1}f(x_1,\dots,x_{k})|$ for all $i\ne k$, hence
$
R = |\Phi^{k-1}f(x_1,\dots,x_{k})|\,.
$

Now let $y\in B_{k+1}$. By the choice of $x_1,\dots,x_{k}$ one has
$$
|\Phi^{k-1}f(x_1,\dots,x_{k})| \ge |\Phi^{k-1}f(x_1,\dots,x_{k-1},y)|\,,
$$
hence, again by (3.4),
$$
|y - x_k|  = |\Phi^{k-1}f(x_1,\dots,x_{k}) -
\Phi^{k-1}f(x_1,\dots,x_{k-1},y)| \le |\Phi^{k-1}f(x_1,\dots,x_{k})| =
R\,.
$$
Consequently $x\in B_k$, which is a contradiction.
    \qed
\enddemo

\subhead{3.6}\endsubhead Now we can proceed with the



\demo{Proof of Theorem 3.2}
We need to show that for any open ball
$B\subset V$  one has
$$
\forall\,\vre > 0\quad \lambda\big(B^{f,\vre}\big) \le k^{3-1/k}
\left(\frac\vre{\|f\|_{B}}\right)^{1/k}  \lambda(B)\,,
$$
whenever $V$ is an open subset of $f$ and $f\in C^k(V)$ satisfies (3.4).

It is clear that the result does not
depend on the normalization of $\lambda$, and it will be convenient
to assume $\lambda(\Cal O) = 1$, so that $\lambda(J) = \diam(J)$ for
any ball $J$.
In view of Proposition 3.3,
it suffices to show that for any ball $J\subset B$ with $\|f\|_{J} < \vre$
one has
$$
   \lambda(J) \le k^{2-1/k}
\left(\frac\vre{\|f\|_{B}}\right)^{1/k} \lambda(B)\,.\tag 3.8
$$
Also, as in the proof of Proposition 3.3, let us replace $f$ by $f/A$ and
thus assume  $A = 1$. Note however that now we have in addition
assumed that $f\in C^k(V)$,
therefore (3.4) implies that
   $$|\bar\Phi^kf(x_1,x_2,\dots,x_{k+1})| = 1 \quad\forall\,
x_1,x_2,\dots,x_{k+1}\in  V\,.\tag 3.9$$

It is easy to see that one can
choose
$x_1,\dots,x_{k+1}\in J$ such that
$$
|x_i - x_j| \ge 
{\lambda(J)}/{k}\qquad\text{for }
i\ne j\,.\tag 3.10$$
After that let $P$ be the Lagrange polynomial of degree $k$
formed by using values of $f$ at these points, i.e.\ given by
$$P(x) = \sum_{i = 1}^{k+1}f(x_i)
\frac{\prod_{j = 1,\,j\ne i}^{k+1}(x - x_j)}{\prod_{j = 1,\,j\ne
i}^{k+1}(x_i - x_j)}\,.
$$
Then we have $
\Phi^{k}(f-P)(x_1,\dots,x_{k+1}) = 0$, that is,
$$
\split
\Phi^{k}f(x_1,\dots,x_{k+1}) &= \Phi^{k}P(x_1,\dots,x_{k+1}) =
\text{the leading coefficient of }
P \\ &= \sum_{i = 1}^{k+1}
{f(x_i)}\Big(\prod_{j = 1,\,j\ne i}^{k+1}(x_i - x_j)\Big)^{-1}\,.
\endsplit
$$
Taking absolute values, one obtains
$$
\aligned
1 &= |\Phi^{k}f(x_1,\dots,x_{k+1})| =
|\Phi^{k}P(x_1,\dots,x_{k+1})|\\ &\le \max_{i}
\Big|{f(x_i)}{\Big(\prod_{j = 1,\,j\ne i}^{k+1}(x_i -
x_j)}\Big)^{-1}\Big| \un{(3.10) and
$\|f\|_{J} < \vre$}<
\vre\frac{  k^k}{\lambda(J)^{k}}\,.\endaligned\tag 3.11
$$

Next, take any $y\in J$ and let $Q$ be the Taylor polynomial of $f$
at $y$ of degree $k-1$. By
Taylor's formula \cite{Sc, Theorem 29.4}, for any $x$ one has
$$
f(x) =  Q(x) + (x-y)^k \bar\Phi^{k}f(x,y,\dots,y)\,,\tag 3.12$$ hence $$
\|f - Q\|_{J} \le \lambda(J)^k\|\bar\Phi^{k}\|_{J} \un{(3.9)}=
\lambda(J)^k\un{(3.11)}
\le k^k\vre\,.$$
This implies
$$
\|Q\|_{J} \le \max\big(\|f\|_{J}, \lambda(J)^k\big) <
\max\big(\vre, k^k\vre\big) = k^k\vre\,.\tag 3.13
$$
Now let us apply Lagrange's formula to reconstruct $Q$ on $B$ by its
values at $x_1,\dots,x_{k}$.
Namely, for $x\in B$ write
$$
\aligned
|Q(x)| = &\left|\sum_{i = 1}^{k}Q(x_i)
\frac{\prod_{j = 1,\,j\ne i}^{k}(x - x_j)}{\prod_{j = 1,\,j\ne
i}^{k}(x_i - x_j)}\right|\\
\un{(3.13),\,(3.10)}< &k^k\vre
{\lambda(B)^{k-1}}\frac{  k^{k-1}}{\lambda(J)^{k-1}} \le k^{2k-1} \vre
\Big(\frac{\lambda(B)}{\lambda(J)}\Big)^{k}\,.\endaligned\tag 3.14
$$
Finally, the difference between $f$ and $Q$ on $B$ is, again in view
of  (3.12) and (3.9),
bounded from above by $\lambda(B)^k$, hence
$$
\split
\|f\|_{B} &\le \max\big(\|Q\|_{B}, \lambda(B)^k\big) \un{(3.14)}<
\max\left(k^{2k-1} \vre
\Big(\frac{\lambda(B)}{\lambda(J)}\Big)^{k}, \lambda(B)^k\right) \\
&= \lambda(B)^k
\max\left( \vre
\frac{k^{2k-1}}{\lambda(J)^{k}}, 1\right)\un{(3.11)}<  k^{2k-1} \vre
\Big(\frac{\lambda(B)}{\lambda(J)}\Big)^{k}\,,
\endsplit
$$
which is equivalent to (3.8).
\qed\enddemo

\heading{\S 4 Nondegenerate, $\mu$-nonplanar and $\mu$-good   maps}\endheading

\subhead{4.1}\endsubhead
In this section
we will consider vector-valued functions of an ultrametric variable.
If $\vf =
(f_1,\dots,f_n)$ is a  map from an open subset  $U$  of
       $F^d$ into $F^n$, for any multiindex $\beta = (i_1,\dots,i_d)$
we let $$\Phi_{\beta}\vf \df
(\Phi_{\beta}f_1,\dots,\Phi_{\beta}f_n)\,,$$ and say that $\vf$ is $C^k$ if so
is each $f_i$.
In the latter case one  denotes by $\bar \Phi_{\beta}\vf$ the
continuous function
extending
  $\Phi_\beta f$
to $U_1^{i_1+1} \times \cdots \times U_d^{i_d + 1}$, so that (3.2) holds with
$f$ replaced by $\vf$.

Let us now take $F$ to be either $\br$ or an ultrametric valued
field, and say that a  map $\vf:U\to F^n$, where $U$ is an open subset
of $F^d$, is  {\sl $k$-nondegenerate at
$\vx_0\in U$\/} if  it  is $C^k$ on a neighborhood
of $\vx_0$, and the space $F^n$ is spanned by all the partial
derivatives $\partial_{\beta}\vf(\vx_0)$ of
$\vf$ at
$\vx_0$  with $|\beta| \le k$. We will say that
$\vf$ is {\sl nondegenerate at
$\vx_0$\/} if it is $k$-nondegenerate at
$\vx_0$ for some $k$. Another way of saying this is as follows: $\vf$ is
   $k$-nondegenerate at
$\vx_0
$ iff for any function $f$ of the form $f = c_0 + \vc\cdot\vf$, where
$c_0 \in F$ and
$\vc\in F^n \nz$
there exists a multiindex $\beta$ with $|\beta| \le k$
such that $\partial_\beta f(\vx_0) \ne 0$.

In particular, it follows from the nondegeneracy of $\vf$ at $\vx_0$
that for any
neighborhood
$B$ of $\vx_0$ the
restrictions of
$1,f_1,\dots,f_{n}$ to
$B$
are linearly independent over ${{F}}$; in other words, $\vf(B)$
   is not contained in any proper affine subspace  of $F^n$.
On the other hand, the converse is true under an additional
assumption that $\vf$ is analytic in a neighborhood of $\vx_0$:
  indeed, if $\vf$ can be written as a Taylor series in a neighborhood
$B$ of $\vx_0$,
and it is known that {\it all\/} partial derivatives of $\vf$ at $\vx_0$
belong to a proper subspace $L$ of $F^n$, then $\vf(B)$ must be contained
in $L + \vx_0$.

In more general situations it will be convenient to use the following
terminology:
if $X$ is a   metric space and $\mu$ a measure on $X$, a  map $\vf =
(f_1,\dots,
f_n)$ from
$X$ to
${{F}}^n$ will be called {\sl
$\mu$-nonplanar at
$x_0\in X$\/}
   if
for any
neighborhood
$B$ of $x_0$ the
restrictions of
$1,f_1,\dots,f_{n}$ to
$B\,\cap \,\supp\,\mu$
are linearly independent over ${{F}}$ $\iff$ $\vf(B \,\cap\,
\supp\,\mu)$  is not contained in any proper affine subspace  of $F^n$.
We will omit the dependence on the measure when it is taken to be
Lebesgue or Haar. Thus the above remark
translates into saying that for a $C^k$ (resp., analytic)
function $\vf:F^d\to F^n$, nondegeneracy implies (resp., is equivalent to)
nonplanarity.



\subhead{4.2}\endsubhead We are now going to discuss another property
of $\vf$ which will also be  implied by nondegeneracy. Namely, if
  $({F},|\cdot|)$ is a valued field, $X$   a metric space,
$\mu$  a measure
on $X$
and  $\vf$ a
      map
from
$X$
to
${F}^n$, let us say that  $\vf$ is {\sl $\mu$-good  at  $x_0\in
X$\/} (cf.\ \cite{K2}) if there exists a neighborhood
$V$ of $x_0$ and positive $C,\alpha$
such that
       any linear combination of
$1,f_1,\dots,f_n$  is
$(C,\alpha)$-good on
$V$ with respect to  $\mu$. Again, the reference to the measure
will be omitted
when $\mu = \lambda$.
For example, it follows from Lemma 2.4 that polynomial
maps are  good at every point.
Similarly, in \cite{KM} Lemma 2.5 was used to show that 
 a
smooth map  $\vf:\br^d\to \br^n$ is  good  at every point
where it is
nondegenerate. 


Our goal in this section is to prove an ultrametric analogue of 
the aforementioned result, using  Theorem 3.2 in place of Lemma 2.4. Namely, we have


\proclaim{Proposition}
Let $F$  be  an ultrametric valued
field, and let $\vf = (f_1,\dots,f_n)$ be a $C^\ell$ map from an open
subset $U$ of $F^d$
to $F^n$ which is $\ell$-nondegenerate at $\vx_0\in U$. Then there exists
a neighborhood
$V\subset U$ of $\vx_0$ such that  any linear combination of
$ 1,f_1,\dots,f_n$ is $\big(d\ell^{3-1/\ell},1/{d\ell}\big)$-good
on $V$.
In particular, the nondegeneracy of
$\vf$   at
$\vx_0$ implies that $\vf$ is good at
$\vx_0$. \endproclaim

\demo{Proof}
Without loss of generality we can put $\vx_0 = 0$, and
consider the family of functions
$$\Cal H \df \left\{h =  c_0 + \left.\sum_{i = 1}^nc_if_i \right|
\max_{i = 0,1,\dots,n}|c_i| = 1\right\}\,.
$$
It is enough (see Lemma 2.1(b)) to find a neighborhood $V$ of $0$ in $F^d$
such that any $h\in \Cal H$ is $\big(d\ell^{3-1/k},1/{d\ell}\big)$-good
on $V$.

 From the nondegeneracy assumption it follows that
for any
  $h \in \Cal H$
one can find a
multiindex $\beta$ with 
$$
1\le |\beta| = k \le \ell \quad\text{and}\quad|\partial_\beta h(0)| =
\big|\sum_{i = 1}^nc_i\partial_\beta f_i(\vx_0)\big| \ne 0\,.\tag 4.1
$$

Now take  $h \in \Cal H$ and
consider the
  functions $h\circ g$, where $g$ runs through the group
$\GL(d,\Cal O)$ of linear isometries of $F^d$.
(We recall that $\Cal O = \big\{x\in F \bigm| |x| \le 1\big\}$,
see \S 1.6.) For any given multiindex $\gamma$,
$\partial_\gamma (h\circ g)(0)$ is a homogeneous polynomial in matrix
elements of $g$ of degree $\le |\gamma|$ with coefficients given by
$\partial_{\gamma'} f(0)$ where  $|\gamma'| = |\gamma|$. It follows from (4.1)
that for any $\gamma$ with $ |\gamma|= k$ this polynomial is nonzero.
Hence it is possible to choose $g$ so that $\partial_\gamma (h\circ g)(0)\ne 0$
for {\it all\/} multiindices $\gamma$ with $ |\gamma|=k$. In fact, we
are only interested
in choosing  $g$ with
$$\partial_i^k (h\circ g)(0)\ne 0\quad\text{for each
  }i = 1,\dots,d\,.\tag 4.2
$$
Using (4.2) and the compactness of
both   $\Cal H$ and $\GL(d,\Cal O)$, 
one can find a ball $V = V_1\times\dots\times V_d\ni 0$ in $F^d$
(here $V_i$ are balls in $F$ of the same radius) such that
  for any
  $h \in \Cal H$ there exist  $1\le k \le \ell$,
$g\in \GL(d,\Cal O)$ and $A_1,\dots,A_d\in F\nz$ such that
(3.3) holds for $f = h\circ g$. Therefore, by Theorem 3.2, for any
  $h \in \Cal H$ one can find   $1\le k \le \ell$ and
$g\in \GL(d,\Cal O)$ such that for any ball $B\subset V$ and $\vre > 0$ one has
$$
  \lambda\big(B^{h\circ g,\vre}\big) \le dk^{3-1/k}
\left(\frac\vre{\|h\circ g\|_{B}}
\right)^{1/{dk} } \lambda(B)\le d\ell^{3-1/\ell}
\left(\frac\vre{\|h\circ g\|_{B}}
\right)^{1/{d\ell} } \lambda(B)\,.
$$
To finish the proof it remains to notice that $g$ leaves $V$
invariant, sends balls to balls,
and one clearly has $B^{h\circ g,\vre} = g(B)^{h,\vre}$ and
$\|h\circ g\|_{B} = \|h\|_{g(B)}$.
\qed\enddemo



\subhead{4.3}\endsubhead For convenience let us summarize the results
of this section as follows:

\proclaim{Theorem}
Let $F$  be either $\br$ or an ultrametric valued
field, and let $\vf$ be a $C^\ell$ map from
an open subset $U$ of $F^d$
to $F^n$. Then $\vf$ is nonplanar and good at every point of $U$ where
it is nondegenerate. \endproclaim

\heading{\S 5. Maps of posets into spaces of good functions}\endheading

\subhead{5.1}\endsubhead The goal of this section is to generalize  a
construction described in
\cite{KM} in order to make it work for functions defined on
arbitrary metric spaces.
More precisely, we will work with  mappings of partially
ordered sets ({\sl posets\/}) $\frak P$
into spaces of functions on a metric space
$X$ with a measure $\mu$. Given  such a mapping, we will {\it
mark\/} certain points  (see the definition below), and prove an upper
estimate (Theorem 5.1) for  the measure of the set of
``unmarked" points\footnote{A possibility of such a generalization is
mentioned in \cite{KM, \S 6.1}. The paper \cite{KLW} contains a slightly
different presentation of the same argument, written in the special
case of $\frak P$ being the poset of nonzero rational subspaces of
$\br^m$.}.

         For a poset $\frak P$, we will denote by $l(\frak P)$ the
{\sl length\/}
of $\frak P$, i.e.~the number of elements in a maximal linearly
ordered subset of
$\frak P$. If $\frak S$ is a subset of $\frak P$, we let $\frak P(\frak S)$ be
the poset of elements of
$\frak P\smallsetminus \frak S$ comparable with any element of $\frak S$. Note
that one always has
$$
l\big(\frak P(\frak S)\big) \le l(\frak P) - l(\frak S)\,.\tag 5.1
$$

We will fix a metric space $X$, and consider posets $\frak P$ together
with  a mapping  $\psi$ from
$\frak P$ to the space $C(B)$ of $\br$-valued
      continuous
functions on some subset $B$ of
$X$, to be denoted by
$s\mapsto \psi_s$. Given such a mapping and positive numbers $\vre \le
\rho$, we will say that a point
$z\in B$  is {\sl  $(\vre,\rho)$-marked relative to $\frak P$\/} if
there exists a linearly ordered subset
$\frak S_z$ of $\frak P$ such that

\roster
\item"(M1)" $\vre\le |\psi_s(z)| \le \rho \quad \forall\,s\in\frak S_z$;

\item"(M2)" $|\psi_s(z)| \ge \rho\ \ \,\qquad\forall\,s\in \frak P(\frak S_z)$.
\endroster
We will denote the set of all such points by $\p(\vre,
\rho,{\frak P})$.
When it does not cause confusion, we will omit the
reference to either $\frak P$ or $(\vre,\rho)$, and will simply say that
$z$  is
$(\vre,\rho)$-marked, or marked relative to $\frak P$.

         \proclaim{Theorem}  Let $X$ be a 
\be\ metric space,
$\mu$ a  uniformly Federer measure  on $X$,
$m\in
\bz_+$ and
$C,\alpha,\rho  > 0$.
Suppose that we are given a poset $\frak P$, a ball $B = B(x,r)$ in $X$, and
          a mapping $\psi:\frak P\to C(\tilde B)$, where $\tilde B\df
B\big(x,3^mr\big)$,  such that the
following holds:

\roster
\item"(A0)" $l(\frak P) \le m$;
\item"(A1)" $\forall\,s\in \frak P\,,\quad \psi_s$ is \cag\ on $\tilde B$ with
respect to
$\mu$;
\item"(A2)" $\forall\,s\in \frak P\,,\quad\|\psi_s\|_{\mu,B} \ge \rho
$;
\item"(A3)"  
$\forall\,y\in \tilde
B \,\cap\,\supp\mu,\quad\#\{s\in \frak P\bigm|
|\psi_s(y)| < \rho\} < \infty$.
\endroster
Then
$
\forall\,\vre \le \rho$ one has
$$
\mu\big(B\smallsetminus \p(\vre,
\rho,{\frak P})\big) \le mC
\big(N_{\sssize X}D_{\mu}^2\big)^m
\left(\frac\vre \rho \right)^\alpha
\mu(B)\,.
$$
\endproclaim

\demo{Proof}  We proceed by induction on $m$. If
$m=0$, the poset $\frak P$ is empty, and for any
$z\in B$  one can take
$\frak S_z = \vrn$ and check that (M1) and (M2) are satisfied for all
$\vre,\rho$; thus all points of $B$ are marked.  Now
take
$m\ge 1$ and  suppose that the claim is proved for all smaller values of $m$.

Fix $C,\alpha,\rho,  \frak P, B = B(x,r)$ and $\psi$  as in the
formulation of
the theorem.   For any  $y\in B \,\cap\,\supp\mu$ define
$$
H(y)\df\{s\in \frak P\bigm| |\psi_s(y)| < \rho \}\,;
$$
this is a finite subset of $\frak P$ in view of (A3). If
$H(y)$ is empty, $y$ is clearly $(\vre,\rho)$-marked  for any
positive
$\vre$: indeed, since  $|\psi_s(y)| \ge
\rho
$ for all $s\in \frak P$, one can again take $\frak S_y$ to be the
empty set and
check that (M1) and (M2) are satisfied. Thus one only needs to consider
points
$y$ from the set
$$
E\df \{y\in B \,\cap\,\supp\mu \mid H(y)\ne\vrn\} = \{y\in B
\,\cap\,\supp\mu\mid
\exists\,s\in \frak P
\text{ with }|\psi_s(y)| < \rho \}\,.
$$

Take $y\in E$ and $s\in H(y)$, and define
$$
r_{s,y}\df\sup\{t > 0\bigm| \|\psi_s\|_{\mu,B(y,t)} < \rho \}\,.\tag 5.2
$$
         It follows from the continuity of functions $\psi_s$ that
for small enough positive $t$ one has $\|\psi_s\|_{\mu,B(y,t)} < \rho $,
hence
$r_{s,y} > 0$. Denote $B(y,r_{s,y})$ by
$B_{s,y}$. From (A1) it is clear
that $B_{s,y}$ does not contain $B$; therefore one has
$r_{s,y} < 2r$.
Note also that (5.2) immediately
implies that
$$\|\psi_s\|_{\mu,B_{s,y}}\le\rho\,.\tag 5.3$$


Now for any $y\in E$ choose an element $s_y$ of $ H(y)$ such that
$r_{s_y,y}
\ge   r_{s,y}$ for all $s\in H(y)$ (this can be done since $H(y)$ is finite).
For brevity let us denote $r_{s_y,y}$ by $r_y$ and $B_{s_y,y}=
\bigcup_{s\in H(y)}B_{s,y}$ by $B_y$. Also let us denote the poset
$\frak P(\{s_y\})$ by $\frak P_y$.

      \medskip
\noindent{\bf 5.1.1.} The next lemma allows one to show  a
point
$z\in B_y$ to be
marked relative to $\frak P$ once it is
marked relative to $\frak P_y$. It is proved by a verbatim repetition
of the proof of \cite{KM, Lemma 4.6}, yet we do it in full detail here
to make the argument self-contained.

\proclaim{Lemma} For  $\vre \le \rho$ and $y\in E$, let
$z\in B_y\,\cap\,\supp\,\mu\,\cap \,\p(\vre,
\rho,{\frak P_y})$ be such that
$|\psi_{s_y}(z)| \ge\vre$; then $z$ belongs to $\p(\vre,
\rho,{\frak P})$. Equivalently,
$$
(B_y\cap\,\supp\mu)\ssm\p(\vre,
\rho,{\frak P}) \subset
\big(B_y\ssm\p(\vre,
\rho,{\frak P_y}) \big) \cup (B_y)^{\psi_{s_y},\vre}\,.\tag 5.4
$$
\endproclaim

\demo{Proof} By definition of $\p(\vre,
\rho,{\frak P_y})$, there exists
a linearly ordered subset $\frak S_{y,z}$ of $\frak P_y$ such that
$$
\vre\le |\psi_s(z)| \le \rho\quad \forall\,s\in\frak S_{y,z}\tag 5.5
$$
and
$$
|\psi_s(z)| \ge \rho\quad
\forall\,s\in \frak P_y(\frak S_{y,z})\,.\tag 5.6
$$
Put $\frak S_{z}\df\frak S_{y,z}\cup\{s_y\}$. Then
$\frak P(\frak S_{z}) = \frak P_y(\frak S_{y,z})$; therefore (M2)
immediately follows from
(5.6), and, in view of (5.5),  it remains to check (M1) for $s = s_y$. The
latter is straightforward:
$|\psi_{s_y}(z)|$ is not less than $\vre$ by the assumption and
is not greater than $\rho $ in  view of (5.3). \qed\enddemo

      \medskip
\noindent{\bf 5.1.2.}
Note that one clearly has $r_{y} <
2r$, which in particular implies that $B_y\subset B(x,3r)$. We are  going
to  fix some
$r'_y$ strictly between
$r_y$ and
$\min(2r, 3r_{y})$, and denote $B(y,r'_{y})$ by $B'_y$. Clearly one
has
$$
\|\psi_s\|_{\mu,B'_{y}}  \ge \rho
\quad\text{for any }y\in E\text{ and }s\in \frak P\,.\tag 5.7
$$
(Indeed, the definition of $r_y$ and
        (5.2) imply the above inequality for any  $s\in H(y)$,
and it obviously holds if $s\notin H(y)$.)

        Now observe that
$\frak P_y$,  $B'_y$ and $\tilde B'_y\df B(y,3^{m-1}r'_y)$ satisfy
properties
\roster
\item"$\bullet$" (A0) with $m$ replaced by $m-1\qquad$ --- $\qquad$in view of
(5.1);
\item"$\bullet$" (A2)
\hskip 1.9in ---
$\qquad$in view of (5.7);
\item"$\bullet$" (A1) and (A3) \hskip 1.31in ---
$\qquad$since
\endroster
$$
\tilde B'_y =  B(y,3^{m-1}r'_y)\subset B(x,3^{m-1}r'_y + r)
\subset B\big(x,(2\cdot 3^{m-1} + 1)r\big)\subset B(x,3^mr) = \tilde B\,.
$$

Therefore one has
$$
\aligned
\mu\big(B_y\ssm\p(\vre,
\rho,{\frak P_y}) \big) &\le \mu\big(B'_y\ssm\p(\vre,
\rho,{\frak P_y}) \big)\\  &\le(m-1)C\big(N_{\sssize
X}D_{\mu}^2\big)^{m-1} \left(\frac{\vre}\rho \right)^\alpha
\mu(B'_y)\\ &\le D_{\mu}(m-1)C\big(N_{\sssize
X}D_{\mu}^2\big)^{m-1} \left(\frac{\vre}\rho \right)^\alpha
\mu(B_y)
\endaligned\tag 5.8
$$
         by the induction assumption and the Federer property of $\mu$. On the
other hand, in view of
$\psi_{s_y}$ being \cag\ on $\tilde B\supset B'_y$,  one can write
$$
\aligned
\mu\left((B_y)^{\psi_{s_y},\vre}\right) &\le
\mu\left((B'_y)^{\psi_{s_y},\vre}\right) \le
C\left(\frac\vre{\|\psi_{s_y}\|_{\mu,B'_y}}\right)^\alpha  \mu(B'_y)
\underset{(5.7)}\to\le C  \left(\frac{\vre}\rho \right)^\alpha
\mu(B'_y)\\ &\le CD_{\mu} \left(\frac{\vre}\rho \right)^\alpha
\mu(B_y)\,.\endaligned\tag 5.9
$$
Recall that we need to estimate the measure of
$E\ssm \p(\vre,
\rho,{\frak P})$. For any $y\in E$, in view of  (5.4), (5.8) and (5.9) one
has
$$
\aligned
\mu\big(B_y\ssm\p(\vre,
\rho,{\frak P})\big) &\le
C\left((m-1)N_{\sssize
X}^{m-1}D_{\mu}^{2m-1} +
D_{\mu} \right)\left(\frac{\vre} \rho \right)^\alpha  \mu(B_y)%
\\ &\le
mCN_{\sssize
X}^{m-1}D_{\mu}^{2m-1}\left(\frac{\vre} \rho
\right)^\alpha  \mu(B_y)\,.
\endaligned
\tag 5.10
$$

Now consider the covering $\{B_y\mid y\in E\}$ of $E$, choose a
countable subset $Y$ of $E$ such that the multiplicity of the subcovering
$\{B_y\mid {y\in Y}\}$ is   at most $N_{\sssize X}$, and write
$$
\sum_{y\in Y}\mu(B_y) \le N_{\sssize X}\mu\big(\bigcup_{y\in Y} B_y\big)
\le N_{\sssize X}\mu\big(B(x,3r)\big) \le N_{\sssize
X}D_{\mu}\mu(B)\,.\tag 5.11
$$
Therefore the measure of $E\ssm \p(\vre,
\rho,{\frak P})$ is bounded from above by
$$
\split
\sum_{y\in Y}\mu
\big(B_y\ssm\p(\vre,
\rho,{\frak P})\big)
&\underset{(5.10)}\to\le  mCN_{\sssize
X}^{m-1}D_{\mu}^{2m-1}\left(\frac{\vre} \rho \right)^\alpha
\sum_{y\in Y}\mu(B_y)\\ &\underset{(5.11)}\to\le  mC\big(N_{\sssize
X}D_{\mu}^2\big)^{m}\left(\frac{\vre} \rho \right)^\alpha
\mu(B)\,.\quad\qed
\endsplit
$$
\enddemo


\heading{6. Primitive submodules of $\Cal D^m$}
\endheading

\subhead{6.1}\endsubhead We start this section by assuming that

\roster
\item"$\bullet$"  $\Cal D$ is an integral domain, that is, a  commutative
ring with
$1$ and without zero divisors;
\item"$\bullet$" $K$ is the quotient field of $\Cal D$;
\item"$\bullet$" ${\Cal R}$ is a commutative ring containing
${K}$ as a subring.
\endroster

We need the following elementary 
lemma:

\proclaim{Lemma} Let $k,m\in\bn$,  $k\le m$, and let
$\gamma_1,\dots,\gamma_k\in {K}^m$
be linearly independent over ${K}$. Then they are   linearly
independent over ${\Cal R}$.
\endproclaim

\demo{Proof} Let $A$ be the $m\times k$-matrix
with columns given by  $\gamma_1,\dots,\gamma_k$. Then there exists at least
one $k\times k$ minor $B$  of $A$ with $\det(B)$ being a nonzero element
of ${K}$, hence invertible in ${\Cal R}$. By Cramer's rule, for any solution
$\beta = (\beta_1,\dots,\beta_k)\in {\Cal R}^k$
of $A\beta = 0$ one must have $\det(B)\beta_i = 0$ for every $i$,
hence $\beta = 0$.
         \qed \enddemo

\subhead{6.2}\endsubhead If $\Delta$ is an ${\Cal D}$-submodule of
${\Cal R}^m$,
let us denote by ${K}\Delta$ (resp.~${\Cal R}\Delta$)
         its ${K}$- (resp.~${\Cal R}$-) linear span inside ${\Cal R}^m$,
and define
the {\sl rank\/} $\rk\Delta$ of $\Delta$ by
$$\rk(\Delta) \df\dim_{K}({K}\Delta)\,.\tag 6.1
$$
For example, one has $\rk({\Cal D}^m) = m$ for any $m\in\bn$.
If $\Lambda$ is an ${\Cal D}$-submodule of ${\Cal R}^m$ and $\Delta$
is a submodule
of $\Lambda$,
say that
$\Delta$ is {\sl primitive in\/}  $\Lambda$ 
if any submodule  of
$\Lambda$ of rank equal to $\rk(\Delta)$ and containing  $\Delta$ is equal to
$\Delta$. It is clear that the set of nonzero
primitive submodules of a fixed ${\Cal D}$-submodule
$\Lambda$ of ${\Cal R}^m$ is a partially ordered set (with respect to
inclusion)
of length equal to $\rk(\Lambda)$.

The next lemma  characterizes primitive submodules of ${\Cal D}^m$:

\proclaim{Lemma} The following are equivalent for a submodule
$\Delta$ of ${\Cal D}^m$:

{\rm (i)} $\Delta$ is primitive;

{\rm (ii)} $\Delta = {K}\Delta \,\cap\, {\Cal D}^m$;

{\rm (iii)} $\Delta = {\Cal R}\Delta \,\cap\, {\Cal D}^m$ for any
commutative ring
${\Cal R}$ containing
${K}$ as a subring.
\endproclaim

\demo{Proof} If $\Delta = \{0\}$, the claim is trivial.
Otherwise, it is obvious that (iii)$\Rightarrow$(ii)$\Rightarrow$(i).
Assuming (i)
and taking $\gamma\in {\Cal R}\Delta \,\cap\, {\Cal D}^m$, let
$\gamma_1,\dots,\gamma_k\in \Delta$
be a basis of ${K}\Delta$, with $k = \rk(\Delta)$.
Then $\gamma,\gamma_1,\dots,\gamma_k$ are linearly
dependent over ${\Cal R}$, hence, in view of Lemma 6.1, over ${K}$.
But $\gamma_1,\dots,\gamma_k$
are linearly independent over ${K}$, thus $\gamma$ belongs to ${K}\Delta$,
therefore the ${\Cal D}$-module $\Delta'$ generated by $\Delta$ and
$\gamma$ has
rank $k$. By primitivity of $\Delta$, $\Delta' = \Delta$,
i.e.~$\gamma\in\Delta$. \qed \enddemo

In fact, Lemma 6.2 implies that for any $\Delta'\subset {\Cal D}^m$ there
exists a unique
primitive $\Delta \supset \Delta'$ of the same rank, namely, $\Delta =
         {K}\Delta' \,\cap\, {\Cal D}^m$.

\subhead{6.3}\endsubhead Let us now assume in addition that ${\Cal
R}$ is a topological ring, and consider the topological group
$\GL(m,{\Cal R})$ of automorphisms of ${\Cal R}^m$, i.e.~the group
of $m\times m$ invertible matrices with entries in ${\Cal R}$. Any
$g\in\GL(m,{\Cal R})$ maps ${\Cal D}$-submodules of ${\Cal R}^m$
to ${\Cal D}$-submodules of ${\Cal R}^m$, preserving their rank
and the inclusion relation. Let us introduce the following
notation:
$$
\frak M({\Cal R},{\Cal D},m) \df \{g\Delta \mid g\in\GL(m,{\Cal R}), \
\Delta\text{ is a submodule of } {\Cal D}^m\}\,,\tag 6.2
$$
and
      $$\frak P({\Cal D},m)\df \text{ the set of all
nonzero primitive submodules  of }
{\Cal D}^m\,.$$
Note that the inclusion relation makes $\frak P({\Cal D},m)$  a poset of length
$m$.
\medskip

  We would like to have a way to measure ``size'' of
submodules from the above collection. Specifically, let us say
         that a function $\nu:
\frak M({\Cal R},{\Cal D},m) \mapsto \br_+$ is {\sl norm-like\/} if
the following three
properties hold:

\roster
\item"(N1)"  For any  $\Delta,\Delta'\in\frak M({\Cal R},{\Cal D},m)$ with
$\Delta'\subset \Delta$
and $\rk(\Delta') = \rk(\Delta)$ one has \newline $\nu(\Delta')\ge
\nu(\Delta)$;

\item"(N2)" there exists $C_\nu > 0$ such that for any
$\Delta\in\frak M({\Cal R},{\Cal D},m)$
and any $\gamma\notin {\Cal R}\Delta$ one has
$\ \nu(\Delta + {\Cal D}\gamma) \le
C_\nu\cdot\nu(\Delta)\nu({\Cal D}\gamma)$;
\item"(N3)"  for every submodule $\Delta$ of ${\Cal D}^m$, the function
$ \GL(m,{\Cal R})\to \br_+$, $g\mapsto \nu(g\Delta)$, is
continuous.
\endroster

If $\nu$ is as above and $\gamma\in {\Cal R}^m$,
we will define $\nu(\gamma)$ to be equal to $\nu({\Cal D}\gamma)$.
The model example is given by taking ${\Cal D}  = \bz$, ${K} = \bq$,
${\Cal R} = \br$. Then
the set $\frak M({\Cal R},{\Cal D},m)$ coincides with the set of all discrete
subgroups of $\br^n$,
and one can define $\nu(\Delta)$ to be the covolume of $\Delta$ in $\br\Delta$,
with $\nu(\vv)$ being equal to
the Euclidean norm of a vector $\vv\in\br^m$;
in that case one can easily
         check that (N1)--(N3) are satisfied,  with $C_\nu = 1$. In the next
section we will do this   in a more general context, when ${\Cal R}$
is not a field anymore.

Now we can apply Theorem 5.1 to the poset $\frak P({\Cal D},m)$.

         \proclaim{Theorem} Let $X$ be a 
\be\ metric space,
$\mu$ a  uniformly Federer measure  on $X$, and let ${\Cal D}\subset {K}
\subset {\Cal R}$ be as above,
${\Cal R}$ being a topological ring. For  $m\in \bn$, let  a
ball $B = B(x_0,r_0)\subset X$ and a continuous map $h:\tilde B \to
\GL(m,{\Cal R})$
be given, where $\tilde B$ stands for $B(x_0,3^mr_0)$. Also let $\nu$ be a
norm-like function on $\frak M({\Cal R},{\Cal D},m)$.
%
%
For any $\Delta\in \frak P({\Cal D},m)$ denote by $\psi_\Delta$ the
function
$x\mapsto \nu\big(h(x)\Delta\big)$ on  $\tilde B$.
Now suppose for some
$C,\alpha > 0$ and $0 < \rho  < 1/C_\nu$  one has

{\rm(i)} for every $\Delta\in \frak P({\Cal D},m)$, the function $\psi_\Delta$
is \cag\ on
$\tilde B$  with respect to
$\mu$;

{\rm(ii)}  for every $\Delta\in \frak P({\Cal D},m)$, $\|\psi_\Delta\|_{\mu,B}
\ge \rho$;

{\rm(iii)}   $\forall\,x\in \tilde
B \,\cap\,\supp\mu,\quad\#\big\{\Delta\in\frak P({\Cal D},m)\bigm|
\psi_\Delta(x) <
\rho\big\} < \infty$.

\noindent Then 
         for any  positive $ \vre \le \rho$ one has
$$
\mu\left(\left\{x\in B\left|\aligned \nu \big(h(x)\gamma\big) <  \vre
\text{ for \ }\\ \text{some }\gamma\in {\Cal D}^m\nz
\endaligned\right.\right\}\right)\le mC
\big(N_{\sssize X}D_{\mu}^2\big)^m
\left(\frac\vre \rho \right)^\alpha
\mu(B)\,.\tag 6.3
$$
\endproclaim

\demo{Proof} For simplicity let us denote $\frak P({\Cal D},m)$ by
$\frak P$. As was observed above, the length of
$\frak P$ is equal to
$m$, and one immediately
verifies that conditions (i)--(iii) imply that $\frak P$, $B$ and $\tilde B$
satisfy properties (A1)--(A3) of Theorem 5.1.  Thus it suffices to
prove that for any  positive
$ \vre \le \rho$ one has
$$
\p(\vre,
\rho,{\frak P}) \subset  \big\{x\in B\mid \nu\big(h(x)\gamma\big) \ge
\vre\text{ for all }\gamma\in {\Cal D}^m\nz\big\}\,.\tag 6.4
$$

Take an $(\vre,
\rho)$-marked point $x\in B$, and let $\{0\} =
\Delta_0 \subsetneq \Delta_1 \subsetneq\dots \subsetneq\Delta_l =
{\Cal D}^m$ be all the elements of $\frak S_x \cup
\big\{\{0\},{\Cal D}^m\big\}$. Pick any $\gamma\in {\Cal D}^m\nz$. Then there
exists $i$, $1 \le i \le l$, such that $\gamma\in \Delta_i\ssm
\Delta_{i-1}$. From the primitivity of $\Delta_{i-1}$ and Lemma 6.2 it
follows that
$\gamma\notin {\Cal R}\Delta_{i-1}$, hence
$g\gamma\notin g{\Cal R}\Delta_{i-1} = {\Cal R}g\Delta_{i-1}$ for any
$g\in\GL(m,{\Cal R})$.
Therefore, if one defines
$\Delta' \df {\Cal D}\Delta_{i-1} + {\Cal D}\gamma$, in view of  (N2) one has
$$\nu\big(h(x)\Delta'\big) \le C_\nu
\nu\big(h(x)\Delta_{i-1}\big)\nu\big(h(x)\gamma\big)\,.\tag 6.5
$$
Further, let $\Delta\df {K}\Delta'\,\cap\, {\Cal D}^m$. It is a primitive
submodule containing $\Delta'$ and of rank equal to $\rk(\Delta')$,
so, by (N1),
$$\nu\big(h(x)\Delta\big) \le \nu\big(h(x)\Delta'\big)\,.\tag 6.6
$$
Moreover, it is also contained in $\Delta_i$, since
$$
\Delta = {K}\Delta'\,\cap\, {\Cal D}^m =  {K}\Delta\,\cap\, {\Cal D}^m \subset
{K}\Delta_i\,\cap\, {\Cal D}^m = \Delta_i\,.
$$
Therefore it is comparable to any element of
$\frak S_x$, i.e.~belongs to $\frak S_x \cup \frak P(\frak S_x)$.
         Then one can use
properties (M1) and  (M2) to deduce
that
$$
|\psi_{\Delta}(x)| = \nu\big(h(x)\Delta\big) \ge \min(\vre,
\rho) =
\vre\,,
$$
and then, in view of (6.5) and (6.6), conclude that
$$
\nu\big(h(x)\gamma\big) \ge
{\nu\big(h(x)\Delta\big)}/{C_\nu \nu\big(h(x)\Delta_{i-1}\big)}  \ge
\vre/C_\nu\rho  \ge
\vre\,.
$$
This shows (6.4) and completes the proof of the theorem. \qed\enddemo

\heading{7. Discrete  submodules of $\bq_S^m$}
\endheading

\subhead{7.1}\endsubhead  The goal of this section is to describe a certain
class of triples ${\Cal D}\subset {K}\subset {\Cal R}$ and construct a
norm-like function on $\frak M(\Cal R,
\Cal D,m)$ which is important in applications to both dynamics and \da.

We
let $\ell \in \bn$ and take $S = \{p_1,\dots,p_{\ell-1},\infty\}$,
where $p_1,\dots,p_{\ell-1}$ are primes. The (possibly empty)
subset $ \{p_1,\dots,p_{\ell-1}\}$ of $S$ will be denoted by
$S_f$.
To every element of $S$ we associate the {\sl normalized\/}
valuation $|\cdot|_v$ of ${\bq}$; in other words, $|\cdot|_v$ is
the usual absolute value if $v = \infty$,
and is defined as in Example 1.6 if $v$ is $p$-adic. We let
${\bq}_S$ be the direct product of all the completions ${\bq}_v$,
$v\in S$, in which ${\bq}$ is diagonally imbedded (here we use the
notation $\bq_\infty = \br$), and  let
$$\bz_S \df \bz\big[\tfrac1{p_1},\dots,\tfrac1{p_{\ell-1}}\big] =
\big\{x\in {\bq}\bigm|
x\in\bz_p\text{ for all primes }p\notin S_f\big\}$$
stand for the ring of $S$-integers of ${\bq}$. We also let
$${\bz}_{S,f} \df
\prod_{i = 1}^{\ell-1}{\bz}_{p_i}\quad\text{and}\quad{\bq}_{S,f} \df
\br \times {\bz}_{S,f}\,.$$

Denote by $\lambda_v$
the {\sl normalized\/} 
Haar measure on ${\bq}_v$ (that is, $\lambda_\infty$ is the usual
Lebesgue measure on $\br$ and $\lambda_{p_i}$ is normalized by
$\lambda_{p_i}(\bz_{p_i}) = 1$), and by $\lambda_S = \prod_{v\in
S} \lambda_v$ the product measure on ${\bq}_S$. Elements of
${\bq}_S$ will be denoted as $x = (x^{(v)})_{v\in S}$ or simply $x
= (x^{(v)})$, where $x^{(v)}\in {\bq}_v$. For $x$ of this form, we
define  the $S$-adic {\sl absolute value\/} $|\vx|$ and the {\sl
content\/} $c(x)$ of $x$ to be the maximum (resp.\ the product) of
all $|x^{(v)}|_v$, $v \in S$.
Since all the valuations are normalized, one has
$$
\lambda_S(xM) = c(x)\lambda_S(M)\,.\tag 7.1
$$
where $x \in {\bq}_S$ and $M$ is a measurable subset of ${\bq}_S$.

If $m$ is a natural number, we preserve the same notation
$\lambda_v$ and $\lambda_S$ to denote the product measure on
${\bq}_v^m$ and ${\bq}_S^m$, respectively. Elements $\vx =
(x_1,\dots,x_m)$ of ${\bq}_S^m$ will be denoted as $\vx =
(\vx^{(v)})$, where $\vx^{(v)}= (x^{(v)}_1,\dots,x^{(v)}_m)\in
{\bq}_v^m$. We denote by $\|\cdot\|_v$ the
the usual (Euclidean) norm on $\br^m$ if $v = \infty$, and the
sup-norm defined by
$$\|(x^{(v)}_1,\dots,x^{(v)}_m)\|_v = \max_i|x^{(v)}_i|_v$$
if $v$ is non-Archimedean. For $\vx = (\vx^{(v)})$ in ${\bq}_S^m$
we define
   the {\sl norm\/}
$\|\vx\|$  and the {\sl content\/}  $c(\vx)$ of $\vx$  to be the
maximum (resp., the product)
   of all the numbers
$\|(\vx^{(v)})\|_v$, $v \in S$. The group
$\GL(m,{\bq}_S)=\prod_{v\in S}\GL(m,{\bq}_v)$ acts naturally on
${\bq}_S^m$,  and one has
$$
\lambda_S(gM)=c\big(\det(g)\big)\lambda_S(M),
$$
where $M\subset {\bq}_S^m$ is any measurable subset of
${\bq}_S^m$, $g=(g^{(v)}) \in \GL(m,{\bq}_S)$, and $\det(g) \df
\big(\det(g^{(v)})\big)$ is an invertible element of $\bq_S$.

\subhead{7.2}\endsubhead    Our goal now is to consider discrete
${\bz}_S$-submodules $\Delta$ of ${\bq}_S^m$. It turns out that
any such $\Delta$ is a finitely generated  free ${\bz}_S$-module:

\proclaim{Proposition} Let $\Delta$ be a  discrete
${\bz}_S$-submodule of ${\bq}_S^m$.
Then
$
\Delta = {\bz}_S\va_1 \oplus \dots \oplus {\bz}_S\va_r
$
for some $\va_1,\dots,\va_r \in {\bq}_S^m$
such that
$$\va_1^{(v)},\dots,\va_r^{(v)}\text{ are linearly
independent over }{\bq}_v\text{  for any }v\in S\,.\tag 7.2$$
Furthermore, there
exists $g \in \GL(m,{\bq}_S)$ such that $\Delta$ is contained in
$g {\bz}_S^m$.
\endproclaim

\demo{Proof} The proposition is trivial if $S = \infty$. So assume
that $S \varsupsetneq \infty$ and denote by $\Delta_0$ the
intersection of
$\Delta$ with $
{\bq}_{S,f}^m$. Let $\pi: {\bq}_{S,f}^m \rightarrow
{\bq}_{\infty}^m$ be the natural projection. Since $\ker (\pi) =
{\bz}_{S,f}^m$ is compact and $ {\bz}_{S,f}^m$ does not contain
nontrivial discrete subgroups, $\pi(\Delta)$ is a free  abelian
group of rank $r \leq m$ and $\pi(\Delta)$ is isomorphic to
$\Delta$. If $\vx$ is any element in $\Delta$, then there exists
$\xi\in{\bz}_S^*$ (the group of $S$-adic units)
  such that $\xi \vx \in \Delta_0$. This
implies that $\Delta$ is a free ${\bz}_S$-module of rank $r$.

Let $\Delta = {\bz}_S\va_1 \oplus \dots \oplus {\bz}_S\va_r$. Suppose that
$\va_1^{(v)} = 0$ for some $v \in S$,   and let $\{\xi_i\}$ be a
sequence of $S$-adic units such that $\lim_{i \rightarrow \infty}
|\xi_i|_w = 0$ for all $w \neq v$. Then $\lim_{i \rightarrow
\infty} \xi_i\va_1 = 0$ which contradicts the discreteness of
$\bz_S \va_1$. Therefore
  $\va_1^{(v)} \neq 0$ for all $v \in S$, which proves (7.2)
   for $r = 1$. To complete the proof we use
  induction on $r$. Assume that $r > 1$ and (7.2)
  is true for free modules of rank $ < r$. Shifting $\Delta$
  by an appropriate automorphism from $\GL(m,\bq_S)$, without loss of
  generality we may and
  will assume that $\va_1 = \ve_{1}$,
the first vector of
  the standard basis of $\bq_S^m$. Let $\varphi: \bq_S^m \rightarrow \bq_S^m
/\bq_S\,\ve_1$
  be the natural homomorphism. Since $\bz_S\,\ve_1$ is a cocompact lattice in
$\bq_S\,\ve_1$,
  we get that $\varphi(\Delta)$ is discrete in $\varphi(\bq_S^m) \cong
  \bq_S^{m-1}$. By the induction hypothesis
  $\varphi(\va_2)^{(v)},\dots,\varphi(\va_r)^{(v)}$ are linearly
  independent over $\bq_v$ for all $v$, which  completes the proof of
(7.2). It remains to
   observe that the last part of the proposition immediately follows from (7.2).
  \qed\enddemo

Note that it follows from Proposition 7.2 that the rank of $\Delta$ as a
free module is equal to $\rk(\Delta)$ as defined  in (6.1) with $K
= \bq$.

\subhead{7.3}\endsubhead Our next goal is to define the normalized
Haar measure  on free ${\bq}_S$-submodules  of ${\bq}_S^m$.
      Let $L$ be a free
${\bq}_S$-module of rank $r$ generated by $\va_1,\dots,\va_r\in
{\bq}_S^m$. Then one can write $L =\prod_{v\in S}L_v$, where $L_v
=
{\bq}_v \va_1^{(v)}\oplus\dots\oplus {\bq}_v\va_r^{(v)}$. 
Let us fix a basis $(\vb_1,\dots,\vb_r)$ of $L$ with the following
properties:
$$
\text{ if }v = \infty\text{, \ \  then }
\text{$(\vb_1^{(v)},\dots,\vb_r^{(v)})$ is an orthonormal basis
of } L_v\,;\tag 7.3a
$$
$$   \text{if $v\in S_f$, then } L_v \cap {\bz}_v^m =
{\bz}_v \vb_1^{(v)}+\dots+ {\bz}_v \vb_r^{(v)}\,.\tag 7.3b
$$
Then
consider the ${\bq}_S$-linear map 
      sending the standard
basis $(\ve_1,\dots,\ve_r)$ of ${\bq}_S^r$ to
$(\vb_1,\dots,\vb_r)$, and
          define the {\sl volume\/}  on $L$ as the pushforward of
$\lambda_S$ (the normalized Haar measure on ${\bq}_S^r$) by
this map.
When it does not lead to confusion, we
will denote this measure on $L$ by $\lambda_S$ as well.

Note that the existence of $\vb_1^{(v)},\dots, \vb_r^{(v)}$ 
in the case (7.3b) easily
follows from the fact that  ${\bz}_v$ is a principal ideal domain.
Note also that the above definition does not depend on the choice
of the basis $(\vb_1,\dots,\vb_r)$ satisfying (7.3ab), because if
$(\vb'_1,\dots,\vb'_r)$ is another such basis and $h\in
\GL(r,\bq_S)$ represents the isomorphism of $L\cong \bq_S^r$
sending the first basis to the second one, then
$c\big(\det(h)\big)= 1$, which implies that $h$ is
measure-preserving.

\medskip

For any $r = 1,\dots,m$ we will also consider  the $r$-th
exterior power
$$\tsize \bigwedge^r {\bq}_S^m\simeq \bigoplus_{v \in S} \bigwedge^r {\bq}_v^m
$$ of 
${\bq}_S^m$, which  is a free ${\bq}_S$-module with the standard
basis
$$
\{\ve_{i_1}\wedge \dots\wedge \ve_{i_r}
\mid 1 \leq i_1 <  \dots < i_r \leq m \}\,,
$$
where $\ve_1,\dots,\ve_m$ is the standard basis of ${\bq}_S^m$. We
will keep the notation $\lambda_v$, $\lambda_S$, $\|\cdot\|_v$
and $c(\cdot)$ to denote the measures, the norms and the content
on the exterior powers $\bigwedge^r {\bq}_v^m$ and $\bigwedge^r
{\bq}_S^m$, respectively.


\subhead{7.4}\endsubhead Recall (cf.\ \cite{W}) that ${\bz}_S$ is
a lattice in ${\bq}_S$ with covolume $1$,  that is, it is discrete
in ${\bq}_S$ and the Haar measure (induced by $\lambda_S$) of the
quotient space is equal to one. Likewise, ${\bz}_S^m$ is a lattice
in ${\bq}_S^m$ with covolume $1$. It follows from Proposition 7.2
that the set of discrete ${\bz}_S$-submodules  of ${\bq}_S^m$ can
be identified with $\frak M(\bq_S,{\bz}_S,m)$ as defined in (6.2).
Furthermore, it also follows that any $\Delta\in \frak M(\bq_S,{\bz}_S,m)$ is a
lattice in $\bq_S\Delta$.
The following lemma shows how one can explicitly compute
covolumes:

\proclaim{Lemma} Let $\Delta = {\bz}_S \va_1 \oplus \dots \oplus
{\bz}_S \va_r\in \frak M(\bq_S,{\bz}_S,m)$,
where
$\va_1,\dots,\va_r \in
\bq_S^m$. Then the covolume  $\cov(\Delta)$  of $\Delta$ in
$\bq_S\Delta$ with respect to the volume on $L = \bq_S\Delta$
normalized as in \S 7.3 is equal to
$$
\cov(\Delta) =
c ( \va_1\wedge\dots\wedge
\va_r)\,.
$$
\endproclaim

\demo{Proof}
Put $L =
\bq_S\Delta$,  define a basis $(\vb_1,\dots,\vb_r)$ of $L$ as in
(7.3ab), and then complete
it to a basis $(\vb_1,\dots,\vb_m)$ of the whole space $\bq_S^m$.
Also let $h_v \in \GL(m,\bq_v)$  be such that $h_v\vb_i^{(v)} =
\ve_i^{(v)}$ for all $1 \leq i \leq m$,
         where $\ve_1^{(v)},\dots,\ve_m^{(v)}$ is the standard basis
of $\bq_v^m$.
It follows from the definition of the measure on $L$ that the map
$h:L \rightarrow h(L)$, where $h = (h_v)\in \GL(m,\bq_S)$,  is
measure preserving. Since the map $\bigwedge^r h: \bigwedge^r
\bq_S^m \rightarrow \bigwedge^r \bq_S^m$ preserves the content $c$
on $\bigwedge^r \bq_S^m$, we may reduce the problem to the case
$\vb_i = \ve_i$ for all $1 \leq i \leq m$. Let $\varphi \in
\GL(r,\bq_S)$ be such that  $\varphi(\ve_i) = \va_i$ for all $1
\leq i \leq r$. Then since $\va_1\wedge\dots\wedge \va_r = \det
(\varphi) \ve_1\wedge\dots\wedge \ve_r$, we get
$$
\cov(\Delta)=\cov\big(\varphi ({\bz}_S^r)\big)=c(\det
\varphi)\cov({\bz}_S^r)=  c ( \va_1\wedge\dots\wedge \va_r)\,.\qed
$$ \enddemo

\subhead{7.5}\endsubhead Lemma 7.4 immediately implies

\proclaim{Corollary} For every
$\Delta\in \frak M(\bq_S,{\bz}_S,m)$,
         the function
$ \GL(m,{\bq}_S)\to \br_+$, $g\mapsto \cov(g\Delta)$, is
continuous.
\endproclaim

\proclaim{7.6.\ Corollary} If $\Delta,\ \Delta'\in \frak
M(\bq_S,{\bz}_S,m)$ are
such that $\bq_S \Delta  \,\cap\, \bq_S \Delta' = \{0\}$, then
$$
\cov(\Delta + \Delta') \leq \cov(\Delta) \cov(\Delta').
$$
\endproclaim

\demo{Proof} Using Proposition 7.2, write $\Delta = \bz_S \va_1
\oplus \dots \oplus \bz_S \va_r$ and $\Delta' = \bz_S \vb_1 \oplus
\dots \oplus \bz_S \vb_r$. Since
$$
\bq_S\Delta + \bq_S\Delta' = \bq_S\Delta \oplus \bq_S\Delta',
$$
in view of Lemma 7.4, it is enough to prove that
$$
c(\va_1 \wedge \dots \wedge \va_r \wedge \vb_1 \wedge \dots \wedge \vb_s)
\leq c(\va_1 \wedge \dots \wedge \va_r) c(\vb_1 \wedge \dots \wedge
\vb_s),
$$
which is  easy to verify using the definition of content and the basic
properties of the exterior product. \qed\enddemo

\subhead{7.7}\endsubhead In the remaining part of the section we
investigate metric properties of discrete submodules of $\bq_S^m$.
Let us   state the following $S$-arithmetic version of the
classical  Minkowski's Lemma:

\proclaim{Lemma
} Let
       $\Delta\in \frak M(\bq_S,{\bz}_S,m)$  be
of  rank $r$,
and let $B$ be a closed ball in $ \bq_S\Delta$ (with respect to
the norm $\|\cdot\|$) centered at $0$ such that
      $\lambda_S(B) \geq 2^{r}\cov(\Delta)$.
Then $ \Delta \cap B \neq
\{0\} $.
\endproclaim

\demo{Proof}
Since the volume $\lambda_S$ on $L =  \bq_S\Delta$  was defined by
identifying $L$ with $ \bq_S^r$ via the basis 
(7.3ab), without loss of generality we can assume that $\Delta$ is
a lattice in $ \bq_S^r$.
Then we can write $B = B_{\infty} \times B_{f}$, where $B_{\infty}
\df B  \,\cap\, \bq_{\infty}^r$ and $B_{f} \df B  \,\cap\, \left(\prod_{v\in
S_f}\bq_{v}\right)^r$. Note that
$$
\lambda_S(\tfrac{1}{2}B_{\infty} \times B_{f}) =
\frac{1}{2^{r}}\lambda_S(B) \geq \cov(\Delta).
$$
Since $\frac{1}{2}B_{\infty} \times B_{f}$ is closed, the above
implies that there exist $\vx, \vy \in \frac{1}{2}B_{\infty} \times
B_{f}$, $\vx \neq \vy$, such that $\vx-\vy \in \Delta$.
This finishes the proof, since clearly $\vx-\vy \in B$ as well.
\qed\enddemo



\subhead{7.8}\endsubhead We will also need the following result:

\proclaim{Lemma} There exists a constant $A > 1$ depending only on
$S$ such that if $\vx \in \bq_S^m$ and $c(\vx) \neq 0$,
then there exists $\xi\in{\bz}_S^*$ 
such that
$$
A^{-1} c(\vx)^{1/\ell} \leq \|\xi \vx \| \leq A c(\vx)^{1/\ell},
$$
where $\ell$ is the cardinality of $S$.
\endproclaim

\demo{Proof} Let
$$
H = \{ (a_1,\dots,a_\ell) \in \br^\ell_+ \mid a_1\cdots a_\ell = 1\}.
$$
Write $S = \{v_1,\dots,v_\ell\}$;
it is easy to see that the group
$$
\big\{(|\xi|_{v_1},\dots,|\xi|_{v_\ell})\bigm|\xi\in{\bz}_S^*\big\}
$$
is a cocompact lattice in the multiplicative group $H$. Therefore
there exists a constant $A
>  1$ such that for any $(a_i) \in H$ one can find an $S$-adic unit $\xi$ with
$$
A^{-1} \leq |\xi|_{v_i} a_i \leq A\tag 7.4
$$
for all $i$.

Let $\vx = (\vx^{(v_i)}) \in \bq_S^m$ and $c(\vx) \neq 0$. Note
that the vector
$\left(\frac{\|\vx\|_{v_i}}{c(\vx)^{1/\ell}}\right)$ is in $H$,
and
one has
$$
\|\xi \vx^{(v_i)} \|_{v_i} = |\xi|_{v_i} \|\vx^{(v_i)}\|_{v_i},
$$
and
$$
c(\vx) = c(\xi \vx)
$$
for all $\xi \in {\bz}_S^* $. This, in view of (7.4), implies  the
claim. \qed\enddemo

\subhead{7.9}\endsubhead Let us denote by $\Omega_{S,m}$ the space
of all lattices in $\bq_S^m$. It follows from Lemmas 7.2 and 7.4
that it can also be defined as
$$\Omega_{S,m} = \{ g{\bz}_S^m \mid g \in \GL(m,\bq_S) \} \cong
\GL(m,\bq_S)/\GL(m,\bz_S)\,.
$$

\proclaim{Corollary}
For any  $\Lambda  \in \Omega_{S,m}$ and any
$\rho > 0$,  the number
of   submodules of
$\Lambda$ with covolume $ \le \rho$ is finite.
\endproclaim

\demo{Proof}
If $\Delta$ is a
${\bz}_S$-submodule with $\rk(\Delta) = r$, then $\bigwedge^r
\Delta$ is a
${\bz}_S$-submodule of rank $1$ of the lattice $\bigwedge^r
\Lambda $ in $\bigwedge^r \bq_S^m$, and, in view of Lemma 7.4,
$$
\cov(\Delta) =
\cov(\tsize \bigwedge^r \Delta).
$$
Therefore it is enough to prove that the number of rank-one 
${\bz}_S$-submodules of $\Lambda$ is finite. If $\Delta = {\bz}_S
\va$ is such a submodule, then in view of Lemma 7.8 the generator
$\va$ can be chosen in such a way that
$$
\|\va\| \leq A\rho^{1/\ell}.
$$
Since $\Lambda $ is discrete in $\bq_S^m$, the set of all $\va \in
\Lambda$ satisfying the above inequality is finite, which proves
the corollary. \qed\enddemo

\subhead{7.10}\endsubhead Because of lack of an appropriate
reference we will prove an
$S$-adic version of Mahler's Compactness Criterion. 
Consider the
group
$$
\GL^1(m,\bq_S) \df \big\{g \in \GL(m,\bq_S)\bigm|
c\big(\det(g)\big)= 1\big\}\,,
$$
consisting of $\lambda_S$-preserving linear automorphisms of
$\bq_S^m$, and let
$$\Omega^1_{S,m} \df \big\{ \Lambda \in \Omega_{S,m} \mid \cov(\Lambda)
=1
\big\}
\,.
$$
Note that ${\bz}_S^m$ is an element of $\Omega^1_{S,m}$, and its
stabilizer in $\GL(m,\bq_S)$ coincides with $\GL(m,{\bz}_S)$ which
is understood to be diagonally imbedded in $\GL(m,\bq_S)$. Since
$c(\xi) = 1$ for any $\xi\in {\bz}_S^*$,
it follows  that
$\GL(m,{\bz}_S)$ is contained in $\GL^1(m,\bq_S)$.
Thus $\Omega^1_{S,m}$ is naturally identified
with
the \hs\ $\GL^1(m,\bq_S)/\GL(m,{\bz}_S)$.
        Since $\SL(m,{\bz}_S)$ is a lattice
in $\SL(m,\bq_S)$ and ${\bz}_S^*$ is a lattice in $\GL(1,\bq_S)$,
$\GL(m,{\bz}_S)$ is a lattice in $\GL^1(m,\bq_S)$.

\medskip

Let us say that a set $Q$ of lattices in  ${\Cal R}^m$, where ${\Cal R}$ is a
topological ring,  is {\sl separated from
$0$
\/} if  there exists a nonempty
neighborhood $B$ of $0$ in  ${\Cal R}^m$
such that $\Lambda  \,\cap\, B = \{0\}$ for all $\Lambda \in
Q$.

\proclaim{Theorem {\rm (Mahler's Compactness Criterion)}} A subset
$Q \subset \Omega^1_{S,m}$ is bounded if and only if it is
separated from $0$.
\endproclaim

\demo{Proof} The implication ($\Longrightarrow$) is trivial. In order
to prove the converse, note that
$$\GL^1(m,\bq_S) =  \bq_S^{1} \ltimes  \SL(m,\bq_S)\quad \text{ and
         }\quad\GL(m,{\bz}_S) =  {\bz}_S^{*} \ltimes  \SL(m,{\bz}_S)\,,$$
where $ \bq_S^{1} = \{ x \in
          \bq_S \mid c(x) = 1 \}$. Since ${\bz}_S^{*}$ is a cocompact lattice
in $ \bq_S^{1}$, it is enough to
         prove the theorem with $\Omega^1_{S,m}$ replaced by  $\SL(m,\bq_S)/
\SL(m,{\bz}_S)$, i.e.\ with the set $$\{ g {\bz}_S^m \mid g \in
\SL(m,{\bz}_S)\}\,.$$

It follows from the strong approximation theorem for classical groups
         \cite{Kn}
that
$$
         \SL(m,\bq_S) = \SL(m, \bq_{S,f}) \SL(m,{\bz}_S)\,.
$$
Thus  every $g \in \SL(m,\bq_S)$ can be represented as $g =
g_{f}g_l$,
where $g_{f} \in  \SL(m,\bq_{S,f})$ 
and $g_{l} \in \SL(m,{\bz}_S)$. One has
$$
g{\bz}_S^{m}  \,\cap\,  \bq_{S,f}^{m} = g_{f}\left({\bz}_S^{m}  \,\cap\,
\bq_{S,f}^{m}\right) = g_{f}{\bz}^m\,.\tag 7.5
$$
      Let $\tilde Q \subset \SL(m,\bq_S)$ be such that the set of lattices
$\{g {\bz}_S^{m} \mid  g \in\tilde Q\}$ is separated from $0$.
It follows from (7.5) that
$$\{g_f {\bz}^m \mid  g \in \tilde Q\}\text{
is separated from $0\,$ in  }\bq_{S,f}^{m}\,. \tag 7.6
$$
Note that $\SL(m,\bq_{S,f}) = \SL(m, \br)\times
\SL(m,{\bz}_{S,f})$, and, therefore, every $g_f\in \SL(m,\bq_{S,f}) $ can
be written as $g_f = g_{\infty} g_c$, where $g_{\infty} \in \SL(m,
\br)$  and $g_{c}$ belongs to the compact group $
\SL(m,{\bz}_{S,f})$. It follows from (7.6) that $\{g_\infty {\bz}^m
\mid  g \in \tilde Q\}$ is separated from $0$.
      This reduces the proof to
the case $S = \infty$, that is, to the original Mahler's Criterion
(see \cite{R, Corollary
10.9}).
\qed\enddemo

In particular, it follows from
Lemma 7.8 and Theorem 7.10 that for all positive $\vre$, the sets
$Q_\vre$ defined as in (0.2) are compact.

\heading{8. $S$-arithmetic quantitative nondivergence}
\endheading

\subhead{8.1}\endsubhead In this section we apply Theorem 6.3 to
the triple $({\Cal D}, K,{\Cal R})$ of a particular type. Namely,
as in \S 7, we let $K = \bq$, choose a finite set $S$ of
valuations $|\cdot|_v$ of $\bq$ containing  the Archimedean one,
and for the rest of this section take  ${\Cal D} = {\bz}_S$ and
${\Cal R} = \bq_S$. Then for any $m\in\bn$, the set $\frak
M(\bq_S,{\bz}_S,m)$, as defined in (6.2), is equal to the set of
all submodules of all lattices $\Lambda\in\Omega_{S,m}$, where
$\Omega_{S,m}$ is as defined in (7.5). Let us now state the
following

         \proclaim{Lemma} The function  $\nu:\frak M(\bq_S,{\bz}_S,m)\to\br_+$
given by $\,\nu(\Delta) = \cov(\Delta)$, with $\,\cov(\cdot)$ as
in \S7.4, is norm-like, with $C_\nu = 1$.
\endproclaim

\demo{Proof} Property (N1) is straightforward since
$\Delta'\subset \Delta$ and $\bq_S\Delta' = \bq_S\Delta$ implies
that $\Delta$ is a subgroup of $\Delta'$ of finite index and
$\cov(\Delta') = [\Delta:\Delta']\cov(\Delta)$.  Property (N2)
with $C_\nu = 1$ follows from Corollary 7.6 with $\Delta' =
{\bz}_S\gamma$. Finally, (N3) has already been mentioned as
Corollary 7.5. \qed\enddemo


\subhead{8.2}\endsubhead Now define a function $\delta:\Omega_{S,m}\to\br_+$ by
$$
\delta(\Lambda) \df \min \big\{c(\vx)\bigm| \vx \in\Lambda\nz
\big\}\,.
$$
Note that the minimum
is well defined
due to Lemma 7.8 and every $\Lambda\in\Omega_{S,m}$ being discrete
in  $\bq_S^m$. We will use the following

         \proclaim{Lemma} There exists a   constant $A > 0$ depending only on
$S$ and $m$ such that the following holds: for $\rho > 0$ and
$\Lambda\in\Omega_{S,m}$ suppose there exists a submodule $\Delta$ of $\Lambda$
with $\cov(\Delta) \le \rho$; then $
\delta(\Lambda) \le A\rho^{1/m}$.
\endproclaim

\demo{Proof} Take $\vre > 0$ and let $B$ be a ball in
$\bq_S\Delta$  centered at $0$ of radius $\vre$ (with respect to
the norm $\|\cdot\|$ introduced in \S 7.1). Then one has
$\lambda_S(B) \le \const\cdot\vre^{r\ell}$, 
  where $r = \rk(\Delta)$,
$\ell$ is the cardinality of $S$, and the
constant depends only on  $S$ and $m$.  By Lemma 7.7, $\Delta$ has a
nontrivial intersection with $B$ whenever $\const\cdot\vre^{r\ell}
\ge 2^{r}\rho$. This shows how one can choose $A$ such that
$\Delta$ (and hence $\Lambda$) is guaranteed to contain a nonzero
vector $\vx$ with $\|\vx\|\le A^{1/\ell}\rho^{1/r\ell}$, which
clearly implies $c(\vx) \le A\rho^{1/r} \le A\rho^{1/m}$.
\qed\enddemo

\subhead{8.3}\endsubhead As in \S 6, let us use the notation
$\frak P({\bz}_S,m)$ for the set of all nonzero primitive
submodules  of ${\bz}_S^m$.


         \proclaim{Theorem} Let $X$ be a 
\be\ metric space,
$\mu$ a  uniformly Federer measure  on $X$, and let
$S$
         be as above.  For  $m\in \bn$, let  a
ball $B = B(x_0,r_0)\subset X$ and a continuous map $h:\tilde B
\to \GL(m,\bq_S)$ be given, where $\tilde B$ stands for
$B(x_0,3^mr_0)$.
Now suppose that for some
$C,\alpha > 0$ and $0 < \rho  < 1$  one has

{\rm(i)} for every $\,\Delta\in \frak P({\bz}_S,m)$, the function
$\cov\big(h(\cdot)\Delta\big)$ is \cag\ on $\tilde B$  with
respect to $\mu$;

{\rm(ii)}  for every $\,\Delta\in \frak P({\bz}_S,m)$,
$\|\cov\big(h(\cdot)\Delta\big)\|_{\mu,B} \ge \rho$.

         Then 
         for any  positive $ \vre \le
\rho$ one has
$$
\mu\left(\big\{x\in B\bigm| \delta \big(h(x){\bz}_S^m\big) < \vre
\big\}\right)\le mC \big(N_{\sssize X}D_{\mu}^2\big)^m
\left(\frac\vre \rho \right)^\alpha \mu(B)\,.\tag 8.1
$$
\endproclaim

Note   that
\cite{KLW, Theorem 4.3} is a special case (corresponding to 
$S =
\{\infty\}$) of the above theorem.

\demo{Proof} To apply Theorem 6.3, one uses Lemma 8.1 which
guarantees the norm-like property of $\cov(\cdot)$, and Corollary
7.9 which implies condition (iii) of Theorem 6.3. To derive (8.1)
from (6.3) it remains to observe that $\delta
\big(h(x){\bz}_S^m\big) < \vre$ amounts to the
existence of a vector 
$\vx\in h(x){\bz}_S^m\nz$ with
        $$\cov({\bz}_S\vx) = c(\vx) <
        \vre\,.\quad\qed
$$ 
\enddemo

\subhead{8.4}\endsubhead In order to interpret the above result,
let us assume, as it will be the case in many applications, that
the function $h$ takes values in the group $ \GL^1(m,\bq_S)$. Then
  $ h(x){\bz}_S^m$ belongs to
$\Omega^1_{S,m}$ for any $x$, and the inequality $\delta
\big(h(x){\bz}_S^m\big) <\vre$
can be equivalently written as $h(x){\bz}_S^m\notin Q_\vre$. 
defined in (7.7).
      This way,
Theorem 8.3 estimates, in terms of $\vre$,
         the relative measure of points $x\in B$ which are mapped,
by $x\mapsto h(x){\bz}_S^m$, to the complement of $Q_\vre$
in $\Omega^1_{S,m}$.

\medskip
As an application, let us take $X$, $\mu$ and $h$ of a special
form. Namely, for every $v\in S$ choose $d_v \in\bn$, and consider
$X = \prod_{v\in S}\bq_v^{d_v}$, $\mu = \lambda$ as defined  in
(0.3),   and a map $h = (h_v)_{v\in S}: X\to \GL^{1}(m,\bq_S)$,
where
        each   $h_v$ is a map from $\bq_v^{d_v}$ to $\GL(m,\bq_v)$. We
        say
that $h$ is {\sl polynomial\/} (or {\sl regular\/}) if for every
$v$ all matrix coefficients of $h_v$ and
its inverse are polynomials (equivalently, if every $h_v$ is the
restriction of a regular map of algebraic varieties
$\bar{\bq}_v^{d_v} \rightarrow \GL(m,\bar{\bq}_v)$, where
$\bar{\bq}_v$ is the algebraic closure of $\bq_v$).

          \proclaim{Theorem}
Let $X$ and $h$ be defined as above.  Then there exists $\alpha >
0$ (depending only on $m$, $d_v$ and  the degrees of the maps)
such that for every compact set $L\subset \Omega^1_{S,m}$ one can find 
positive 
           $C_0$ and $\tau$ (depending only
on $m$, $d_v$, the degrees of the maps $h_v$, and $L$) such that $Q_\tau
\supset L$,
and the following property holds: for any
positive $\vre$ and any ball $B\subset X$ one has
$$
\lambda\left(\big\{x\in B\bigm| h(x){\bz}_S^m \notin Q_\vre
\big\}\right) \le C_0\vre^\alpha \lambda(B)\tag 8.2
$$
whenever
$
h(B){\bz}_S^m  \,\cap\, Q_\tau\ne\varnothing
$. Furthermore, if $ h(X){\bz}_S^m  \,\cap\, Q_\tau = \varnothing$,
then there exists a proper $\,\Delta\in \frak P({\bz}_S,m)$ such
that $h(x)(\bq_S\Delta) = h(0)(\bq_S\Delta)$ for all $x \in X$.
\endproclaim

\demo{Proof} It was mentioned in Example 1.6 that $X$ is
Besicovitch, and in Example 1.7 that $\lambda$ is uniformly
Federer. Using the exterior power representation of Lemma 7.4, one
can easily show that for every $\Delta\in \frak P({\bz}_S,m)$ the
function
$\cov\big(h(\cdot)\Delta\big)$ 
has the form $\prod_{v\in S}\|f_v \|_v$,  where each $f_v$ is a
polynomial map from $\bq_v^{d_v}$ to another vector space over
$\bq_v$. Every such function is $(C',\alpha')$-good with uniform
$C'$ and $\alpha'$ due to Lemma 2.4, Lemma 2.1(cd)
          and Corollary 2.3. Thus
condition (i) of Theorem 8.3 is satisfied with some $C,\alpha  > 0$.

Now let us take $\tilde C\df \max\left(mC
\big(N_{\sssize X}D_{\lambda}^2\big)^m,1\right)$.
%
It follows from Lemma 8.2 and Theorem 7.10 that there
exists  $\tau = \tau(L) > 0$ such that for any submodule $\Delta$
of any $\Lambda \in
L$ one has
$\cov\big(\Delta\big) \ge
           \tau$. Without loss of generality we can assume that $\tau <
\frac1{(2\tilde C)^{1/\alpha}}$. 
Note that $L$ is contained in $Q_\tau$, since by definition of
$\tau$ one has $c(\vx) = \cov({\bz}_S\vx) \ge \tau$ for any
nonzero element $\vx$ of any $\Lambda \in L$.

If $B$ is such that $ h(B){\bz}_S^m  \,\cap\, Q_\tau\ne\varnothing$, it
follows from Lemma 8.2 that
       condition (ii) of Theorem 8.3 is satisfied with $(\tau/A)^m$ in
place of $\rho$ (where $A$ is as in Lemma 8.2).  Therefore one has
$$
\lambda\left(\big\{x\in B\bigm| h(x){\bz}_S^m \notin Q_\vre
\big\}\right) \le  \tilde C \left(\frac{\vre
A^m}{\tau^m}\right)^\alpha
\lambda(B)
$$
for all $\varepsilon \leq \tau^m\,/ A^m$. Replacing, if necessary,
$\tilde C$ by a larger number $C$, we conclude that (8.2) is valid
for all positive $\varepsilon$.

Now assume that $h(x){\bz}_S^m \notin Q_\tau$ for all $x\in X$.
Take $\rho  \df (2\tilde C)^{1/\alpha}\tau < 1$, write $X =
\bigcup_{i = 1}^\infty B_i$ where $B_i$ are balls centered at $0$
with  $B_i\subset B_{i+1}$ for all $i$, and consider
$$
\frak P_i \df\big\{\Delta\in \frak P({\bz}_S,m) \bigm|
\|\cov\big(h(\cdot)\Delta\big)\|_{B_i} < \rho\big\}\,.
$$
Then clearly $ \frak P_{i+1} \subset \frak P_i$ for all $i$, and
all these sets are finite due to Corollary 7.9. We claim that
$\bigcap_{i = 1}^\infty \frak P_i$ must be nonempty. Indeed,
otherwise one obtains a nonempty ball $B$ such that
$\|\cov\big(h(\cdot)\Delta\big)\|_{B} \ge \rho$ for every
$\Delta\in \frak P({\bz}_S,m)$.  Thus Theorem 8.3 can be applied,
and one can conclude that
$$
\lambda(B) = \lambda\left(\big\{x\in B\bigm| h(x){\bz}_S^m \notin
Q_\tau \big\}\right)\le \tilde C \left(\frac{\tau} \rho
\right)^\alpha \lambda(B) = \frac12\lambda(B)\,,
$$
a contradiction.

Consequently, there exists  a proper $\Delta\in \frak
P({\bz}_S,m)$ such that $\cov\big(h(x)\Delta\big) < \rho$ for all
$x\in X$. It follows that each of the polynomials  $f_v$ in the
aforementioned representation  for
      $\cov\big(h(\cdot)\Delta\big)$ is bounded.
Therefore  $f_v \equiv \const$ for each $v$,  which implies that
$h(x)(\bq_S\Delta)$ does not depend on $x$. \qed\enddemo

We note that to derive Theorem 0.2 from the above theorem, one
needs to take $L = \{\Lambda\}$ and $h$ of the form $x\mapsto
h(x)g$, where $g\in \GL^{1}(m,\bq_S)$ is such that  $\Lambda =
g{\bz}_S^m$, and observe  that $ h(B){\bz}_S^m  \,\cap\, Q_\tau$ is
nonempty whenever $B$ contains $0$.

\heading{9.  Invariant locally finite measures for actions of
unipotent groups on homogeneous spaces}
\endheading

\subhead{9.1}\endsubhead Theorem 8.4 implies results closely
related to the structure of
       orbits and  invariant measures on $S$-adic homogeneous
spaces under the action of subgroups generated by unipotent
elements (see \cite{Rt4}, \cite{MT2}  and \cite{To}).

Let us recall some definitions from \cite{MT2}. As in the previous
sections, $S$ is  a finite set of normalized valuations of $\bq$
containing the archimedean one and $\bq_S$ is the direct product
of all $\bq_v$, $v \in S$. By a {\sl $\bq_{S}$-algebraic group\/}
$\bold{G}$ we mean a (formal) direct product $\prod_{v \in S}
\bold{G}_v$ of $\bq_{v}$-algebraic groups $\bold{G}_v$. The group
$\prod_{v \in S} \bold{G}_v(\bq_{v})$ will be denoted by
        $\bold{G}(\bq_{S})$ and called {\sl the group of
$\bq_{S}$-rational points\/} of $\bold{G}$. We will also use the
simpler notations
        $G$ for $\bold{G}(\bq_{S})$ and $G_v$ for $\bold{G}_v(\bq_{v})$.
If $\bold{H}$ is another
        $\bq_{S}$-algebraic group then a homomorphism $\varphi:
\bold{G} \to \bold{H}$ is called $\bq_{S}$-{\sl homomorphism\/} if
$\varphi$ is a product of $\bq_{v}$-homomorphisms of algebraic
groups $\varphi_v: \bold{G}_v \to \bold{H}_v$, $v \in S$. We
preserve the same terminology for the restriction map $\varphi: G
\to H$. By the {\sl Zariski topology\/}  on $\bold{G}$
(respectively, $G$) we mean the formal product of the Zariski
topologies on $\bold{G}_v$ (respectively, $G_v$). By a
$\bq_{S}$-{\sl algebraic subgroup\/} of $G$ (or simply an
algebraic subgroup of $G$) we mean a Zariski closed subgroup of
$G$. An element $g = (g_v)$ in $G$ is {\sl unipotent\/} if each
component $g_v \in G_v$ is unipotent. A subgroup of $G$
    consisting of unipotent elements is called {\sl
unipotent\/}.

Up to the
        end of this section we will denote by $\Gamma$ a lattice in $G$.
        Any subgroup
of $G$ acts on the homogeneous space $G/\Gamma$ by left
        translations.

Let $U$ be a unipotent algebraic subgroup of $G$. Then 
$U = \prod_{v \in S} U_v$, where $U_v$ are unipotent algebraic
subgroups of $G_v$. Given $v \in S$, we denote by $\exp_v:
\Lie(U_v) \rightarrow U_v$ the exponential map and by $\log_v =
\exp_v^{-1}$ the logarithmic map. Also we denote by $\Lie(U)$ the
direct product of the Lie algebras $\Lie(U_v)$ of $U_v$, $v \in
S$, and  by $\exp: \Lie(U)\rightarrow U$ the direct product of the
maps $\exp_v$, $v \in S$. By a (rational) {\sl parametrization \/}
of $U$ we mean a product $\phi = (\phi_v)_{v \in S}$ of surjective
maps $\phi_v: \bq_{v}^{d_v} \rightarrow U_v$, $v \in S$, such that
for every $v \in S $ the map $\log_v \circ \phi_v : \bq_{v}^{d_v}
\rightarrow \Lie(U_v)$   is polynomial and $\phi_v(0) = e$. If
$d_v$ is the degree of $\log \circ \phi_v$ then $d = \max\{d_v | v
\in S\}$ is called the {\sl degree\/} of the parametrization
$\phi$. Clearly $\exp$ is a parametrization of $U$
which we call trivial. 
We will denote by $\lambda$ the Haar measure on $X = \prod_{v\in
S}\bq_{v}^{d_v}$. We fix a metric on $X$ and by a ball in $X$ we
mean a ball with respect to this metric.

The following theorem generalizes earlier results, which  for
one-parameter real groups were proved in \cite{
     D1--2}, and for
one-parametric ultrametric groups were  announced, with indications of
the  proof,
in \cite{MT2, Theorem 11.4} and \cite{Rt3, Theorem
9.1}.

\proclaim{Theorem} Let $G$  and $\Gamma$  be as above, let  $d$ be
a natural number, and let $L$ be a compact  subset of $G/\Gamma$. Then
$L$ is contained in a  compact $L_0$ with the following property:
given $\beta
>      0$ there exists a compact $M_\beta \subset G/\Gamma$  such that
for any $y \in G/\Gamma$ and any parametrization $\phi: X\to U$ of
degree $\leq d$ of a unipotent algebraic subgroup $U$ of $G$ one
of the following holds:

{\rm(i)} If $Uy  \,\cap\, L_0 \neq \varnothing$ and $B$ is a ball in $X$
with $\phi(B)y  \,\cap\, L_0 \neq \varnothing$ then
$$
\frac{\lambda \left(\big\{x\in B\bigm| \phi(x)y \in M_\beta
\big\}\right)}{\lambda(B)} \geq 1 -\beta  \,.\tag 9.1
$$
In particular, {\rm (9.1)} is satisfied if $y \in L$ and $B$ contains the
origin.

{\rm(ii)} If $Uy  \,\cap\, L_0 = \varnothing$, then there exists a
proper closed subgroup $H$ of $G$ containing $U$ such that the orbit
$Hy$ is closed and carries a finite $H$-invariant Borel measure.

\endproclaim

\subhead{9.2}\endsubhead Before proving the theorem we will
establish the following result which is well known in the real
case.

\proclaim{Proposition} Let $G$  and $\Gamma$  be as in Theorem
9.1, and  let $R(G)$ be the solvable radical of $G$ (i.e.\ $R(G)$
is the direct product of the solvable radicals $R(G_v)$ of $G_v$,
$v \in S$). Then $R(G)  \,\cap\, \Gamma$ is a cocompact lattice in
$R(G)$.

\endproclaim

\demo{Proof} Let $\widetilde{\Gamma}$ be the Zariski closure  of
$\Gamma$ in $G$. In view of the Borel Density Theorem for
$\bq_{S}$-algebraic groups
    (see  \cite{MT2, Lemma 3.1}) $\widetilde{\Gamma}$ contains all
unipotent algebraic
    subgroups and all $S$-split tori of $G$. Therefore
$G/\widetilde{\Gamma}$ is
    compact, and without loss of generality we may (and will) assume
    that
    $\Gamma$ is Zariski dense in  $G$. The theorem will be proved
    in two steps.

{\bf{Step 1.}}    First we will prove that $R(G)  \,\cap\, \Gamma$
is a lattice in $R(G)$.
    Let $R_u(G)$ be the unipotent radical of $G$ (i.e.\ $R_u(G)$ is the
    group of all unipotent elements in $R(G)$. Denote by  $\varphi: G
\rightarrow G/R_u(G)$, $\psi: G \rightarrow G/R(G)$ and $\chi:
G/R_u(G) \rightarrow G/R(G)$
     the natural $S$-rational homomorphisms. (We have
$\chi\circ\varphi=\psi$.)  Using verbatim
     an argument from the proof of a theorem of Zassenhaus (see \cite{R,
     Section 8.14}), one proves the following 

     {\bf{Claim}}: There exists a
neighborhood $\Omega$ of $e$ in $G/R_u(G)$ such that if $K$ is a
bounded subset in $\varphi^{-1}(\Omega)$, then $K^{(n)} \rightarrow
e$, where $K^{(0)}=K$ and $K^{(n)} = [K,K^{(n-1)}]$ for all $n
\geq 1$.
\medskip

 Now let $\Delta$ be the Hausdorff closure in $\varphi(G)$ of
the subgroup generated by $\varphi(\Gamma)  \,\cap\, \Omega$. Since
$\varphi(\Gamma)  \,\cap\, \Omega$ is dense in
$\overline{\varphi(\Gamma)}  \,\cap\, \Omega$, the group  $\Delta$ is
open in $\overline{\varphi(\Gamma)}$. (As usually in this paper,
here and hereafter $\overline{X}$ stands for the closure of $X
\subset G$ with respect
    to the Hausdorff  topology.) Let $K
\subset \Gamma$ be a finite set such that $\varphi(K) \subset
\Omega$. In view of the above claim there exists $n_0 > 0$ such
that $K^{(n_0)} = e$. Therefore the group generated by $K$ is
nilpotent \cite{R, Lemma 8.17}, which implies that $\Delta$ is
solvable \cite{R, Lemma 8.4}. Since $\varphi(\Gamma)$ is Zariski
dense in $\varphi(G)$,  the Lie algebra of
$\overline{\varphi(\Gamma)}$ is solvable and
$\Ad\big(\varphi(G)\big)$-invariant. Therefore
$Z\big(\varphi(G)\big)  \,\cap\, \overline{\varphi(\Gamma)}$ is open in
$\overline{\varphi(\Gamma)}$, where $Z\big(\varphi(G)\big)$
denotes the center of $\varphi(G)$. Let $H$ be a maximal
semisimple subgroup of $\varphi(G)$. Then $\varphi(G)$ is an
almost direct product of $H$ and $Z\big(\varphi(G)\big)$. Since
$Z\big(\varphi(G)\big)  \,\cap\, \overline{\varphi(\Gamma)}$ is open in
$\overline{\varphi(\Gamma)}$, $H  \,\cap\, \overline{\varphi(\Gamma)}$
is a discrete normal subgroup of  $\overline{\varphi(\Gamma)}$.
Also, $H  \,\cap\, \overline{\varphi(\Gamma)}$ is Zariski dense in $H$
because the commutator of $\overline{\varphi(\Gamma)}$ is a
Zariski dense subgroup of $H$. Therefore $\chi\big(H  \,\cap\,
\overline{\varphi(\Gamma)}\big)$ is discrete in $\psi(G)$ and
normalized by $\psi(\Gamma)$.

Assume by contradiction that $\psi(\Gamma)$ is not discrete. Then
there exists a neighborhood $W$ of $e$ in
$\overline{\psi(\Gamma)}$ and an infinite Zariski closed subgroup
$L$ of $\psi(G)$ such that if $W'$ is any neighborhood of $e$
contained by $W$ then the Zariski closure of the subgroup
generated by $W'$ coincides with $L$. Now for every $g \in
\chi\big(H  \,\cap\, \overline{\varphi(\Gamma)}\big)$ there exists
a neighborhood $W_g$ of $e$ in $W$ which centralizes $g$.
Therefore $L$ centralizes $\chi\big(H  \,\cap\,
\overline{\varphi(\Gamma)}\big)$, which implies that $L$ is central
(because $\chi\big(H  \,\cap\, \overline{\varphi(\Gamma)}\big)$ is
Zariski dense in  $\psi(G)$). This contradicts the fact that
$\psi(G)$ has finite center. Therefore $\psi(\Gamma)$ is discrete
and, in view of \cite{R, Theorem 1.13}, $R(G)  \,\cap\, \Gamma$
is a lattice in $R(G)$.

\medskip

{\bf{Step 2.}} In order to complete the proof of the proposition,
it is enough to prove that if $G$ is solvable and $\Gamma$ is a
closed subgroup of $G$ such that $G/\Gamma$ admits finite
$G$-invariant measure, then $G/\Gamma$ is compact.

For  $G$ and $\Gamma$ as above, denote $G = G_\infty \times
G_f$, where $G_f = \prod_{v \in S_f} G_v$. Let
$\alpha_\infty :G \rightarrow G_\infty$ and $\alpha_f:G
\rightarrow G_f$ be the natural projections. As in Step 1, using
the Borel Density Theorem we reduce the proof to the case when
$\Gamma$ is Zariski dense in $G$. Let $G_{f}^*$ be an open compact
subgroup of $G_f$. (The group $G_{f}^*$ exists because $G_{f}$ is
a direct product of $p$-adic Lie groups.) Then $(G_\infty \times
G_{f}^*)/\big(\Gamma \,\cap\, (G_\infty \times G_{f}^*)\big)$ has finite
$(G_\infty \times G_{f}^*)$-invariant measure,
   which, in view of the compactness of
$G_f^*$ and the cocompactness of  lattices in  real
solvable Lie groups \cite{R, Ch.\,3}, implies that
$G_\infty/\alpha_\infty(\Gamma  \,\cap\, (G_\infty \times G_f^*))$ is
compact and, therefore, $(G_\infty \times G_{f}^*)/
\big(\Gamma \,\cap\, (G_\infty \times G_{f}^*)\big)$
    is compact. Assume for a moment that
$G_f/\overline{\alpha_f(\Gamma)}$ is compact. Then $G/G_1$ is
compact, where $G_1 \df
\alpha_f^{-1}\left(\overline{\alpha_f(\Gamma)}\right)$. Since $G_1
= \big(G_1  \,\cap\, (G_\infty \times G_f^*)\big)\Gamma$, we get
that $G_1/\Gamma$ is compact and, therefore, $G/\Gamma$ is
compact.

So, it remains to prove that if $G = G_f$ and $\Gamma$  is Zariski
dense in $G$, then $G/\Gamma$ is compact. Let $\Gamma_u$ denote
the group  generated by all unipotent algebraic subgroups
contained in $\Gamma$. Then $\Gamma_u$ is normal in $G$, and
replacing $G$ by $G/\Gamma_u$ we reduce the proof to the case when
$\Gamma_u =e$. Let $P$ be an open subgroup of $G$ containing
$R_u(G)$ and such that $P/R_u(G)$ is compact. Then $P/(P  \,\cap\,
\Gamma)$ admits a finite $P$-invariant measure and, in view of
\cite{MT3, Lemma 1.10}, $P  \,\cap\, \Gamma \supset R_u(G)$.
Therefore $R_u(G) = e$ and $G$ is an abelian group. This proves
that the quotient $G/\Gamma$ is a locally compact topological
group with finite Haar measure. Therefore $G/\Gamma$ is compact.
\qed\enddemo

   \demo{{\bf 9.3.} Proof of Theorem 9.1}
   By  Proposition 9.2, $R(G)\cap \Gamma$ is a cocompact lattice in $R(G)$. Let $N$ be a
maximal subgroup in the class of all normal algebraic subgroups of
$G$ such that $N \cap \Gamma$ is a cocompact subgroup of $N$. Let
$\bold{N}$ be the Zariski closure of $N$ in $\bold{G}$. By the
general structure theory of algebraic groups \cite{Bo}, $\bold{H}
= \bold{G}/\bold{N}$ is a $\bq_S$-algebraic group, $\bold {H} =
\prod_{v \in S} \bold{H}_v$ where each $\bold{H}_v$ is a
semisimple
      group and $\bold{H}_v(\bq_{v})$ has no compact  factors, and
there exists a
$\bq_S$-homomorphism $\varphi: G \rightarrow H$ (where $H =
\bold{H}(\bq_S)$)  such that $\varphi(G)$ has finite index in $H$.
Denote $\varphi(\Gamma)$ by $\Sigma$. Since $N \cap \Gamma$ is
cocompact in $N$, $\Sigma$ is a lattice in $H$ and the natural map
$\widetilde{\varphi}:G/\Gamma \rightarrow H/\Sigma$ is proper.
With $\phi$ as in the formulation of the theorem,  note that
$\varphi \circ \phi$ is a rational parametrization for
$\varphi(U)$ of degree depending only on the degree of $\phi$ and
the homomorphism $\varphi$. Therefore without loss of generality
we can reduce the proof to the case when  every $G_v$ is a
semisimple group without compact factors. Furthermore, we may
assume that
        $\Gamma$ is an irreducible lattice in $G$.

        With the above assumptions, let $\rank_{S}G = \sum_{p \in
        S}\rank_{\bq_v}G_v$ be the $S$-rank of $G$. If $\rank_S G = 1$, then
        either $G$ is a real rank-one semisimple group
        and the theorem is proved in \cite{D2} (see also \cite{D3, Remark
3.7}), or $G$ is a $p$-adic Lie
        group and, therefore, $\Gamma$ is a cocompact lattice  \cite{T} in $G$
        and there is nothing to prove.
        It remains to consider the case $\rank_S G  > 1$.   In
view of Margulis' Arithmeticity Theorem
        \cite{Zi, Theorem 10.1.12} we may assume that $\bold{G}$ is a
$\bq$-algebraic subgroup of
        $\bold{SL}(m)$, $\Gamma = G  \,\cap\, \SL(m, {\bz}_S)$,
         and, after eventually replacing $\bold{G}$ by its image under a
$\bq$-irreducible representation,
         we may also assume that $\bold{G}(\bq)$ acts irreducibly on $\bq^m$.

Writing $y$ in the form $g\Gamma$ for some $g\in G$ and
      applying  Theorem 8.4 to the map $h(x) = \phi(x)g$,
         we get a  compact $L_0 = Q_\tau\supset L$ and  constants $C,\alpha>0$
         such that for any $\tau >0$ and any ball $B \subset X$
one has
$$
\lambda\left(\big\{x\in B\bigm| \phi(x)g{\Cal
O}_S^m \in Q_\vre \big\}\right) \geq (1 - C\vre^\alpha
)\lambda(B),
$$
whenever $\phi(B)y  \,\cap\, L_0 \neq \varnothing$.
        Choosing $\vre$ such that $\beta \geq C\vre^\alpha$,  we get a
compact $M_\beta \df Q_\vre$ satisfying (9.1).

      If $Uy  \,\cap\,
L_0 = \varnothing$, then in view of Theorem 8.4 there exists a
proper nonzero vector subspace $V \subset \bq^m$ such that $\bq_S
V$ is invariant under the action of $g^{-1}Ug$. (We consider
$\bq^m$ diagonally imbedded in $\bq_S^m$, which justifies the
expression $\bq_SV$.) Let $\bold{P}$ be the Zariski closure of the
stabilizer of $V$ in $\bold{G}(\bq)$ under the natural action of
$\bold{G}(\bq)$ on $\bq^m$. Since $\bold{G}(\bq)$ acts irreducibly
on $\bq^m$, $\bold{P}$ is a proper  $\bq$-algebraic subgroup of
$\bold{G}$. Therefore so  is the subgroup $\bold{P}_u$ of
$\bold{P}$ spanned by all unipotent elements of $\bold{P}$. It is
easy to see that $g^{-1}Ug \subset \bold{P}_u(\bq_S)$. By the
Borel--Harish-Chandra theorem $\bold{P}_u(\Cal {O}_S)$ is a
lattice in $\bold{P}_u(\bq_S)$. Therefore the group $H =
g\bold{P}_u(\bq_S)g^{-1}$ satisfies the requirements  of the
formulation of the theorem. \qed\enddemo

We remark that to derive Theorem 0.2 from the above theorem it
suffices to take $K = \bq$, $L = \{y\}$ and, in the case $Uy \, \cap\, L_0 \ne
\varnothing$, choose $R$ such that the intersection of
$\phi\big(B(0,R)\big)$ with $L_0$ is nonempty.

\subhead{9.4}\endsubhead Let $\sigma$ be a Haar measure on $U$.
   Given a bijective parametrization $\phi:
X\to U$, let $\phi_*\lambda$ be the pushforward of $\lambda$ to
$U$ via $\phi$. Note that $\phi$ can be chosen in such a way that
$\sigma = \phi_*\lambda$. If $U$ is abelian,
      $\phi$ can be taken to be the trivial parametrization (because $\exp$ is
an isomorphism of locally compact groups). In the general case,
since $U$ is unipotent and algebraic, there exist abelian
algebraic subgroups $U_1,\dots,U_s$ of $U$ such that the map
$\psi: U_1 \times\dots\times U_s \rightarrow U$,
$\psi(u_1,\dots,u_s) = u_1\cdots u_s $, is bijective, and for any
$i$, $1 \leq i <s $, the product  $U_1\cdots U_i$ is a normal
subgroup of $U$ and $U_1\cdots U_{i+1}/U_1\cdots U_i$ is a central
subgroup of $U/U_1\cdots U_i$. For any $j$ let $\phi_j$ be the
trivial parametrization of the abelian group $U_j$. Then a simple
computation shows that $\phi =\psi\circ(\phi_1,\dots,\phi_s) $ is
such  that $\sigma = \phi_*\lambda$, proving the claim. Using the
bijective map $\phi$ we translate the metric from $X$ to $U$. In
view of the preceding discussions, Theorem 9.1 immediately implies
the following

\proclaim{Theorem} Let $G$  and $\Gamma$  be as in Theorem
9.1. Then every compact $L \subset G/\Gamma$ is contained in a
compact $L_0$ with the following property: given $\beta
>      0$ there exists a compact subset $M_\beta$ of $G/\Gamma$  such that
for any $y \in G/\Gamma$ and any unipotent algebraic subgroup $U$
of $G$  the following is satisfied:

{\rm(i)} If $Uy  \,\cap\, L_0 \neq \varnothing$ and $B$ is a ball in $U$
such that  $\phi(B)y  \,\cap\, L_0 \neq
\varnothing$, then
$$
\frac{\sigma \left(\big\{u\in B\mid uy \in M_\beta
\big\}\right)}{\sigma(B)} \geq 1 -\beta  \,,\tag 9.2
$$
where $\sigma$ is a Haar measure on $U$. In particular, {\rm (9.2)}
holds if $y \in L$ and $e\in B$.

{\rm(ii)} If $Uy  \,\cap\, L_0 = \varnothing$, then there exists a closed
proper subgroup $H$ of $G$ containing $U$ such that the orbit
$Hy$ is closed and carries a finite $H$-invariant Borel measure.
\endproclaim

\subhead{9.5}\endsubhead In order to prove Theorem 0.3 we will need
a  version of the Birkhoff
ergodic theorem. Let
$G, \Gamma, H$ and $\mu$ be the same as in the formulation of
Theorem 0.3, and let $U$ be a unipotent algebraic subgroup of  $G$
contained in $H$.  Let $\mu = \int_{(A,\rho)}\mu_a\, d\rho$ be the
decomposition of
        $\mu$ into $U$-invariant ergodic locally finite measures, where
        $(A,\rho)$ is a measure space parametrizing the ergodic
        components $\mu_a$. For almost every $y \in G/\Gamma$ there
exists a well
defined ergodic component $\mu_{a(y)}$, where $a(y)\in A$, whose
support contains $y$.

Fix an imbedding $G \subset \GL(m, \bq_{S})$ and a maximal
unipotent subgroup $W = \prod_{v \in S} W_v \subset \GL(m,
\bq_{S})$ which contains $U$. There exists an element $g = (g_v)
\in \GL(m, \bq_{S}) $ such  that $gWg^{-1} = W$, $W = \{x \in
\GL(m, \bq_{S}) | \lim_{n \rightarrow \infty}g^{-n}xg^n = e \}$,
   and
for every $v \in S$ there exists $\pi_v \in \bq_{v}$ such that
$|\pi_{v}| > 1$ and all eigenvalues of $g_v$ are powers of
$\pi_v$. It is easy to see \cite{MT2, Proposition 2.2} that the
sequence $\Lie ( g^{-n}Ug^n )$ has a limit $\Lie(U_0)$, where
$U_0$ is a $\bq_{S}$-algebraic subgroup of $W$,  in the
Grassmannian variety $\Gr (\Lie W)= \prod_{v \in S} \Gr (\Lie W_v)
$ of $\Lie W$. Also it is known \cite{MT2, Proposition 2.8} that
there exists an $\Int(g)$-invariant Zariski closed subset $V
\subset W$ such that the maps $U \times V \rightarrow W,
(u,v)\rightarrow uv$, and $U_0 \times V \rightarrow W, (u_0,v)
\rightarrow uv$ are bijective. Denote by $p:U \rightarrow U_0$ the
projection of $U$ onto $U_0$ parallel to $V$. Note that $p$ is
bijective. Put $\psi = p^{-1} \circ \Int(g) \circ p $ . Then
$\psi:U \rightarrow U$ and $\psi$ acts as an expansion on $U$, in
particular, if $B$ is a relatively compact neighborhood of $e$ in
$U$ then

$$
U = \bigcup_{i \geq 0} B_i \,,\tag 9.3
$$
    where $B_0 =B$ and $B_{i+1} = \psi
(B_i)$ for all $i \geq 0$. Note that $\lim_{i \rightarrow
\infty}\psi^{-i}(u) = e$ for all $u \in U$. Further on we fix a
set $B$ as above. Since the Jacobian of $\psi$ is  constant, for
every $u \in U$ we have

$$
\lim_{i \rightarrow \infty}\frac{\sigma(B_i \triangle
uB_i)}{\sigma(B_i)} = \lim_{i \rightarrow \infty}\frac{\sigma(B
\triangle \psi^{-i}(u)B)}{\sigma(B)} =0\,\tag 9.4
$$
and
$$
\sup_{i}\frac{\sigma(B_i^{-1}B_i)}{\sigma(B_i)} =
\frac{\sigma(B^{-1}B)}{\sigma(B)} < \infty \,,\tag 9.5
$$
where $\sigma$ is the Haar measure on $U$ and $C \triangle D$
denotes the symmetric difference between two sets $B$ and $D$.

In view of {\rm (9.3)}, {\rm (9.4)} and {\rm (9.5)}, the following
result directly follows from \cite{Te, Corollary 3.2, Ch.\,6} (see
also \cite{MT2, Proposition 7.1}):

\proclaim{Proposition} With the above notation and assumptions,
let $f$ be a continuous $\mu$-integrable function on $G/\Gamma$.
Then

$$
\lim_{n \rightarrow \infty}\frac{1}{\sigma(B_n)}\int_{B_n}f(gy)
\,d\sigma(g) = \int_{G/\Gamma}f(z)\,d\mu_{a(y)}(z)
$$
for almost all $y \in G/\Gamma$. Furthermore, the limit function
$$
f^*(y) \df \int_{G/\Gamma}f(z)\,d\mu_{a(y)}(z) \tag 9.6
$$
is $\mu$-integrable and $U$-invariant.
\endproclaim

\subhead{9.6}\endsubhead   We  now proceed with the

\demo{Proof
of Theorem 0.3}  In order to prove the theorem, as in
the real case considered in
        \cite{D1}, it is enough to find
        a function $h \in L^1(G/\Gamma, \mu)$ which is $H$-invariant and
$h(y) > 0$ for
        $\mu$-almost all $y \in G/\Gamma$. Indeed, if $h$ is such a function,
        then
        the sets $$Y_i = \left\{y \in G/\Gamma \mid h(y) \geq
\tfrac{1}{i}\right\}$$
satisfy the conditions
        in the formulation of the theorem.

        Among the  unipotent algebraic subgroups of $G$  contained in
$H$ we fix a maximal one
and denote it by $U$.
        It is well known that the minimal normal subgroup of $H$
        containing $U$ coincides with $H$ itself (see \cite{Bo}).
Therefore, in view of the
Mautner
        phenomenon for products of real and $p$-adic Lie groups
        \cite{MT3, Proposition 2.1}, if $f \in L^1(G/\Gamma, \mu)$ is
        $U$-invariant, then there exists an $H$-invariant $\mu$-integrable
function on
$ G/\Gamma$
        which coincides with $f$ almost everywhere. So, it is
        enough to prove the theorem for $H = U$.

Let $f$ be a positive continuous $\mu$-integrable function on
$G/\Gamma$. Applying Proposition 9.5, we get a $U$-invariant
function $f^*$ defined by  formula {\rm (9.6)}. It is enough to
prove that  $f^* > 0$ $\mu$-a.e. Note that
$$
\frac{1}{\sigma(B_n)}\int_{B_n}f(gy)\, d\sigma(g) =
\frac{1}{\sigma(B)}\int_{B}f\big(\psi^n(g)y\big)\, d\sigma(g)\,,
$$
where $\psi$ is as in \S 9.5. It follows from Theorem 9.1 and the
facts that $\sigma = \phi_*\lambda$ for some rational
parametrization $\phi$ of $U$ (see \S 9.4) and all $\psi^{n} \circ
\phi$ have the same degree (because $\psi$ is linear), that there
exists a compact $M \subset G/\Gamma$ such that for any positive
$n$
$$
\frac{\sigma\big(\{g \in B \mid \psi^n(g)y \in M
\}\big)}{\sigma(B)}
>\frac{1}{2}.
$$
Since $f$ is positive, the above formula implies that $f^*(y) > 0$
$\mu$-a.e., which completes the proof of the theorem. \qed\enddemo

\example{9.7.\ Remark} Using the methods from \cite{MT3} it is
easy to see that
Theorems 9.1, 9.3, 0.3 and Proposition 9.5 remain valid for the
larger class of so-called {\sl almost linear groups\/}, that is,
when $G$ is a finite direct product of a connected real Lie group
and finite central extensions of closed linear $p$-adic groups.
\endexample


\heading{10.  $S$-arithmetic \da}
\endheading


\subhead{10.1}\endsubhead
In this section we present a motivation for the definitions
of VWA and VWMA vectors given in the introduction, and state
the main \di\ result of the
present paper, of which Theorem 0.4 is a special case.
Here and for the rest of the paper we fix
a set $S$ 
of cardinality  $\ell$ consisting of distinct normalized
valuations of $\bq$ (and not necessarily containing $\infty$) ,
and let
let $\bq_S = \prod_{v\in S}\bq_{v}$.


        We
will  interpret elements $$\vy\ =\ (\vy^{(v)})_{v\in S}\ =\
(y_1,\dots,y_n)$$ of $\bq_S^n$, where $\vy^{(v)}=
(y^{(v)}_1,\dots,y^{(v)}_n)\in \bq_{v}^n$ and $y_i =
(y^{(v)}_i)\in \bq_S$,
        as linear forms on $\bq_S^n$, and
       will study their values  $\vy\cdot\vq = y_1q_1
+\dots+y_nq_n$ at integer points $\vq = (q_1,\dots,q_n)$.  The approximation
properties of our interest will be related to these values being close
(in terms of the $S$-adic absolute value $|\cdot|$ on $\bq_S$) to integers.
Alternatively, one could consider  a dual case when one approximates (in
terms of the
$S$-adic norm)
$\vy\in\bq_S^n$ by
rational
vectors. See \cite{KLW} and
\cite{KM, Remark 6.2} for a discussion of this set-up in the real case.

For $\vy$, $\vq$ as above and for $q_0\in\bz$, it will be convenient to
use the notation
$$\tilde \vy\df(1,y_1,\dots,y_n)\quad\text{and}\quad
\tilde \vq \df (q_0,q_1,\dots,q_n)\,,$$
      so that $q_0 + y_1q_1 +\dots+y_nq_n = q_0
+
\vy\cdot\vq$ is written as $\tilde\vy\cdot\tilde\vq$.
Also, by the absolute value $|\cdot|$
of integers and the norm $\|\cdot\|$ of integer vectors we will always
mean those coming from the infinite valuation. Hopefully
it will cause no confusion.

\subhead{10.2}\endsubhead A natural starting point in the theory of
simultaneous
\da\ is usually  a Dirichlet-principle-type result. Let us work it out.
       The goal is to find the optimal exponent $\beta$
such that for any $\vy\in \bq_S^n$ and any
       $N>0$ one is guaranteed to have
two different
integer vectors
$\tilde\vq_1,\tilde\vq_2$  of 
      norm $\le
N$  such that $\tilde\vy\cdot\tilde\vq_1$ and
$\tilde\vy\cdot\tilde\vq_2$ are at most
$\const(\vy) N^{-\beta}$ apart.

It turns out that the answer depends on
whether or not $S$ contains the Archimedean valuation $v = \infty$.
Indeed, note that one has
$$|\tilde\vy^{(v)}\cdot\tilde\vq|_v \le \max(\|\vy^{(v)}\|_v,1)\tag 10.1$$
for any ultrametric $v$ and any
integer vectors
$\tilde\vq$. Therefore:

\roster
\item"$\bullet$" If all the valuations in $S$ are ultrametric,
the $(2N+1)^{n+1}$ values of $\tilde\vy\cdot\tilde\vq$ for all $\tilde\vq$ of
       norm $\le N$   are in the
ball of radius
$\max(\|\vy\|,1)$ in $\bq_S$.
       The latter can be partitioned into
$\const\cdot N^{n+1}$ balls of radius $\const(\vy)
N^{-\frac{n+1}\ell}$.
Thus
for any
$\vy\in \bq_S^n$ and any
        $N>0$ one can find $\tilde\vq \in\bz^{n+1}\nz$
with $\|\tilde\vq\|\le N$ and $$
|\tilde\vy\cdot\tilde\vq| \le
\const(\vy)\cdot N^{-\frac{n+1}\ell}\,.
$$
       \item"$\bullet$" If  $v = \infty$,  there is clearly no
universal upper bound similar to (10.1). However one can
for any  given $\vq\in\bz^n$ choose $q_0\in\bz$ such that
$$|\tilde\vy^{(\infty)}\cdot\tilde\vq|_\infty = |q_0 +
\vy^{(\infty)}\cdot\vq|_\infty \le 1\,.
$$
Thus, taking all $\vq$ of
       norm $\le N$, one is only guaranteed to have $(2N+1)^{n}$ values of
$\tilde\vy\cdot\tilde\vq$  in the
ball of radius
$\max(\|\vy\|,1)$ in $\bq_S$. Partitioning it into
$\const\cdot N^{n}$ balls of radius $\const(\vy)
N^{-\frac{n}\ell}$, one gets a nonzero integer vector
$\tilde\vq$
with $\|\vq\|\le N$ and $$
|\tilde\vy\cdot\tilde\vq| = |q_0 + \vy\cdot\vq|\le
\const(\vy) N^{-\frac{n}\ell}\quad \text{for some }q_0\in\bz\,.
$$
Note that the absolute value of $q_0$ above, and hence the norm of $\tilde\vq$,
is bounded from above by
$\const(\vy) N$.
\endroster

It will be convenient to define
$$
i_{\sssize S}\df\cases 1\quad\text{ if }\infty\notin S\\
0\quad\text{ if }\infty\in S\,.\endcases
$$
Then it follows that for any $\vy\in\bq_S^{n}$
the supremum of
$w > 0$ for which
there exist infinitely many $\tilde\vq\in\bz^{n+1}$ with
$$|\tilde\vy\cdot\tilde\vq|^\ell \le \|\tilde\vq\|^{-w}
$$
      is not less than
$n +i_{\sssize S}$. On the other hand, it
can be easily shown using the Borel-Cantelli Lemma that
the above supremum is equal to
$n +i_{\sssize S}$ for
almost every $\vy\in\bq_S^{n}$ (with respect to Haar measure $\lambda_S$
on
$\bq_S^{n}$).    
Thus it is natural to say that
$\vy\in \bq_S^n$ is
{\sl \vwa\/}, or
{\sl VWA\/},  if the above supremum is strictly bigger  than
$n +i_{\sssize S}$; in other words, if for some $\vre > 0$ there are
infinitely many solutions
$\tilde\vq\in\bz^{n+1}$ to
$$|\tilde\vy\cdot\tilde\vq|^\ell \le \|\tilde\vq\|^{-(n +i_{\sssize S})(1 +
\vre)}\,.\tag 10.2
$$
Note  that  in
the  case when $\infty\in S$,  any solution $\tilde\vq$ of (10.2)
automatically satisfies
$$|q_0| \le 1 + n\|\vy^{(\infty)}\|_\infty\|\vq\|\,,
\tag 10.3
$$
and hence $\|\tilde\vq\|$ in the right-hand side of (10.2) 
definition of VWA vectors
can be
replaced by  $\|\tilde\vq\|$,
agreeing with (0.4).
\medskip

\subhead{10.3}\endsubhead The next step is to define {\sl \vwma\/}, or
{\sl
VWMA\/}, vectors $\vy\in \bq_S^n$.  To do this,
one would like to
   replace
the left hand side of (10.2) by the product of norms of all the
components
of $\tilde\vy\cdot\tilde\vq$,
and the norm of
$\tilde\vq$ in (10.2)
with
the geometric mean of its coordinates. However one needs to be
careful and keep in
mind  the dichotomy in the Dirichlet-principle
argument.

Namely, if  $\infty\notin S$ (when $i_S = 0$) one can indeed replace
$\|\tilde\vq\|^{n + 1}$ by
$\Pi_+(\tilde\vq)$, and thus define $\vy\in \bq_S^n$ to be
   VWMA  if for some $\vre > 0$ there are
infinitely many solutions
$\tilde\vq\in\bz^{n+1}$ to
$$c(\tilde\vy\cdot\tilde\vq)  \le \Pi_+(\tilde\vq)^{-(1 +
\vre)}\,.\tag 10.4
$$
       On the other hand,  if
$\infty\in S$ it seems tempting to 
define $\vy$ to be  VWMA
if for some $\vre > 0$ there are
infinitely many 
$\vq$ such that
$$c(\tilde\vy\cdot\tilde\vq)  = c(q_0 + \vy\cdot\vq)  \le \Pi_+(\vq)^{-(1 +
\vre)}
\tag 10.4$'$
$$
holds for some $q_0\in \bz$. Indeed, this coincides with the
standard definition when $S =
\{\infty\}$, cf.\ \cite{KM}. However it is not hard to see that, whenever
$S$ contains both finite and infinite valuations,
for any
$\vre > 0$ the set of
$\vy\in \bq_S^n$ for which (10.4$'$) admits infinitely many solutions
has full measure. Indeed, the trouble here comes from the fact that an
upper estimate for $c(\tilde\vy\cdot\tilde\vq)$ does not imply a bound similar
to (10.3), that is, a bound on $|q_0|$ in terms of $\vq$. And one can
easily show that for any fixed $\vq$ the set of $\vy\in \bq_S^n$ for
which there exists
$q_0$ satisfying (10.4$'$) has full measure.

It follows that in order to achieve a multiplicative analogue of (10.4)
in the case $\{\infty\}\subsetneq S$, one needs to take
special precautions in the case when $|q_0|$ is much bigger than the norm
of $\vq$.
Namely,
in the case $\infty\in S$ we will define $\vy\in \bq_S^n$ to be {\sl
VWMA\/}  if for some $\vre > 0$ there are
infinitely many solutions
$\vq$ to
$$c(\tilde\vy\cdot\tilde\vq)    \le \Pi_+(\vq)^{-(1 +
\vre)}|q_0|_+^{-\vre}\,.\tag 10.4$_\infty$
$$
Put together with (10.4), the latter inequality can be written
in the form  (0.5), or, equivalently, in the unified form as
$$c(\tilde\vy\cdot\tilde\vq)    \le \Pi_+(\vq)^{-(1 +
\vre)}|q_0|_+^{-(i_S + \vre)}\,.\tag 10.5
$$

Several remarks are in order. First,
note that  in the case $S =
\{\infty\}$ (10.4$_{\infty}$) can be replaced by (10.4$'$),
perhaps with a slightly different value of $\vre$: indeed, any solution
of (10.5) will satisfy (10.3), hence $|q_0|_+$ is bounded
from above by some  power of $\Pi_+(\vq)$.
Similarly, it can be easily seen that infinitely many solutions
to (10.2) imply infinitely many solutions
to (10.5) (with the same $\vre$ if $\infty\notin S$, and,
in view of (10.3) and $\Pi_+(\vq) \ge \|\vq\| $,
   perhaps with a different  $\vre$
if $\infty\in S$).
And yet, VWMA as defined above happens to be a zero measure condition.
This can be shown directly using a Borel-Cantelli argument, and it
will also be an implication of Theorem 10.4 below.

\subhead{10.4}\endsubhead
Recall
      that   a measure
$\mu$
on $\bq_S^{n}$ is called
   extremal (resp.,  strongly extremal) if
$\mu$-almost every point of $\bq_S^{n}$
is
not VWA (resp., not VWMA).
Here is the main theorem of the section:

\proclaim{Theorem}  For $v\in S$, let $X_v$ be a metric space
with a measure $\mu_v$
such that $X = \prod_{v\in S}X_v$ is \be\
and $\mu= \prod_{v\in S}\mu_v$ is   Federer, and let
$\vf = (\vf^{(v)})_{v\in S}$, where
$\vf^{(v)}$ are
continuous maps from
$X_v$ to
$\bq_{v}^n$ which are $\mu_v$-good
and $\mu_v$-nonplanar at $\mu_v$-almost every point of
$X_v$.  Then $\vf_*\mu$
is strongly extremal.
\endproclaim

       It is clear  from Theorem 4.3, as well as from Examples 1.6 and 1.7,
that
Theorem 0.4
is a special case of the above result.

\subhead{10.5}\endsubhead Let us also remark that Theorem 10.4  generalizes the
main result of
\cite{KLW}. Indeed, in the latter paper  a certain class of measures
on $\br^n$ 
was introduced, and it was proved that
measures from that class are strongly extremal.
Specifically, following \cite{KLW} let us say that a measure  $\mu$ on
$F^n$, where $F$ is a locally compact field, is

\roster
\item"$\bullet$" {\sl nonplanar\/} if $\mu(L) = 0$ for any proper affine
subspace of $F^n$;

       \item"$\bullet$" {\sl decaying\/} if for  for $\mu$-a.e.\
$\vy \in F^n$  there exist a neighborhood $V$ of $\vy$ and $C,\alpha >
0$ such that all affine functions are \cag\ on $V$ with
respect to
$\mu$;

      \item"$\bullet$" {\sl friendly\/} if it is Federer, nonplanar and
decaying.
\endroster

Comparing this with \S 4.2, one easily observes that $\mu$ is  decaying
if and only if the identity map
$F^n\to F^n$ is $\mu$-good at $\mu$-almost every point.  It is also not
hard to see that the nonplanarity of  $\mu$
forces  the
aforementioned identity map
      to be $\mu$-nonplanar at $\mu$-almost every point
(converse is true under the additional
assumption that
$\mu$ is decaying).

It is now clear that Theorem 10.4 immediately implies

\proclaim{Corollary}  Let $\mu= \prod_{v\in S}\mu_v$, where $\mu_v$ is a
friendly measure on $\bq_{v}^n$ for every
$v\in S$.  Then $\mu$ is strongly extremal.
\endproclaim

Thus \cite{KLW, Theorem 1.1} is a special case of Theorem
10.4. (As was mentioned before, our proof is also a generalization of
the argument from
\cite{KLW}, which, in turn, generalizes the one from \cite{KM}.)

\subhead{10.6}\endsubhead It is not hard to see that many examples of friendly
measures on
$\br^n$ exhibited in \cite{KLW} can be constructed on a vector space over
arbitrary locally compact valued field $F$. For instance, fix a valuation
$|\cdot|$ on $F$ inducing the metric ``$\dist$'' on $F^n$, and
      say that a  map $\vh : F^n \to F^n$ is a  {\sl contracting similitude\/}
with {\sl contraction rate\/} $\rho$ if
$0 < \rho < 1$ and
$$\dist\big(\vh ( \vx ),\vh ( \vy )\big) = \rho \dist( \vx, \vy )
  \quad\forall\,\vx,\vy \in F^n\,.
$$
It is known, see \cite{H, \S 3.1}, that for any finite family
$\vh_1,\dots, \vh_m$ of  contracting similitudes
there exists a unique nonempty compact set $Q$,
called the {\sl limit set\/} of the family,
such that
$$
Q = \bigcup_{ i = 1 }^m\vh_i ( Q ).
$$
Say that
$\vh_1, \dots,  \vh_m$ as above  satisfy the {\sl open set
condition\/} 
if there exists an open subset $U \subset F^n$ such that
$$\vh_i ( U ) \subset U  \text{ for  all  } i=1, \dots,
m\,,$$
          and
$$i \ne j \Longrightarrow \vh_i ( U )  \,\cap\, \vh_j(U) =
\varnothing\,.$$

J.~Hutchinson \cite{H, \S 5.3} proved\footnote{Hutchinson stated his
results for
the case $F = \br$, but the proofs apply verbatim to the case of arbitrary
locally compact valued field.} that if
$\vh_i$, $i = 1 , \dots , m$, are contracting similitudes with contraction
rates $\rho_i$ satisfying the open set
condition, and if
$s > 0$ is the unique solution of
$
\sum_i \varrho_i^s = 1\,,
$
called the {\sl similarity dimension\/}
of the family $\{\vh_i\}$,
then the $s$-dimensional Hausdorff measure $\Cal H^s$ of $Q$ is
positive and finite. Let us also say that the family $\{\vh_i\}$
is {\sl irreducible\/} if there does not exist a finite
$\{\vh_i\}$-invariant collection of proper affine subspaces
of
$F^n$.
The proof of \cite{KLW,
Theorem 2.3} applies verbatim  and yields

\proclaim{Proposition} For any completion $\bq_v$ of $\bq$,
          let $\{\vh_1 , \dots , \vh_m\}$ be an irreducible family of
contracting similitudes of $\bq_v^n$ satisfying the
open set condition, $s$ its similarity dimension, $\mu$ the restriction
of $\Cal H^s$ to its limit set. Then $\mu$ is
friendly  (and hence strongly extremal).
\endproclaim

Measures on $\br^n$ obtained via the above construction
have been thoroughly studied; perhaps the simplest
example is given by the $\frac{\log2}{\log 3}$-dimensional Hausdorff measure on
the Cantor ternary set. Similarly one can consider ultrametric analogues of
the Cantor set, for example let
$$
      Q = \left\{\left.\sum_{k = 0}^\infty a_k 3^k\right|a_k =
1,2\right\} \subset
\bz_3\,.
$$
It is a 3-adic version of the Cantor ternary set, which also has \hd\
$s= \frac{\log2}{\log 3}$, and it follows that almost all  numbers in $Q$ (with
respect to the $s$-dimensional Hausdorff measure) are not VWA.

\medskip

\subhead{10.7}\endsubhead We conclude this section with the following
modification of  Theorem 10.4:

\proclaim{Theorem}  For every $v\in S$, let $X_v$ be a metric space
with a measure $\mu_v$
such that $X = \prod_{v\in S}X_v$ is \be\
and $\mu= \prod_{v\in S}\mu_v$ is uniformly  Federer, and let $\vf =
(\vf^{(v)})_{v\in S}$, where
$\vf^{(v)}$ are
continuous maps from
$X_v$ to
$\bq_{v}^n$ such that
       for
$\mu_v$-a.e.\
$x_v\in X_v$
one can find a ball $B_v = B(x_v,r)\subset X_v$  with the following
properties:
$$
\aligned
\text{for some
$C_v,\alpha_v > 0$, any linear combination of }
1,f^{(v)}_1,\dots&,f^{(v)}_n\\\text{  is
$(C_v,\alpha_v)$-good on }
B(x_v,3^{n+1}r) \text{ with respect to  $\mu_v$}\,,\endaligned \tag 10.6
$$
        and
$$
\aligned
\text{   the restrictions of $1,f^{(v)}_1,\dots,f^{(v)}_{n}$ to
}B_v\,\cap\,&\supp\,\mu_v \\ \text{
are linearly independent over $\bq_{v}$\,.}
\endaligned\tag 10.7
$$
       Then $\vf_*\mu$
is strongly extremal.
\endproclaim

\demo{Reduction of Theorem 10.4 to Theorem 10.7} Let $X$ and $\vf$ be as in
Theorem 10.4.  First note that, replacing
$X_v$ by appropriate neighborhoods of its $\mu_v$-generic points for each
$v$, one can without loss of
generality assume that $\mu$ is uniformly Federer. Then for any $v$,
since  $\vf^{(v)}$ is $\mu_v$-good at $\mu_v$-a.e.\ point, one can for
$\mu_v$-a.e.\
$x_v\in X_v$ choose   a neighborhood
$U_v$ of $x_v$ and $C_v,\alpha_v > 0$ such that
       any linear combination of
$1,f^{(v)}_1,\dots,f^{(v)}_n$  is
$(C_v,\alpha_v)$-good on
$U_v$ with respect to  $\mu_v$.
Further,
since  $\vf^{(v)}$ is $\mu_v$-nonplanar $\mu_v$-almost everywhere,  one can
(after throwing away points from a null set) take
a ball $B_v = B(x_v,r)$ such that $B(x_v,3^{n+1}r) \subset U_v$ and
(10.7) holds, and the conclusion follows.
\qed\enddemo

In the next two sections we present the proof of Theorem 10.7,
separately considering the cases of $S$ containing or not containing
the Archimedean valuation. In both  cases the core of the proof is
a generalization of the correspondence between real \da\ and
dynamics on real \hs s.

\heading{11.
Proof of Theorem 10.7 for
$\infty\notin S$}
\endheading

\subhead{11.1}\endsubhead In order to prove Theorem 10.7, we are going to
dynamically interprete the  approximation properties of $S$-adic vectors
defined in the previous section, similarly to the approach of \cite{KM}.
In this section   we suppose that all the valuations in $S$ are
ultrametric, that is, $S = \{p_1,\dots,p_\ell\}$ where $p_1,\dots,p_\ell$
are distinct primes. Up to the end of this section we  will work with
$$S^+ \df S\cup
\{\infty\},\quad {\Cal R} \df \bq_{S^+} = \bq_S\times \br,\quad
\text{and}\quad{\Cal D} \df {\Cal O}_{S^+} =
\bz[\tfrac1{p_1},\dots,\tfrac1{p_\ell}]\,.$$ Then to any $\vy\in
\bq_S^{n}$ we associate  a lattice $u_\vy {\Cal D}^{n+1}$ in
${\Cal R}^{n+1}$, where $u_\vy\in \GL^1(n+1,{\Cal R})$ is defined
by
$$
u_\vy^{(p_j)} = \pmatrix 1 & \vy^{(p_j)}\\ 0 & I_n\endpmatrix,\ {j =
1,\dots,\ell},\quad u_\vy^{(\infty)} = I_{n+1}\,,
$$
with $I_k$ standing for the
$k\times k$ identity matrix. Note that
the $p_j$-adic components of vectors from
$u_\vy {\Cal D}^{n+1}$ are of the form $
\pmatrix\tilde\vy^{(p_j)}\cdot\tilde \vq \\ \vq\endpmatrix$,
where
$\tilde \vq\in {\Cal D}^{n+1}$.

We need to introduce some more notation.
For a vector $\tilde \vt \df
(t_0,t_1,\dots,t_n)\in\br^{n+1}$ we denote $
(t_1,\dots,t_n)$ by $\vt$, and let $$\tilde t = \sum_{i = 0}^nt_i\quad\text{
and
}\quad t =
\sum_{i = 1}^nt_i\tag 11.1 $$
(this convention will be used throughout the next two sections,
so that whenever $t$ and $\vt$, or $\tilde t$ and $\tilde \vt$, appear in the
same context, (11.1) will be assumed). Then, given $\tilde \vt$ as above
and another vector
$\vs = (s_1,\dots,s_\ell)\in\bz_+^\ell$,
define $g_{\vs,\tilde \vt}\in \GL(n+1,{\Cal R})$ by
$$
(g_{\vs,\tilde \vt})^{(p_j)} = \pmatrix p_j^{-s_j} & 0\\ 0 &
        I_n\endpmatrix,\ {j =
1,\dots,\ell},\quad (g_{\vs,\tilde \vt})^{(\infty)} =\,
\text{\rm diag}(e^{-t_0},e^{-t_1},\dots,e^{-t_n})\,.
\tag 11.2
$$

The next lemma shows how a good
approximation for $\vy$ in the sense of (10.4) gives rise to a translation
of $u_\vy {\Cal D}^{n+1}$ by  $g_{\vs,\tilde \vt}$ for some $\vs,\tilde \vt$,
so that
$\delta(g_{\vs,\tilde \vt}u_\vy {\Cal D}^{n+1})$ is small. This
allows one to use
Theorem 8.3 to derive the needed measure estimate.

\proclaim{Lemma} Let $\vre > 0$, $\vy\in \bq_S^{n}$ and
$\tilde \vq\in\bz^{n+1}$ be such that {\rm (10.4)} holds.
For $i = 0,1,\dots,n$ define  $t_i> 0$ by 
$$
|q_i|_{\sssize +} = \Pi_{\sssize +}(\tilde \vq)^{-\frac{\vre}{n+1}}
e^{t_i}\,,\tag 11.3a
$$
and let 
$$
\gamma = \frac{\vre}{(n+1)(1 + \vre)}\,.\tag 11.3b
$$
Then there exists $\vs =
(s_1,\dots,s_\ell)\in\bz_+^\ell$ such that
$$
\delta(g_{\vs,\tilde \vt}u_\vy {\Cal D}^{n+1}) \le \sqrt{n+1}e^{-\gamma\tilde
t}\tag 11.3c
$$
and 
$$
\prod_{j = 1}^\ell p_j^{s_j} \le e^{\tilde t} < \prod_{j = 1}^\ell
p_j^{s_j+1}\tag
11.3d
$$
\endproclaim

\demo{Proof}
Multiplying equalities
(11.3a),  we get
$$
e^{\tilde t} = \Pi_{\sssize +}(\tilde \vq)^{1 + \vre}\tag 11.3e
$$
and
$$
e^{-t_i}|q_i|\le e^{-t_i}|q_i|_{\sssize +}\un{(11.3a)}\le \Pi_{\sssize
+}(\tilde \vq)^{-\frac{\vre}{n+1}}\un{(11.3be)}=  e^{-\gamma\tilde  t},\ i =
0,1,\dots,n\,,
$$
hence $
\|(g_{\vs,\tilde \vt})^{(\infty)}u_\vy^{(\infty)}\tilde \vq\|_\infty =
\|(g_{\vs,\tilde \vt})^{(\infty)}\tilde \vq\|_\infty\le \sqrt{n+1}e^{-\gamma
\tilde t}$.

Now let us  
define
$s_j$, $j = 1,\dots,\ell$, inductively by
$$p_j^{s_j} \le \min\left(\frac{e^{\tilde t}}{\prod_{i = 1}^{j-1}
p_i^{s_i}},
\frac1{|\tilde\vy^{(p_j)}\cdot\tilde \vq|_{p_j}}\right) < p_j^{s_j +
1}\tag 11.4
$$
(where if $j = 1$ we set $\prod_{i = 1}^{j-1}
p_i^{s_i} = 1$). This, in particular, implies that
$$|p_j^{-s_j}\tilde\vy^{(p_j)}\cdot\tilde \vq|_{p_j}  =
p_j^{s_j}|\tilde\vy^{(p_j)}\cdot\tilde \vq|_{p_j} \le
1
       $$ for each
$j$. Taking into account that $|q_i|_{p_j} \le 1$ for all $i$ and
$j$, one concludes
that
        $c(g_{\vs,\tilde \vt}u_\vy
\tilde \vq) \le \sqrt{n+1}e^{-\gamma \tilde t}
$.

It remains to check that
inequalities (11.3d) are satisfied. Taking $j = \ell$ in (11.4)
immediately implies the lower estimate. To prove the upper
estimate, let us consider two cases:

\roster
\item"$\bullet$" If for some $j$  the minimum in
(11.4) is equal to the first of the quantities compared, then clearly
$${e^{\tilde t}} <  p_j^{s_j +
1}{\prod_{i = 1}^{j-1}
p_i^{s_i}} \le \prod_{j = 1}^\ell
p_j^{s_j+1}\,.
$$

       \item"$\bullet$" Otherwise, it follows that
$|\tilde\vy^{(p_j)}\cdot\tilde \vq|_{p_j} > p_j^{-(s_j+1)}$ for all $j$, and to
derive
the desired estimate it remains to notice that (10.4),  in view of (11.3e),
can be rewritten as $
c(\tilde\vy\cdot\tilde\vq)  = \prod_{j =
1}^{\ell}|\tilde\vy^{(p_j)}\cdot\tilde \vq|_{p_j}\le e^{-\tilde t}$.
       \qed
\endroster
\enddemo

\proclaim{11.2.\ Corollary} Assume that $\vy\in\bq_S^n$ is VWMA. Then
         for some $c,\gamma > 0$ there are infinitely many
$\tilde \vt\in\bz_{\sssize +}^{n+1}$ and
$\vs\in\bz_{\sssize +}^{\ell}$ such that
$$
e^{-(n+1)} \prod_{j = 1}^\ell p_j^{s_j} \le e^{\tilde t} < \prod_{j = 1}^\ell
p_j^{s_j+1}\tag 11.5a
$$
and
$$
\delta(g_{\vs,\tilde \vt}u_\vy {\Cal D}^{n+1}) \le ce^{-\gamma\tilde  t}\,.\tag
11.5b
$$
\endproclaim

\demo{Proof} By definition, for some $\vre > 0$ there are
infinitely many
solutions $\tilde \vq\in\bz^{n+1}$ of (10.4). Therefore, by the above lemma and
with
$\gamma$ as in (11.3b), there exists an unbounded set of
$\tilde \vt\in\br_{\sssize +}^{n+1}$  such that  (11.3c) holds for some
$\vs\in\bz_{\sssize
+}^{\ell}$ satisfying (11.3d). Denote by $[\tilde \vt]$ the vector consisting
of  integer
parts of
$t_i$, then clearly the ratio of $\delta(g_{\vs,\tilde \vt}u_\vy
{\Cal D}^{n+1})$
and $\delta(g_{\vs,[\tilde \vt]}u_\vy {\Cal D}^{n+1})$ is bounded from above by
some uniform constant.
        Thus, replacing $\tilde \vt$ by $[\tilde \vt]$, for some $c >
0$ one gets
infinitely many solutions
$\tilde \vt\in\bz_{\sssize +}^{n+1}$ of (11.5b), with $e^{\tilde t}$ being
smaller than before by at
most a factor of
$
e^{n+1}$, hence  (11.5a).
\qed\enddemo


\proclaim{11.3.\ Corollary}  Let $X$ be a 
\be\ metric space and
$\mu$ a  uniformly  Federer measure  on $X$. Suppose
we are given
a continuous map $\vf:X\to
\bq_S^n$  such that for $\mu$-a.e.\ $x_0\in X$ there exist a  ball
$B = B(x_0,r)$
and constants
$C,\alpha,\rho$ with the following property: for any
$\Delta\in \frak P({\Cal D}, n+1)$ and any $\vs\in\bz_{\sssize
+}^{\ell}$,
$\tilde \vt\in\bz_{\sssize
+}^{n+1}$ satisfying {\rm (11.5a)}, one has
$$\text{the function
$x\mapsto \cov\big(g_{\vs,\tilde \vt}u_{\vf(x)} \Delta\big)$ is \cag\ on
$B(x_0,3^{n+1}r)$  w.r.t.\
$\mu$}\,,\tag 11.6a
$$
and
$$
\sup_{x\in
B\,\cap\, \supp \mu}\,\cov\big(g_{\vs,\tilde
\vt}u_{\vf(x)}\Delta\big) \ge \rho\,.\tag
11.6b$$ Then $\vf_*\mu$ is strongly extremal.
\endproclaim

\demo{Proof} Applying Theorem 8.3,
with $h(x) = g_{\vs,\tilde \vt}u_{\vf(x)}$ and $m = n + 1$, we
conclude that
$$
\mu\left(\big\{x\in B\bigm| \delta(g_{\vs,\tilde \vt}u_{\vf(x)}
{\Cal D}^{n+1}) < ce^{-\gamma \tilde t}
\big\}\right)\le (n+1)C
\big(N_{\sssize X}D_{\mu}^2\big)^{n+1}
\left(\frac{ce^{-\gamma \tilde t}} \rho \right)^\alpha
\mu(B)
$$
whenever $ce^{-\gamma \tilde t} \le {\rho}$ and (11.5a) holds. Note that for
fixed
$\tilde \vt$,
the number of
different $\vs\in\bz_{\sssize +}^{\ell}$ satisfying
{\rm (11.5a)} is at most $\const\cdot  \tilde t\,^{\ell - 1}$. Therefore the
sum  (over all integer $\vs,\tilde \vt$ for which
{\rm (11.5a)} holds) of
measures of sets
$\big\{x\in B\bigm|
\delta(g_{\vs,\tilde \vt}u_{\vf(x)} {\Cal D}^{n+1}) < ce^{-\gamma \tilde t}
\big\}
$ is finite for every
$c,\gamma > 0$. An application of the Borel-Cantelli Lemma
shows that  for every
$c,\gamma > 0$ and $\mu$-a.e.\
$x\in B$, and hence for $\mu$-a.e.\
$x\in B$ and all
$c,\gamma > 0$,
there are at most finitely many integer solutions $\vs,\tilde \vt$ to
       (11.5ab).
Corollary 11.2 then implies that $\vf(x)$ is not VWMA for $\mu$-a.e.\
$x\in B$.
        \qed\enddemo

\subhead{11.4}\endsubhead We are now ready for the

\demo{Proof of Theorem 10.7, the case $\infty\notin S$}
Recall that we are given the balls $B_v \subset X_v$, $v\in S$, which will
be referred to as
$B_1,\dots,B_\ell$, and measures $\mu_v$ on $X_v$, which we will
call $\mu_1,\dots,\mu_\ell$. We will take
$B$ to be equal to
$\prod_{j = 1}^\ell B_j$ (recall that we are using the
product metric on $X$) and
show that it satisfies the assumptions of Corollary 11.3.
Thus we need to have explicit expressions for functions
$x\mapsto \cov\big(g_{\vs,\tilde \vt}u_{\vf(x)} \Delta\big)$.

Using Proposition 7.2 and Lemma 7.4, one can associate to
any nonzero submodule $\Delta \subset {\Cal D}^{n+1}$
of rank $r$ an element $\vw$  of $\bigwedge^r(\Cal D^{n+1})$
such that
$ \cov(\Delta) =  c (\vw)$ and $ \cov(g_{\vs,\tilde \vt}u_{\vf(x)}\Delta)
  =  c (g_{\vs,\tilde \vt}u_{\vf(x)}\vw)$.
It will be convenient to use  the standard
basis $\ve_0,\ve_1,\dots ,\ve_n$ of $\Cal R^{n+1}$, where $$\ve_i
= \left(\ve_i^{(v)}\right)_{v\in S^+}=
\left(\ve_i^{(p_1)},\dots,\ve_i^{(p_\ell)},\ve_i^{(\infty)}\right)$$
for each $i = 0,1,\dots,n$. Similarly,  we will use the standard
basis $\big\{\ve_{\sssize I} \mid I\subset \{0,1,\dots,n\}\big\}$
of $\bigwedge \Cal R^{n+1}$, where we let $\ve_{\sssize I} \df
\ve_{i_1}\wedge\dots\wedge \ve_{i_r}\in \bigwedge^r(\Cal R^{n+1})$
for $I = \{i_1,\dots,i_r\}\subset \{0,\dots,n\}$, $i_1 < i_2 <
\dots < i_r$. Thus we can write $\vw$ as above in the form $\vw =
\sum_{I\subset \{0,\dots,n\}}w_{\sssize I}\ve_{\sssize I}$, where
$w_{\sssize I}\in {\Cal D}$.
      \medskip

Now let us see how the coordinates of $\vw$ as above change under the
action of
$g_{\vs,\tilde \vt}u_{\vf(x)}$.
Note that:

\roster
\item"$\bullet$" $u_{\vf(x)}^{(\infty)}$ is trivial, and each
$\ve_{\sssize I}^{(\infty)}$ is an eigenvector of $(g_{\vs,\tilde
\vt})^{(\infty)}$ with eigenvalue
$e^{-t_{\sssize I}}$, where $t_{\sssize I} \df \sum_{i\in I}t_i$\,;

\item"$\bullet$"
the action of
$u_{\vf(x)}^{(p_j)}$ leaves $\ve_0^{(p_j)}$ invariant and sends
$\ve_{i}^{(p_j)}$ to
$\ve_{i}^{(p_j)} + f^{(p_j)}_i(x)
\ve_0^{(p_j)}$, $i = 1,\dots,n$, and  each
$\ve_{\sssize I}^{(p_j)}$ is an eigenvector of $(g_{\vs,\tilde \vt})^{(p_j)}$
with eigenvalue $1$ if $0\notin I$ and $p^{-{s_j}}$ otherwise; in other words,
$$
(g_{\vs,\tilde \vt}u_{\vf(x)}\ve_{\sssize I})^{(p_j)} =  \cases
&p_j^{-{s_j}}\ve^{(p_j)}_{\sssize I}
\text{ if }0\in I\\ &\ve_{\sssize I}^{(p_j)} + p_j^{-{s_j}}\sum_{i\in I} \pm
f^{(p_j)}_i(x)
\ve^{(p_j)}_{\sssize I \cup \{0\}\ssm\{i\}}\text{
otherwise}\,.\endcases\tag 11.7
$$
\endroster
Therefore one has \ $
(g_{\vs,\tilde \vt}u_{\vf(x)}\vw)^{(\infty)} = \sum_{ I}e^{-t_{\sssize
I}}w_{\sssize I}\ve_{\sssize
I}^{(\infty)}$ and
$$
(g_{\vs,\tilde \vt}u_{\vf(x)}\vw)^{(p_j)} = \sum_{0\notin I}w_{\sssize
I}\ve_{\sssize I}^{(p_j)} + p_j^{-{s_j}}\sum_{0\in I}\Big( w_{\sssize I} +
\sum_{i\notin I} \pm
w_{\sssize I\cup\{i\}\ssm\{0\}}f^{(p_j)}_i(x) \Big)
\ve_{\sssize I}^{(p_j)}\tag 11.8
$$
for $j =
1,\dots,\ell$.

        In particular,  real components of  all the
coordinates of $ g_{\vs,\tilde \vt}u_{\vf(x)}\vw$ are constant, and $p_j$-adic
components are linear combinations of
$1,f^{(p_j)}_1,\dots,f^{(p_j)}_n$. Condition (11.6a)
then immediately
follows
from Lemma 2.1(bc), (10.6) and Corollary 2.3. On the other hand,
for any $j = 1,\dots,\ell$ one can use (10.7) and the
compactness of the unit sphere in $\bq_{p_j}^{n+1}$ to find
$\rho_j > 0$ such that for any $\va = (a_0,a_1,\dots,a_n)\in\bq_{p_j}^{n+1}$
one has $$\sup_{x\in
B_j \,\cap\, \supp \mu_j}\,|a_0 +
a_1f^{(p_j)}_1(x)+\dots+a_nf^{(p_j)}_n(x)|_{p_j}
\ge
\rho_j\|\va\|_{p_j}\,.\tag 11.9$$
It remains to notice that all the components of $\vw$ necessarily appear
in the second sum in (11.8) (that is, the sum of terms with $0\in I$).
Therefore (11.8) and (11.9) imply
$$\sup_{x\in
B_j\,\cap\, \supp \mu_j}\,\|(
g_{\vs,\tilde \vt}u_{\vf(x)}\vw) ^{({p_j})}\|_{p_j}\ge \rho_j
{p_j}^{{s_j}}\max_{\sssize
I}|w_{\sssize I}|_{p_j}\,,\tag 11.10
$$
and hence
$$
\split
\sup_{x\in
B\,\cap\, \supp \mu}\,c(
g_{\vs,\tilde \vt}u_{\vf(x)}\vw) &\ge \left(\prod_{j = 1}^\ell \rho_j
{p_j}^{{s_j}}\max_{\sssize I}|w_{\sssize I}|_{p_j} \right) \max_{\sssize
I}{e}^{-t_{\sssize I}}|w_{\sssize I}|_{\infty}\\  &\ge \left(\prod_{j =
1}^\ell \rho_j
        \right) {e}^{-\tilde t}\left(\prod_{j = 1}^\ell
{p_j}^{{s_j}}\right) \max_{\sssize
I}c(w_{\sssize I}) \un{(11.5a)}\ge  \prod_{j = 1}^\ell \frac{\rho_j}{p_j}
        \,.
\endsplit
$$
Condition (11.6b) is thus established, and the theorem
follows. \qed
\enddemo

\heading{12.
Proof of Theorem 10.7 for
$\infty\in S$}
\endheading

\subhead{12.1}\endsubhead In this section we  suppose that  $S =
\{p_1,\dots,p_{\ell-1},\infty\}$, where
$p_1,\dots,p_{\ell-1}$ are distinct primes. In this case there is no
need
to artificially add the infinite valuation to $S$; that is,
we will now work with $${\Cal R} \df \bq_{S} \quad \text{and} \quad {\Cal D}
\df {\Cal O}_{S} =
\bz[\tfrac1{p_1},\dots,\tfrac1{p_{\ell-1}}]\,.$$ The element $u_\vy\in
\GL^1(n+1,{\Cal R})$ can be simply defined by
$
u_\vy \df \pmatrix 1 & \vy\\ 0 & I_n\endpmatrix
$,
       so that
$$u_\vy {\Cal D}^{n+1} = \left\{\left.
\pmatrix\tilde\vy\cdot\tilde\vq \\ \vq\endpmatrix\right|\tilde \vq\in {\Cal
D}^{n+1}\right\}\,.
$$
The definition (11.2) of the diagonal element $g_{\vs,\tilde \vt}\in
\GL(n+1,{\Cal R})$ given in the previous section
will still be valid, except for 
       $\vs$
       now having $\ell -1$  components.

Now let us
split the set of VWMA vectors into two parts:
say that a \vwma\
$\vy\in \bq_S^n$ is
{\sl VWMA$_\le$\/}  if for some positive $\vre$ there are
infinitely many solutions
$\tilde \vq$ to (10.4$_\infty$) satisfying
$$
|q_0| \le
\left(1 +
n\|\vy^{(\infty)}\|_\infty\right)\|\vq\|\,,\tag 12.1
$$
and that it is {\sl VWMA$_>$\/}
otherwise.
Our  strategy will be as follows: we will modify
the dynamical approach of the previous section to treat the first case,
and use the conclusion of the ``$\infty\notin S$'' case of Theorem 10.7 to
take care of the second case.

\medskip

Here is a replacement for Lemma 11.1.

\proclaim{Lemma} Let $\vre > 0$, $\vy\in \bq_S^{n}$ 
and
$\tilde \vq\in\bz^{n+1}$ be such that {\rm (10.4$_\infty$)} and {\rm
(12.1)} hold.
For $i = 1,\dots,n$ define  $t_i> 0$ by 
$$
|q_i|_{\sssize +} = \Pi_{\sssize +}(\vq)^{-\frac{\vre}{n+1}}
e^{t_i}\,,\tag 12.2a
$$
and let 
$$
\gamma = \frac{\vre}{n +1 + n\vre}\,.\tag 12.2b
$$
Then there exist
$\vs = (s_1,\dots,s_{\ell-1})\in\bz_+^{\ell-1}$ and $t_0\in\br$
such that
$$
\delta(g_{\vs,\tilde \vt}u_\vy {\Cal D}^{n+1}) \le
\sqrt{n+1}e^{-\gamma t}\,,\tag 12.2c
$$
$$
-t \le t_0 \le  t + \ln \left(1 +
2n\|\vy^{(\infty)}\|_\infty\right)\,,\tag 12.2d
$$
and 
$$
\prod_{j = 1}^{\ell-1} p_j^{s_j} \le e^{\tilde t} < \prod_{j = 1}^{\ell-1}
p_j^{s_j+1}\,.\tag
12.2e
$$
\endproclaim

\demo{Proof}
As in the proof of Lemma 11.1, we consider the product of
equalities (12.2a), namely
$$
e^{t} = \Pi_{\sssize +}(\vq)^{1 + \frac n {n + 1}\vre}\,,\tag 12.2f
$$
and then write
$$
e^{-t_i}|q_i|\le e^{-t_i}|q_i|_{\sssize +}\un{(12.2a)}\le \Pi_{\sssize
+}(\vq)^{-\frac{\vre}{n+1}}\un{(12.2bf)} =  e^{-\gamma t},\ i =
1,\dots,n\,.\tag 12.2g
$$
After that define
$t_0$ by
$$e^{-t_0} \df \min\left({e^{t}},
\frac{e^{-\gamma t}}{|\tilde\vy^{(\infty)}\cdot\tilde
\vq|_{\infty}}\right)\,.\tag 12.3
$$
It follows that $|e^{-t_0}\tilde\vy^{(\infty)}\cdot\tilde\vq|_{\infty} \le
e^{-\gamma t}$, hence $
\|(g_{\vs,\tilde \vt})^{(\infty)}u_\vy^{(\infty)}\tilde\vq\|_\infty \le
\sqrt{n+1}e^{-\gamma t}$. The lower estimate in (12.2d) is immediate
from (12.3), while the upper estimate clearly holds if the minimum in
(12.3) is equal to
$e^{t}$, and otherwise one has
$$
e^{-t_0}  =
\frac{e^{-\gamma t}}{|q_0 + \vy^{(\infty)}\cdot\vq|_{\infty}}\
\un{(12.1)}\ge \
\frac{e^{-\gamma t}}{\left(1 +
2n\|\vy^{(\infty)}\|_\infty\right)\|\vq\|}\ \un{(12.2g)}\ge \
\frac{e^{- t}}{1 +
2n\|\vy^{(\infty)}\|_\infty}\,.
$$
Now that all the components of $\tilde\vt$ are chosen, we can define
$s_j$, $j = 1,\dots,\ell-1$, as in (11.4). After that one can  verify,
following the lines of the proof of Lemma 11.1, that   $$
\|(g_{\vs,\tilde \vt})^{(p_j)}u_\vy^{(p_j)}\tilde\vq\|_{p_j} =
\max\big(\,p_j^{s_j}|\tilde\vy^{(p_j)}\cdot\tilde\vq|_{p_j}, 1\big) \le
1$$
for each $j$, so that $c(g_{\vs,\tilde \vt}u_\vy
\tilde\vq) \le \sqrt{n+1}e^{-\gamma t}
$, and  that
inequalities (12.2e) are satisfied.
        \qed
\enddemo

\proclaim{12.2.\ Corollary} Assume that  $\vy\in\bq_S^n$ is VWMA$_\le$.
Then
         for some $c_0,c,\gamma > 0$ there are infinitely many
$\tilde\vt\in\bz\times\bz_{\sssize +}^{n}$  and
$\vs\in\bz_{\sssize +}^{\ell-1}$ satisfying $$
-t - 1 \le t_0 \le  t + c_0\,,\tag 12.4a
$$
$$
e^{-(n+1)} \prod_{j = 1}^{\ell-1} p_j^{s_j} \le e^{\tilde t} < \prod_{j =
1}^{\ell-1} p_j^{s_j+1}\,,\tag 12.4b
$$
and
$$
\delta(g_{\vs,\tilde \vt}u_\vy {\Cal D}^{n+1}) \le ce^{-\gamma t}\,.\tag 12.4c
$$
\endproclaim

\demo{Proof}
By definition,
      for some $\vre > 0$ there are
infinitely many solutions
$\tilde \vq$ to {\rm (10.4$_\infty$)} and {\rm (12.1)}.
         Therefore, by the
above lemma and with
$\gamma$ as in (12.2b), there exists
an unbounded set of
$\vt\in\br_{\sssize +}^{n}$  such that  inequalities (12.2cde) hold for
some
$t_0\in\br$ and $\vs\in\bz_{\sssize
+}^{\ell-1}$.
The rest of the proof of Corollary 11.2 applies verbatim.
\qed\enddemo


\proclaim{12.3.\ Corollary}  Let $X$ be a 
\be\ metric space and
$\mu$ a  uniformly  Federer measure  on $X$. Suppose
we are given
a continuous map $\vf:X\to
\bq_S^n$  with the follolwing property:   for $\mu$-a.e.\
$x_0\in X$ there exist a  ball
$B = B(x_0,r)$
and constants
$C,\alpha,\rho, c_0$ such that conditions {\rm (11.6ab)} hold for any
$\Delta\in \frak P({\Cal D}, n+1)$ and any $\vs\in\bz_{\sssize +}^{\ell-1}$,
$\vt\in\bz\times\bz_{\sssize +}^{n}$ satisfying {\rm (12.4ab)}.
       Then $\vf(x)$ is not VWMA$_\le$ for $\mu$-a.e.\ $x\in X$.
\endproclaim

\demo{Proof} An application of Theorem 8.3,
again with $h(x) = g_{\vs,\tilde \vt}u_{\vf(x)}$ and $m = n + 1$, yields
$$
\mu\left(\big\{x\in B\bigm| \delta(g_{\vs,\tilde \vt}u_{\vf(x)}
{\Cal D}^{n+1}) < ce^{-\gamma t}
\big\}\right)\le (n+1)C
\big(N_{\sssize X}D_{\mu}^2\big)^{n+1}
\left(\frac{ce^{-\gamma t}} \rho \right)^\alpha
\mu(B)
$$
whenever $ce^{-\gamma t} \le {\rho}$ and (12.4ab) hold. Now observe that
for fixed
$\vt$,
the number of
different $t_0\in\bz$ and $\vs\in\bz_{\sssize +}^{\ell-1}$ satisfying
{\rm (12.4ab)} is  at most $\const\cdot t^{\ell - 1}$. Therefore
for any
$c,\gamma > 0$ the sum (over all integers $\vs,\tilde \vt$ for which
inequalities (12.4ab) hold) of measures of sets
$\big\{x\in B\bigm|
\delta(g_{\vs,\tilde \vt}u_{\vf(x)} {\Cal D}^{n+1}) < ce^{-\gamma t}
\big\}$ converges. As before, an application of the Borel-Cantelli Lemma
shows that
        for $\mu$-a.e.\
$x\in B$ there are at most finitely many integer solutions
$\vs,\tilde \vt$ to (12.4abc) for any
$c,\gamma > 0$.
Corollary 12.2 then implies that $\vf(x)$ is not VWMA$_\le$ for
$\mu$-a.e.\
$x\in B$.
        \qed\enddemo

\subhead{12.4}\endsubhead  Now let us state a lemma showing that $\vy$ being
VWMA$_>$ has some implications to the \di\ properties of its ``finite part''
$
\left(\vy^{(v)}\right)_{v\in S_f}$.

\proclaim{Lemma}  Assume that $\vy = \left(\vy^{(v)}\right)_{v\in
S}$ is VWMA$_>$; then $\left(\vy^{(v)}\right)_{v\in S_f}$
      is VWMA.
\endproclaim

\demo{Proof} By assumption, there exist infinitely many solutions
$\tilde\vq$ of {\rm (10.4$_\infty$)} for which {\rm (12.1)} fails. For each
of them one can write
$$
\split
|\tilde\vy^{(\infty)}\cdot\tilde\vq|_\infty  &=  |q_0
+ \vy^{(\infty)}\cdot\vq|_\infty  \ge |q_0| -
n\|\vy^{(\infty)}\|_\infty\|\vq\|\\
      &\ge |q_0| -\frac{ n\|\vy^{(\infty)}\|_\infty}{
1 + n\|\vy^{(\infty)}\|_\infty} |q_0|  = \frac{|q_0|_+}{
1 + n\|\vy^{(\infty)}\|_\infty}  \,.
\endsplit
$$
Therefore one has
$$
\split
\prod_{v\in S_f}|\tilde\vy^{(v)}\cdot\tilde\vq|_v
&=
\frac{c(\tilde\vy\cdot\tilde\vq)}{|\tilde\vy^{(\infty)}\cdot\tilde\vq|_\infty
}   \le
\Pi_+(\vq)^{-(1 +
\vre)}|q_0|_+^{-\vre}\cdot \left({
1 + n\|\vy^{(\infty)}\|_\infty}\right)|q_0|_+^{-1}\\
&\le \left({
1 + n\|\vy^{(\infty)}\|_\infty}\right) \Pi_+(\tilde \vq)^{-(1 +
\vre)}\,,\quad
\endsplit
$$
which finishes the proof modulo a slight change of $\vre$. \qed
\enddemo

\subhead{12.5}\endsubhead  Finally we are  ready for the

\demo{Proof of Theorem 10.7, the case $\infty \in S$} Applying the case
``$\infty\notin S$'' of Theorem 10.7 to the map $\left(\vf^{(v)}\right)_{v\in
S_f}$, we obtain that the pushforward of $\prod_{v\in
S_f}\mu_v$ by $\left(\vf^{(v)}\right)_{v\in S_f}$  must be
strongly extremal, which, in view  of Lemma 12.4, implies that $\vf(x)$ is not
VWMA$_>$ for $\mu$-a.e.\ $x\in X$. Thus, as before,  it suffices to
take $$B = \prod_{v\in S}B_v = \prod_{j=1}^{\ell-1}B_j\times
B_\infty\quad\text{and}\quad \mu = \prod_{v\in S}\mu_v =
\prod_{j=1}^{\ell-1}\mu_j\times \mu_\infty\,,$$ and  check that the assumptions
of Corollary 12.3
are satisfied  by means of writing down explicit expressions for functions
$x\mapsto \cov\big(g_{\vs,\tilde \vt}u_{\vf(x)} \Delta\big)$. This again boils
down to the computation of components of $g_{\vs,\tilde \vt}u_{\vf(x)}\vw$,
where $\vw = \sum_{I\subset \{0,\dots,n\}}w_{\sssize
I}\ve_{\sssize I}\in \bigwedge^r({\Cal D}^{n+1})$.

Since there was no change in the ultrametric components of
$g_{\vs,\tilde \vt}$ and $u_{\vf(x)}$, formula (11.7) is still valid.
Furthermore, an expression for the Archimedean components turns out to
be similar to (11.7):
$$
(g_{\vs,\tilde \vt}u_{\vf(x)}\ve_{\sssize I})^{(\infty)} =  \cases
&e^{-{t_I}}\ve^{(\infty)}_{\sssize I}
\text{ if }0\in I\\ &e^{-{t_I}}\ve^{(\infty)}_{\sssize I} +
\sum_{i\in I} \pm e^{-t_{I\cup \{0\}\ssm\{i\}}}f^{(\infty)}_i(x)
\ve^{(\infty)}_{\sssize I \cup \{0\}\ssm\{i\}}\text{
otherwise}\,.\endcases
$$
Therefore one has (11.8) for $j =
1,\dots,\ell-1$, and, in addition,
$$
(g_{\vs,\tilde \vt}u_{\vf(x)}\vw)^{(\infty)} = \sum_{0\notin I}w_{\sssize
I}e^{-{t_I}}\ve_{\sssize I}^{(\infty)} + \sum_{0\in I}e^{-{t_I}}\Big(
w_{\sssize I} +
\sum_{i\notin I} \pm
w_{\sssize I\cup\{i\}\ssm\{0\}}f^{(\infty)}_i(x) \Big)
\ve_{\sssize I}^{(\infty)}\,.
$$

We see that  real (resp.\ $p_j$-adic) components of  all the
coordinates of $ g_{\vs,\tilde \vt}u_{\vf(x)}\vw$ are linear combinations of
$1,f^{(\infty)}_1,\dots,f^{(\infty)}_n$ (resp.\
$1,f^{(p_j)}_1,\dots,f^{(p_j)}_n$). Condition (11.6a)
then immediately follows
from Lemma 2.1(bc), (10.6) and Corollary 2.3. On the other hand,
an argument identical to that of the previous section shows
that for every $j = 1,\dots,\ell-1$ there exists $\rho_j > 0$
such that (11.10) holds, and also that there exists $\rho_\infty > 0$
such that
$$\sup_{x\in
B_\infty\,\cap\, \supp \mu_\infty}\,\|(
g_{\vs,\tilde \vt}u_{\vf(x)}\vw) ^{({\infty})}\|_{\infty}\ge
\rho_\infty \min_{0\in
I}e^{-{t_I}} \max_{\sssize
I}|w_{\sssize I}|_{\infty}\ge \rho_\infty e^{-{\tilde t}} \max_{\sssize
I}|w_{\sssize I}|_{\infty}\,.
$$
Therefore
$$
\split
\sup_{x\in
B\,\cap\, \supp \mu}\,c(
g_{\vs,\tilde \vt}u_{\vf(x)}\vw) &\ge \left(\prod_{j = 1}^{\ell-1} \rho_j
{p_j}^{{s_j}}\max_{\sssize I}|w_{\sssize I}|_{p_j} \right)
\rho_\infty e^{-{\tilde t}}\max_{\sssize I}|w_{\sssize I}|_{\infty}\\  &\ge
\rho_\infty \left(\prod_{j = 1}^{\ell-1} \rho_j
        \right) {e}^{-\tilde t}\left(\prod_{j = 1}^{\ell-1}
{p_j}^{{s_j}}\right) \max_{\sssize
I}c(w_{\sssize I}) \un{(12.4b)}\ge \rho_\infty  \prod_{j = 1}^{\ell-1}
\frac{\rho_j}{p_j}
        \,.
\endsplit
$$
This implies (11.6b) and shows that $\vf(x)$ is not VWMA$_\le$ for
$\mu$-a.e.\
$x\in X$, thus finishing the proof of the theorem. \qed
\enddemo


\heading{13. More on $S$-arithmetic \da}\endheading

\subheading{13.1.  Extensions of $\bq$}
 It seems to be a natural task to 
 extend the  metric Diophantine approximation results
proved in this paper to the 
framework of an arbitrary number field
$K$. Indeed, 
the main quantitative nondivergence estimate of the paper
(Theorem 8.3) can be rather straightforwardly generalized to the setting 
of maps from \be\ metric spaces into $\GL(m,K_S)$, where $K$ is
a finite extension of $\bq$, $S$ is a finite set of its normalized
valuations containing all the Archimedean ones,  and $K_S$ 
is the direct product of
completions $K_v$ of $K$ over $v\in S$. See the earlier version 
\cite{KT} of the present paper for more detail.
Similarly  one can mimic the presentation of \S\S 10--12 to define 
{\sl \vwa\/} elements of $K_S^n$, and  prove that those form a null 
set with respect to  pushforwards of Haar measure 
by products of nondegenerate maps.

However, understanding multiplicative approximation over an arbitrary
number field turns out to be more
complicated\footnote{except  when 
$K$ is an  imaginary quadratic extension of $\bq$ (cf.\ \cite{DK}),
in which case the proof of the analogs of our results
 can
  be carried out without major changes}.
Indeed, if $K$ has more than one infinite valuation, the group of units of 
the ring of
integers of $K$ is infinite, which complicates the definition of
VWMA vectors  and makes   proofs of the corresponding results
more delicate. 
The case of an arbitrary number field will be treated in a forthcoming
paper.

\subheading{ 13.2.  Khintchine-type theorems}  Another way
to generalize the \di\ set-up of this paper would be to 
 replace the right hand side of (0.4)  by an
arbitrary function of
$\|\vq\|$ or $\|\tilde\vq\|$.
With the notation of \S 0.4, let us introduce the following definition: 
for a non-increasing function
$\psi:\bn\to (0,\infty)$, say that $\vy\in \bq_S^n$ is 
{\sl $\psi$-approximable\/} if there are infinitely many solutions $\tilde\vq
= (q_0,q_1,\dots,q_n) \in\bz^{n+1}$ to
$$|q_0 + \vq\cdot\vy|^\ell \le \cases \psi(\|\tilde\vq\|)
\quad\text{ if }\infty\notin S\\
\psi(\|\vq\|)\quad\text{ if }\infty\in S\,.\endcases\
%
$$
As in the case of (0.4), it is easy to check using the 
Borel-Cantelli lemma that $\lambda$-a.e.\
 $\vy\in \bq_S^n$  is not $\psi$-approximable whenever the series
$$\cases \sum_{k = 1}^\infty k^n\psi(k)
\quad\text{ \ \ \  if }\infty\notin S\\
\sum_{k = 1}^\infty k^{n-1}\psi(k)\quad\text{ if }\infty\in S\endcases\
\tag
13.1
$$
converges, and using the methods of \cite{L} it should be possible to
prove that $\lambda$-a.e.\
 $\vy\in \bq_S^n$  is  $\psi$-approximable if the above series diverges.
Similar questions then arise regarding measures  other than $\lambda$,
in particular,  volume measures on nondegenerate smooth
manifolds or, in the case $\ell > 1$, their products.

In recent years the case  $S = \{\infty\}$ has been 
completely understood,
see \cite{BKM, Be} for the convergence case and 
\cite{BBKM} for the divergence case. That is, the convergence 
(divergence)
of (13.1) was shown to imply that almost no (almost all)
points on nondegenerate submanifolds of $\br^n$ are $\psi$-approximable.
Combining the approach 
of the present paper with the methods of \cite{BKM}  and 
\cite{BBKM} respectively, it seems plausible that both convergence 
and 
 divergence cases can be proved for $\vf_*\lambda$ as in Theorem 0.4. 
Note that when $S = \{p\}$, both cases were recently established for the curve
  (0.6) \cite{BBK} and  for $\lambda$-a.e.\ nondegenerate 
$\vf:\bz_p\to\bz_p^2$ which is normal in the
sense of Mahler \cite{BK}.  The convergence
case for nondegenerate curves in $\bz_p^3$ was treated by E.\ Kovalevskaya in
 \cite{Ko1--2}, and in another paper  \cite{Ko3} she extended the method of
 \cite{BK} to obtain a  result involving both
$p$-adic and infinite valuations. Note also that
 the paper \cite{BKM} contains a more general (in
particular, multiplicative) version of the convergence case for
nondegenerate submanifolds of $\br^n$, 
and it would be 
interesting to see whether  the $S$-arithmetic
 set-up can be treated in a similar way.

\subheading{ 13.3.  Analogues of other results over $\br$}
Since the introduction \cite{KM} of the dynamical approach 
to \da\ on manifolds, various extensions and 
generalizations of the method have been found.
 We expect that many of the 
ideas developed recently for \da\ over $\br$ 
 can be applied in the non-Archimedean
 setting. Specifically we would like to propose two conjectures,
in which  $\bq_v$ stands for an arbitrary completion of $\bq$.

\proclaim{Conjecture IS (Inheritance for Subspaces)} Let $\Cal L$ be an 
affine subspace of $\bq_v^n$ and 
let $\vf:\bq_v^d\to \Cal L$ be a $C^k$ map which is 
{\sl nondegenerate in}\footnote{that is, the linear part of $\Cal L$ is spanned
by partial derivatives of $\vf$} $\Cal L$ at $\lambda$-a.e.\ point. 
Suppose that the volume measure on $\Cal L$ is extremal 
(resp.\ strongly extremal); then so is $\vf_*\lambda$.
\endproclaim

This was proved in \cite{K2} for $v = \infty$, and in addition
explicit necessary and sufficient conditions, involving
coefficients of linear functions parametrizing  $\Cal L$, were
found  for the volume measure on $\Cal L$ to be extremal 
(strongly extremal). There should be no major difficulties in 
extending these results to the non-Archimedean case.

\proclaim{Conjecture FP (Friendliness of Pushforwards)} 
Let $\mu$  be a self-similar measure on the limit set 
 of an irreducible family of
contracting similitudes of $\bq_v^d$ satisfying the
open set condition (see \S 10.6), and let $\vf:\bq_v^d
\to\bq_v^n $ be a 
smooth enough map which is nondegenerate at $\mu$-a.e.\ point. Then 
$\vf_*\mu$ is strongly extremal.
\endproclaim

The case  $v = \infty$ of the above conjecture
is one of the main results of  \cite{KLW}. Note that a key step
of the proof, see \cite{KLW, Proposition 7.3}, crucially involves the
Mean Value Theorem, and for its non-Archimedean analogue one
would need to come up with a  replacement, perhaps similarly to
our approach to Proposition 3.3.




\enddocument